\documentclass[12pt,a4paper]{article}

\usepackage{GoetzArticle}
\usepackage{amssymb}
\usepackage{verbatim}
\usepackage{amssymb}
\usepackage{amsmath}
\usepackage{graphicx}
\usepackage{epsfig}
\usepackage{color}

\usepackage{rotating}

\usepackage{rawfonts}
\usepackage{bm}
\input prepictex
\input pictexwd
\input postpictex

\newcommand{\mysmallmatrix}[2]{
            \left(\footnotesize
                \setlength{\arraycolsep}{4pt}
                \renewcommand{\arraystretch}{.7}
                    \begin{array}{*{#1}{r}} #2
                    \end{array}
            \renewcommand{\arraystretch}{1}
    \right)
    }

\newcommand{\mysmallmatrixLEFT}[2]{
            \left(\footnotesize
                \setlength{\arraycolsep}{4pt}
                \renewcommand{\arraystretch}{.7}
                    \begin{array}{*{#1}{l}} #2
                    \end{array}
            \renewcommand{\arraystretch}{1}
    \right)
    }

\newcommand{\mysmallmatrixCENTER}[2]{
            \left(\footnotesize
                \setlength{\arraycolsep}{1pt}
                \renewcommand{\arraystretch}{.7}
                    \begin{array}{*{#1}{c}} #2
                    \end{array}
                \renewcommand{\arraystretch}{1}
            \right)
    }

\newcommand{\fhat}[1]{\widehat{#1}}
\newcommand{\Z}{\mathbb{Z}}

\newcommand{\supp}{{\rm supp}\, }

\newcommand{\R}{\mathbb{R}}

\newcommand{\N}{\mathbb{N}}

\newcommand{\C}{\mathbb{C}}

\newcommand{\Dcal}{{\mathcal D}}

\newcommand{\rank}{{\rm rank}\, }

\newcommand{\adj}{\text{adj }}

\newcommand{\vindex}{\vspace{-.7cm}}
\newcommand{\vpropend}{\vspace{-.6cm}}

\title{Uncertainty in time--frequency representations on finite Abelian groups and applications}
\author{Felix Krahmer \footnotemark[2] \and G\"otz E. Pfander \footnotemark[3] \and Peter
Rashkov \footnotemark[3]}

\begin{document}
\maketitle
\date{}

\begin{abstract}

Classical and recent results on  uncertainty principles for
functions on finite Abelian groups relate the cardinality of the
support of a function to the cardinality of the support of its
Fourier transforms. We use these results and their  proofs to obtain
similar results relating the support sizes of functions and their
{\it short--time Fourier transforms}. Further, we discuss
applications of our results. For example, we use our results to
construct a class of equal norm tight Gabor frames that are
maximally robust to erasures and we discuss consequences of our
findings to the theory of recovering and storing signals which have
sparse time--frequency representations.
\end{abstract}
\renewcommand{\thefootnote}{\fnsymbol{footnote}}

\footnotetext[3]{School of Engineering and Science,
    International University Bremen, 28759 Bremen, Germany.
    }

\footnotetext[2]{
    Courant Institute of Mathematical Science, New York University, New York USA
    }
\renewcommand{\thefootnote}{\arabic{footnote}}
\pagenumbering{arabic} \maketitle



\section{Introduction}
The uncertainty principle establishes restrictions on how well
localized the Fourier transform of a well localized function can be
and vice versa.  In the case of a function defined on finite Abelian
groups, localization can be expressed through the cardinality of the
support of the function. This case has recently drawn renewed
interest. This is due in part to their relevance for compressed
sensing and, in particular, for the recovery of lossy signals under
the assumption of restricted spectral content \cite{CRT04}.

A classical result on the uncertainty principle for functions
defined on finite Abelian groups states that the product of the
number of nonzero entries in a nontrivial vector, i.e., nontrivial
function on a finite set, and the number of nonzero entries in its
Fourier transform is not smaller than the order of the group
\cite{DS89}. This result can be improved for cyclic groups of prime
order: the sum of the number of nonzero entries in a vector and the
number of nonzero entries in its Fourier transform exceeds the order
of the group \cite{Tao05}. Further, it has recently been shown that
the classical bound can be refined  for almost any finite Abelian
group \cite{Mes05}.

The objective of this paper is to establish results similar to those
discussed above for joint time--frequency representations, that is,
to obtain restrictions on the cardinality of the support of joint
time--frequency representations of functions defined on finite
Abelian groups. For example, let us consider the simplest
time--frequency representation of a function, namely the one that is
given by the tensor product of a function and its Fourier transform.
In this case, the classical result on the uncertainty principle for
nontrivial functions on finite Abelian groups states that the
cardinality of the support of this tensor is at least the order of
the group.

In the following though, we shall be mostly interested in
time--frequency representations given by  short--time Fourier
transforms. It is well-known that, again, the cardinality of the
support of any short--time Fourier transform of a nontrivial
function defined on a finite Abelian group is bounded below by the
order of the group. As seen below, we can improve this bound by
using the subgroup structure of the groups and/or by allowing only
well-chosen window functions. For example, we show that for a group
with prime order and for almost every window function, the sum of
the cardinality of the support of the analyzed function and the
cardinality of its short--time Fourier transform exceeds the square
of the order of the group (see Theorem~\ref{theorem:p2+1}).

In addition to the above, we shall give applications of our results
to the theory of so-called Gabor frames and the theory of sparse
signal recovery. For example, the results on the cardinality of the
support of short--time Fourier transforms can be translated into
criteria for the recovery of encoded signals from a channel with
erasures.

The paper is organized as follows. In
Section~\ref{section:background} we give a brief but self-contained
account of the Fourier transformation and of the short--time Fourier
transformation for functions defined on finite Abelian groups.
Section~\ref{section:FourierUncertainty} discusses uncertainty
principles which relate the cardinality of the support of functions
with the cardinality of the support of their Fourier transforms.  We
start Section~\ref{section:FourierUncertainty} with a classical
result which is based on standard norm estimates \cite{DS89}. In
Section~\ref{section:taobiroFT} we state results based on the minors
of Fourier transform matrices and which apply only to functions
defined on cyclic groups of prime order \cite{Tao05}. Finite Abelian
groups of any order are analyzed in
Section~\ref{section:meshulamFT}. There, the underlying subgroup
structure of finite Abelian groups is used to obtain improvements to
the classical uncertainty result discussed above~\cite{Mes05}.

Section~\ref{section:STFTuncertainty} is devoted to uncertainty in
the short--time Fourier transformation. Following the organization
of Section~\ref{section:FourierUncertainty}, a discussion of general
results is followed by results for functions defined on cyclic
groups in Section~\ref{section:taobiroSTFT}. Other finite Abelian
groups are covered in Section~\ref{section:meshulamSTFT}. We
conclude our discussion of the cardinality of the support set of
short--time Fourier transforms in Section~\ref{section:outlook} with
a conjecture on the possible cardinalities of the support of
short--time Fourier transforms with respect to a random window
function. In fact, one of the major difficulties to obtain
uncertainty principles for the short--time Fourier transform is its
dependence on the chosen window function.

Section~\ref{section:applications} is devoted to applications of our
findings. In Section~\ref{subsection:applications-Gabor} we give
applications of the results of Section~\ref{section:STFTuncertainty}
to communications engineering. There, we discuss the
identification/measurement problem for time--varying
operators/channels and the transmission through channels with
erasures. In addition, we show the existence of a large class of
equal norm tight frames of Gabor type. In
Section~\ref{subsection:applications-sparsity} we briefly discuss
connections of our work to the recovery of signals which have a
sparse representation in a given dictionary.


\section{Background and Notation} \label{section:background}

For any finite set $A$ we set $\displaystyle
\C^A=\{f:A\longrightarrow \C\}$. For $|A|=|B|=n$, $\C^A\cong
\C^B\cong\C^n$ as vector spaces, where $|A|$ denotes the cardinality
of the set $A$. Further, for $A\subseteq B$, we write
$A^c=B{\setminus} A$ and we define the embedding operator
$i_A:\C^A\longrightarrow \C^{B}$ where $i_Af(x)=f(x)$ for $x\in A$
and $i_Af(x)=0$ for $x\in A^c$. Correspondingly, we define the
restriction operator $r_{A}:\C^B\longrightarrow \C^A$. Similarly,
every map $S:A\longrightarrow B$ induces a map $\widetilde
S:\C^B\longrightarrow \C^A$, $\left(\widetilde S f\right)(a)=
f\left(S(a)\right)$. If $S$ is bijective, then $\widetilde S$ is bijective as well.%

For $M \in \C^{m\times n}$ and $A\subseteq\{0,1,\ldots,n{-}1\}$ and
$B\subseteq\{0,1,\ldots,m{-}1\}$ we let $M_{A,B}$ denote the
$|B|{\times} |A|$--submatrix of $M$ which represents $r_B \circ M
\circ i_A$.

For $f\in \C^A$, we use the now customary notation $\|f\|_0=|\supp
f|$ where $\supp f=\{a\in A:\ f(a)\neq 0\}$. Clearly, $\|\cdot\|_0$
is not a norm.

\subsection{Fourier transforms on finite Abelian groups}

Throughout this paper, $G$ denotes a finite Abelian group. The
identity element of $G$ is denoted by $e$ or by $0$ in case that $G$
is cyclic, i.e., if  $G=\Z_n$ for some $n\in\N$.  The dual group of
characters $\widehat G$ of  $G$ is the set of continuous
homomorphisms $\xi\in \C^G$ which map $G$ into the multiplicative
group $S^1=\{z\in \C: \ |z|=1\}$. The set $\widehat G$ is an Abelian
group under pointwise multiplication and, as is customary, we shall
write this commutative group operation additively. Note that $G$ is
isomorphic to $\widehat G$. Further, Pontryagin duality implies that
$\widehat{\widehat G}$ can be canonically identified with $G$, a
fact which is emphasized by writing $\langle \xi,x\rangle=\xi(x)$.

The Fourier transform ${\mathcal F} f=\widehat{f}\in\C^{\widehat G}$
of $f\in \C^G$ is given by
$$
    \widehat{f}(\xi)=\sum_{x\in G} f(x)\, \overline{\xi(x)}=
    \sum_{x\in G} f(x)\, \overline{\langle\xi,x\rangle},\quad
    \xi \in  \widehat{G}\,.
$$
The inversion formula for the Fourier transformation allows us to
reconstruct the original function from its Fourier transform.
Namely, for $f\in \C^G$ we have
$$
    f(x)=\tfrac 1 {|G|} \sum_{\xi\in \widehat G} \widehat f(\xi)\, \langle \xi,x\rangle, \quad
    x \in G\,.
$$
The inversion formula implies that
\begin{eqnarray}
    \|f\|_2^2=\tfrac 1 {|G|} \sum_{\xi\in \widehat G} |\widehat f(\xi)|^2=\tfrac 1
        {|G|}\|\widehat f\|_2^2, \label{equation:FourierInversion}
\end{eqnarray}
where $\|f\|_2:= (\sum_{t \in G}|f(t)|^2)^{\frac{1}{2}}$. Further,
(\ref{equation:FourierInversion}) together with
$\|\xi\|_2=|G|^{\frac{1}{2}}$ for all $\xi\in\widehat G$ implies
that the normalized characters in
$\{|G|^{-\frac{1}{2}}\xi\}_{\xi\in\widehat G}$ form an orthonormal
basis for $\C^G$, and $\sum_x \langle \xi,x\rangle=0$ if $\xi\neq 0$
and $\sum_\xi \langle \xi,x\rangle=0$ if $x\neq 0$.


Fourier transformations are linear maps and we turn now to a
discussion of their matrix representations.

For $n\in\N$ and $\omega=e^{2\pi i /n}$, the {\em discrete Fourier
matrix} $W_{\Z_n}$ of the cyclic group $\Z_n$ is defined by
$W_{\Z_n}=(\omega^{rs})_{r,s=0}^{n{-}1}$. Identifying $\C^{\Z_n}$
with $\C^n$, we have $\fhat f = W_{\Z_n}\cdot f$.

For an arbitrary finite Abelian group $G$, we can always choose a
representation of $G$ as  direct product of cyclic groups $G \cong
\Z_{d_1} \times \Z_{d_2} \times \ldots \times \Z_{d_m}$ where $d_1,
\ldots, d_m$ can be chosen to be powers of prime numbers. A
character in the dual group $\fhat G$ is then given by
\begin{equation*}
\langle (\xi_1, \xi_2, \ldots, \xi_m),(x_1, x_2, \ldots ,
x_m)\rangle
    = \langle \xi_1, x_1\rangle \langle \xi_2, x_2 \rangle \ldots\langle \xi_m, x_m
    \rangle\,,
\end{equation*}
where $(\xi_1, \xi_2, \ldots, \xi_m)\in \widehat
\Z_{d_1}\times\widehat \Z_{d_2} \times \ldots \times \widehat
\Z_{d_m}\cong \widehat{G}. $ The discrete Fourier matrix $W_G$ for
$G = \Z_{d_1} {\times} \Z_{d_2}{\times} \ldots {\times} \Z_{d_m}$ is
chosen to be the Kronecker product of the Fourier matrices for the
groups $\Z_{d_1}, \Z_{d_2}, \ldots , \Z_{d_m}$, i.e., $W_G=W_{d_1}
\otimes W_{d_2} \otimes \ldots \otimes W_{d_m}$. For example, we
have
$$
W_{\Z_4}= \mysmallmatrix{4}{
1&1&1&1\\
1&i&-1&-i\\
1&-1&1&-1\\
1&-i&-1&i\\} \quad \text{and}\quad  W_{\Z_2 \times \Z_2}=
\mysmallmatrix{4}{
1&1&1&1\\
1&-1&1&-1\\
1&1&-1&-1\\
1&-1&-1&1\\
}.
$$
Note that for appropriately chosen bijections $S_1:\{ 0,1,\ldots,
|G|{-}1\} \longrightarrow G $ and $S_2: \{ 0,1,\ldots, |G|{-}1\}$
$\longrightarrow \widehat G$ we have $\widehat{f}\circ  S_2 =W_G (f
\circ S_1)$ for $f\in\C^{G}$.

\subsection{Short--time Fourier transforms on finite Abelian groups and Gabor frames}

For any $x\in G$, we define the translation operator $T_x$ as the
unitary operator on $\C^G$ given by $T_x f(y)=f(y{-}x)$, $y\in G$.
Similarly, we define the modulation operator $M_{\xi}$ for $\xi\in
\widehat{G}$ as the unitary operator defined by $M_\xi f=f \cdot\xi
$, where here and in the following $f \cdot g$ denotes the pointwise
product of $f,g\in\C^G$. Since $\widehat{M_\xi f}=T_\xi \widehat f$,
we refer to $M_\xi$ also as a frequency shift operator.

We set $\pi(\lambda)= M_\xi\circ T_x$ for $\lambda=(x,\xi)\in
G{\times}\widehat G$. The unitary operators $\pi(\lambda)$, $\lambda
\in G{\times}\widehat G$ are called time--frequency shift operator.

\begin{definition} The {\em short--time Fourier transformation}
$V_g:\C^G\longrightarrow \C^{G{\times}\widehat{G}}$ with respect to
the window $g\in \C^G{\setminus}\{0\}$ is given by
$$
    V_g f (x,\xi) = \langle f, \pi(x,\xi)g\rangle
        = \displaystyle\sum_{y \in G} f(y)\overline{g(y{-}x)}\overline{\langle \xi, y \rangle},\quad (x,\xi)\in
        G{\times}\widehat{G},
$$

\vindex \noindent
 where $f \in \C^G$.
\end{definition}

The inversion formula for the short--time Fourier transform is
\begin{eqnarray}
  f(y)=\tfrac 1 {|G|\, \|g\|_2^{2}} \sum_{(x,\xi)\in G{\times}\widehat{G}} V_gf(x,\xi)\,  g(y{-}x) \langle \xi, y
  \rangle\, , \quad y\in G \label{equation:STFTinversion},
\end{eqnarray}
i.e., $f$ can be composed of time--frequency shifted copies of any
$g\in\C^G{\setminus}\{0\}$. Further, $ \|V_g f\|_2=
\sqrt{|G|}\,\|f\|_2\|g\|_2$. This equation resembles
(\ref{equation:FourierInversion}), but the so-called Gabor system
$\{\pi(x,\xi)g\}_{(x,\xi)\in G{\times}\widehat G}$ is clearly not an
orthonormal basis if $|G|\neq 1$ since it consists of $|G|^2$
vectors in a $|G|$ dimensional space. As a matter of fact, such a
Gabor system is an equal norm tight frame which is defined below.

\begin{definition}
Let $G$ be a finite Abelian group and let $K$ be a finite or
countably infinite index set.  A family of functions
$\{\varphi_k\}\subset \C^G$ with
$$
    A\|f\|_2^2\, \leq\,  \sum_k |\langle f,\varphi_k \rangle|^2 \,
    \leq \,
    B\|f\|_2^2\, , \quad f\in\C^G,
$$

\vindex \noindent for positive $A$ and $B$ is called a {\em frame}
for $\C^G$. $A$ is called an {\em lower frame bound} and $B$ is
called a {\em upper frame bound} of the frame $\{\varphi_k\}$.

A frame is called {\em tight} if we can choose $A=B$. If we can
choose $A=B=1$, then the frame is called {\em Parseval tight frame}.
If $\|\varphi_k\|=C>0$ for all $k$, then the frame $\{\varphi_k\}$
is called {\em equal norm frame} and if in addition $C=1$, then we
have a {\em unit norm frame.}
\end{definition}

A direct consequence of (\ref{equation:STFTinversion}) is
\begin{proposition}\label{finite gabor frame}
For any $g\in \C^G{\setminus}\{0\}$, the collection
$\{\pi(\lambda)g\}_{\lambda\in G{\times}\widehat G}$ is an equal
norm  tight frame for $\C^G$ with frame bound $A=B=|G|\,\|g\|_2^2$.
\end{proposition}

The usefulness of frames stems largely from the existence of an
reconstruction formula similar to (\ref{equation:FourierInversion})
and (\ref{equation:STFTinversion}).
\begin{proposition}
  Let $\{\varphi_k\}$ be a frame for $\C^G$. Then exists a  so-called
  dual frame $\{\widetilde{\varphi}_k\}$, with
  \begin{equation}
    f=\sum_k \langle f,\varphi_k \rangle\widetilde{\varphi}_k
    =\sum_k \langle f,\widetilde{\varphi}_k \rangle\varphi_k\, , \quad f\in\C^G\,
    .\label{equation:FrameInversion}
  \end{equation}
\end{proposition}
Note that Parseval frames are self dual, i.e., we can choose
$\widetilde{\varphi}_k= \varphi_k$ for all $k$.

For additional material on frames and, in particular, Gabor frames
we refer to the excellent expositions \cite{C03,G01,KC06}. The
geometry of finite frames is discussed in \cite{BF03}.

For a given group $G$, we shall use again the previously defined
enumerations $S_2: \{ 0,1,\ldots, |G|{-}1\}$ $\longrightarrow
\widehat G$ and $S_1:\{ 0,1,\ldots, |G|{-}1\} \longrightarrow G$
which gave rise to the Fourier matrix $W_G$.  For $g\in \C^G$ and
$x\in G$, we define the $|G|{\times}|G|$--diagonal matrix
$$
D_{x,g} = \mysmallmatrix{4}{
  g(S_1(0)+x) &         &           &     0   \quad \quad \quad    \\
      & g(S_1(1)+x) &           &                  \\
      &         & \ddots    &                   \\
   0 \quad \quad  &         &           & g(S_1(|G|{-}1)+x)         \\
  }.
$$
Then, the $|G|{\times}|G|^2$--{\em full Gabor system matrix}
with respect to $g$ is given by
\begin{equation}\label{eqn:matrixA}
A_{G,g} = (D_{S_1(0),g}\cdot W_G\,|\,D_{S_1(1),g} \cdot W_G\,|
\,\cdots\, |\,D_{S_1(|G|{-}1),g} \cdot W_G)^\ast,
\end{equation}
where $M^\ast$ denotes the adjoint of the matrix $M$. For example,
for  $G=\Z_4$,
$$
    A_{\Z_4,(1,2,3,4)}:=
    \mysmallmatrix{16}{
            1& 1 &1  &1   & 2 &2 & 2&2 &    3 &3&3&3   &  4&4  &4&4      \\
            2& 2i&-2 &-2i & 3 &3i&-3&-3i & 4 &4i&-4&-4i &1&i&  -1  &-i \\
            3&-3 & 3 & -3 & 4 &-4&4&-4 &  1&-1&1&-1 &    2  &-2  &  2&-2  \\
            4&-4i&-4 &4i  & 1 &-i&-1&i& 2&-2i&-2&2i   & 3&-3i&  -3& 3i\\
       }^\ast\, .$$
Similarly, for the group $G= \Z_2 \times \Z_2$ we have
$$
    A_{\Z_2{\times}\Z_2,(1,2,3,4)}:=
    \mysmallmatrix{16}{
        1& 1 &1&1  &2 &2 & 2&2 &    3 &3&3&3   &  4&4  &4&4      \\
        2& -2&2 &-2 & 1 &-1&1&-1 & 4 &-4&4&-4 &3&-3&  3  &-3 \\
        3&3 & -3& -3 & 4 &4&-4&-4 &  1&1&-1&-1 &    2  &2  &  -2&-2  \\
        4&-4 &-4&4&  3 &-3&-3&3& 2&-2&-2&2   & 1&-1&  -1& 1\\
    }^\ast\, .
$$

Using the enumeration $S:\{0,1,\ldots,|G|^2{-}1\}\longrightarrow
G{\times}\widehat G$ which is given by the lexicographic order that
is induced by $S_1$ and $S_2$ on $G{\times}\widehat G$, we have $V_g
f\circ S=A_{G,g} f $. Therefore, we shall refer to $A_{G,g}$ as
short--time Fourier transform matrix with respect to the window $g$.
Clearly, the rows of $A_{G,g}$ represent the vectors in the Gabor
system $\{\pi(\lambda)g\}_{\lambda\in G{\times}\widehat G}$, and
(\ref{equation:STFTinversion}) implies that $A_{G,g}^\ast A_{G,g}$
is a multiple of the identity matrix.

\section{Uncertainty principles for the Fourier transform on finite Abelian groups}
\label{section:FourierUncertainty}


The following uncertainty theorem for functions defined on finite
Abelian groups is  the natural starting point for our discussion
\cite{DS89}.
\begin{theorem}\label{theorem:classicalUncertainty}
Let $f\in \C^G {\setminus} \{0\}$, then $\|f\|_0\cdot\|\widehat
f\|_0\geq|G|$.
\end{theorem}
\begin{proof}
 For $f\in \C^G$, $f\neq 0$, and without loss of generality
 $\|\widehat{f}\|_\infty=1$, we compute
 \begin{eqnarray*}
|G|&=& |G| \|\widehat{f}\|_\infty^2
    \ \leq \  |G| \left(\sum_{x\in G}|f(x)|\right)^2
    \  \leq \  |G| \|f\|_0\,\sum_{x\in G}|f(x)|^2\\
    &=& |G| \|f\|_0\,\frac 1{|G|} \sum_{\xi\in \widehat{G}}|\widehat f (\xi)|^2
     \ \leq \   \|f\|_0 \|\widehat f\|_0 \|\widehat{f}\|_\infty^2 \ = \ \|f\|_0\,\|\widehat f\|_0.
\end{eqnarray*}

\vspace{-1cm}
\end{proof}

A complementary result characterizes those $f$ for which the bound
in Theorem~\ref{theorem:classicalUncertainty} is sharp
\cite{DS89,MOP04}.

\begin{proposition}\label{proposition:DonohoStark}

\vspace{-.3cm}
\begin{enumerate}
  \item If $k$ divides $|G|$, then there exists $f\in \C^G$ with $\|f\|_0=k$ and
  $\|\widehat f\|_0=\frac {|G|}  k$.
  \item If $\|f\|_0\|\widehat f\|_0 = |G|$ and $e \in \supp f$, then $\supp f$ is a subgroup of $G$.
\end{enumerate}
\end{proposition}


\subsection{Groups of prime order}\label{section:taobiroFT}

The geometric mean of two positive numbers is dominated by their
arithmetic mean; hence, Theorem~\ref{theorem:classicalUncertainty}
implies the weaker inequality
\begin{equation}\label{equation:classicaluncertaintygeommean}
    \|f\|_0+\|\widehat f\|_0\geq 2 \sqrt{|G|}.
\end{equation}

If $|G|$ is prime, i.e., if $G$ is a cyclic group of prime order,
then \eqref{equation:classicaluncertaintygeommean} and also
Theorem~\ref{theorem:classicalUncertainty} can be improved
significantly \cite{Fre04,Tao05}.

\begin{theorem}\label{theorem:tao}
Let $G=\Z_p$ with $p$ prime. Then $\|f\|_0+\|\widehat{f}\|_0\geq
|G|{+}1$ holds for all $f\in \C^G {\setminus} \{0\}$.
\end{theorem}

This result is a direct consequence from Chebotarev's Theorem which
states that every minor of the  Fourier transform matrix $W_{\Z_p}$,
$p$ prime, is nonzero \cite{EI76,SL96,Tao05,Fre04}. In fact, to
obtain Theorem~\ref{theorem:tao} we only need to combine
Chebotarev's Theorem with

\begin{proposition}\label{proposition:nominors}
  Let $M\in\C^{m{\times}n}$. Then $\displaystyle \|f\|_0+\|Mf\|_0\geq m {+}1$ for all $f\in \C^n$
  if and only if every minor of $M$ is nonzero.
  Moreover, if every minor of $M\in \C^{m{\times}n}$ is nonzero
  and $k,l$ are given with $k+l\geq m{+}1$, then there exists $f\in \C^n$
  with $\|f\|_0=k$ and $\|Mf\|_0=l$.
\end{proposition}

\begin{lemma}\label{lemma:minorranks}
For $M\in \C^{m\times n}$ and $1\leq k\leq m$, $1\leq l\leq n$,
there exists $f\in \C^n$ with  $\|f\|_0=k$ and $\|Mf\|_0=l$ if and
only if there exist sets $A\subseteq\{0,\ldots,n{-}1\}$ and
$B\subseteq\{0,\ldots,m{-}1\}$ with $|A|=k$, $|B|=m-l$, and for all
$a\in A$ and $y\in B^c$, we have
    \begin{eqnarray}
      \rank M_{A{\setminus} \{a\},B}\ =\ \rank M_{A,B}\ =\ \rank
      M_{A,B\cup\{y\}}-1\ <\ |A|\label{equation:minorranks}\, .
    \end{eqnarray}
\end{lemma}

\vspace{-.5cm}

{\it Proof of Proposition~\ref{proposition:nominors}.}
 If $f$ has no zero minors, then (\ref{equation:minorranks}) in
Lemma~\ref{lemma:minorranks} is equivalent to $|B|<|A|$, implying
that there exists $f\in \C^n$
  with $\|f\|_0=k$ and $\|Mf\|_0=l$ if and only if $k+l\geq m{+}1$.

It remains to show that $\displaystyle \|f\|_0+\|Mf\|_0\geq m {+}1$
for all $f$ implies that $M$ has no zero minors.  To this end,
assume that there is a $d\times d$ submatrix $M_{A,B}$ of $M$ with
$\det M_{A,B}=0$. Then there exists a nonzero vector   $f'\in \C^A$
such that $M_{A,B} f'=0$. For $f=i_A f'$, $\|Mf\|_0\leq m-d$ and
therefore $\|f\|_0+\|Mf\|_0\leq d+m-d=m<m{+}1$. \hfill $\Box$

%
%


%

Theorem~\ref{theorem:tao} is a clear improvement to
Theorem~\ref{theorem:classicalUncertainty} but it applies only to
cyclic groups of prime order. In fact, any other finite Abelian
group $G$ has proper subgroups which lead to zero minors in $W_G$.
As example, we display in Table~\ref{table:ranksOfMinorsFourier}
counts on the ranks of square submatrices of $W_{\Z_5}$ and
$W_{\Z_6}$. Due to their role in obtaining
Theorem~\ref{theorem:tao},   we shall now collect facts regarding
zero and nonzero minors of Fourier matrices in general.

\begin{table}[th]
\begin{center}
\begin{minipage}[t]{4.8cm}
   {\footnotesize     \begin{tabular}[t]{|c||c|c|c|c|c|}
 \hline
            & 1 & 2 &3&4 &5   \\
          \hline
          \hline
           1& 25 &0  & 0 &  0&    0 \\
          \hline
           2&  0& 100 & 0 &0  &   0 \\
          \hline
           3& 0 &  0&  100& 0 &   0 \\
          \hline
           4& 0 & 0 & 0 &  25 &  0   \\
          \hline
           5&  0&  0 &  0& 0 & 1   \\
          \hline
        \end{tabular} }
        \end{minipage}
\hspace{2cm}
        \begin{minipage}[t]{6cm}
   {\footnotesize     \begin{tabular}[t]{|c||c|c|c|c|c|c|}
\hline
          & 1 & 2 &3&4 &5 &6  \\
          \hline
          \hline
           1& 36 & 36 & 0 &  0& 0 &  0 \\
          \hline
           2& 0 & 189 & 48 & 0 & 0&0   \\
          \hline
           3& 0 & 0 &  352& 36  & 0&  0 \\

          \hline
           4&0  &  0&  0&  189 & 0 & 0  \\
          \hline
           5& 0 & 0  & 0 & 0 & 36& 0  \\
          \hline
           6& 0 & 0 &  0&  0 & 0& 1  \\
           \hline
        \end{tabular} }
        \end{minipage}
    \caption{\color{black} Counts of the numerically computed rank of submatrices of  $W_{\Z_5}$ and $W_{\Z_6}$.
     The column index is the size of square submatrices considered, and
     the row index corresponds to their ranks.
     \label{table:ranksOfMinorsFourier}}\color{black}
\end{center}
\end{table}

Let $M\in\C^{n{\times}n}$ and let $A, B \subset \{1,2,\ldots, n\}$
such that $|A|=|B|$. Then $\det M_{A,B}$ defines a minor of $M$, and
 $\det M_{A^c, B^c}$ is called its {\it complementary minor}.

\begin{proposition}\label{proposition:DFTminors}

\vspace{-.2cm}
\begin{enumerate}
\item The complementary minor of any zero minor in a Fourier matrix $W_G$ is also zero.

\item Let $d_0 > 1$ be the smallest divisor of $|G|$. Then for all $d_0 \le r \le n{-}d_0$,
there exists an $r{\times}r$ zero minor of the Fourier matrix
$W_{G}$. In particular, if $|G|$ is even, then there exist
$r{\times}r$ zero minor for $r=2,3,\ldots, |G|{-}2$.

\item Any minor of the Fourier matrix $W_{\Z_n}$, $n\in\N$, that contains only adjacent rows
or columns is nonzero.

\end{enumerate}
\end{proposition}
\begin{proof}\quad {\it 1.}
The adjoint of a matrix $M=(m_{kl})$ is $\adj M =(M_{kl})$, where
$M_{kl}=(-1)^{k+l}\det M_{\{k\}^c , \{l\}^c}$ is the cofactor of the
element $m_{kl}$. Then for any sets $A,B$ of cardinality $r$,
Jacobi's theorem states that
\begin{equation}\label{eq:jacobi}
\det M_{A,B}=(-1)^r \det (\adj M)_{A^c,B^c}\cdot (\det M)^{r{-}1}
\end{equation}
Furthermore, $\adj M\cdot M =\det M \cdot I$ \cite{Pra94}.

For any zero minor of $M=W_G$ on the left hand side of
\eqref{eq:jacobi}, Jacobi's theorem implies that the right hand
side, representing a minor in $\adj W_G$, is zero as well. Since
$W_G \cdot \overline{W_G}=|G| \cdot I$, we have $\adj W_G =
\frac{\det(W_G)}{|G|} \cdot \overline{W_G}$. Thus the corresponding
minor in $\overline{W_G}$ is zero, which implies that also the
corresponding minor in $W_G$ is zero.

{\it 2.} Let $d$ divide $|G|$. Part {\it 1} in Proposition
\ref{proposition:DonohoStark} allows us to choose $f_d$ such that
$\|f_d\|_0= d$ and $\|\fhat f_d\|_0= \tfrac {|G|}  d$. Hence, for
any $r$ with $d\leq r \leq |G| {-}\tfrac {|G|} d$ we can pick sets
$A \supseteq \supp f_d$ and $B \subseteq (\supp\fhat f_d)^c$ such
that $|A|=|B|=r$. Then $r_A f_d\in \ker M_{A,B}$ and the $r {\times}
r$-minor $\det M_{A,B}$ is zero.

This way, we obtain $r {\times} r$ zero minors for $d_0\leq r \leq
\tfrac {|G|} {d_0} (d_0-1)$ and for $\tfrac {|G|} {d_0} \leq r \leq
|G|- {d_0}$, where $d_0$ is the smallest nontrivial divisor of
$|G|$. The result follows since $d_0-1\geq 1$.

%

{\it 3.}  A minor with adjacent columns is a determinant of the type
\begin{eqnarray*}
\det\mysmallmatrix{4}{
\omega^{k_1 l} & \omega^{k_1(l{+}1)} & \cdots & \omega^{k_1(l+m)} \\
\omega^{k_2 l} & \omega^{k_2(l{+}1)} & \cdots & \omega^{k_2(l+m)} \\
\vdots& \vdots & \cdots & \vdots\\
\omega^{k_m l} & \omega^{k_m(l{+}1)} & \cdots & \omega^{k_m(j+m)} }
&=&\omega^{k_1 l+k_2 l+\dots+k_m l} \ \det\mysmallmatrix{4}{
1 & \omega^{k_1} & \cdots & \omega^{m k_1} \\
1 & \omega^{k_2} & \cdots & \omega^{m k_2} \\
\vdots& \vdots & \cdots & \vdots\\
1 & \omega^{k_m} & \cdots & \omega^{m k_m}
}\\
&=&\omega^{k_1 l+k_2 l+\dots+k_m l}\prod\limits_{i<j\leq
m}(\omega^{k_{j}}-\omega^{k_{i}}) \neq 0
\end{eqnarray*}
The second determinant was evaluated using the formula for
Vandermonde determinants and the result does not equal $0$, as
always $i<j$ and $\omega$ is a primitive $n$-th root of unity.
\end{proof}

\subsection{Groups of non-prime order}\label{section:meshulamFT}

Meshulam improved the bound in the classical uncertainty relation
presented in Theorem~\ref{theorem:classicalUncertainty} for most
finite Abelian groups of non-prime order \cite{Mes05}. He  defines
for $0<k\leq |G|$ the function
$$ \theta(G,k)=\min\big\{\|\fhat{f}\|_0:\ f\in \C^G\text{ and } 0<  \|f\|_0\leq k\big\}\,.$$
Note that Theorem~\ref{theorem:tao} implies that
$\theta(\Z_p,k)=p-k+1$. The main result in \cite{Mes05} is

\begin{theorem}
\label{theorem:meshulam}
For
$k\leq |G|$, let $d_1$ be the largest divisor of $|G|$ which is less
than or equal to $k$ and let $d_2$ be the smallest divisor of $|G|$
which is larger than or equal to $k$. Then
\begin{eqnarray}
\theta(G,k) \geq
\frac{|G|}{d_1d_2}(d_1+d_2-k).\label{equation:meshulamMainResult}
\end{eqnarray}
\end{theorem}

\vpropend Tao realized that this theorem simply states that all
possible lattice points $(\|f\|_0, \|\widehat f\|_0)$ lie in the
convex hull of the points $(|H|,|G/H|)$, where $H$ ranges over all
subgroups of $G$ \cite{Mes05}. To see this, recall that for any
divisor $d$ of $|G|$  exists a subgroup $H$ of  $G$ with $d=|H|$.
Furthermore, the right hand side of expression
 \eqref{equation:meshulamMainResult} is linear between two successive divisors
 and the slope is increasing when $k$
increases. Hence \eqref{equation:meshulamMainResult} characterizes
the convex hull of the  points $(|H|, |G|/|H|)$.
Proposition~\ref{proposition:DonohoStark}, part {\it 1}, implies
that the vertex points $(|H|, |G|/|H|)$ are attained.

The proof of Theorem~\ref{theorem:meshulam} in \cite{Mes05} is
inductive and uses three facts: first, it uses
Theorem~\ref{theorem:tao} as induction seed, and second, it uses the
submultiplicativity of the right hand side of
(\ref{equation:meshulamMainResult}). That is, if we denote this
right hand side by $u(n,k)$ for $n=|G|$, then it uses that
$u(n,k)\leq u(\tfrac n d, t)u(d,s)$ for $d$ dividing $n$ and $st \le
k$. The third ingredient is

\begin{proposition}
\label{proposition:meshulam-induction-argument}
Let
$H$ be a subgroup of $G$. For $k\leq |G|$ there exist $s\leq q$,
$t\leq p$ with $st \le k$ and
$$
\theta (G,k) \ge \theta (H,s)\,\theta(G/H, t)\,.
$$
\end{proposition}

\vpropend Meshulam's proof of
Proposition~\ref{proposition:meshulam-induction-argument} is heavy
on algebraic notation and does not give good insight from the point
of view of Fourier analysis. For this reason, and for completeness
sake, we give a streamlined version of Meshulam's proof of
Proposition~\ref{proposition:meshulam-induction-argument}. See also
\cite{LM05} for an elegant and non-inductive proof of
Theorem~\ref{theorem:meshulam}.

But first, note that if $G \cong H \times G/H$, then
Proposition~\ref{proposition:meshulam-induction-argument}  can be
proven using the fact that then $\widehat G \cong \widehat H \times
\widehat{G/H}$, and, therefore, $\widehat{f}$ can be calculated by
performing two partial Fourier transforms. For example, such
argument can be applied to $G= \mathbb{Z}_m \times \mathbb{Z}_n
\cong \mathbb{Z}_{mn}, \,\gcd(m,n)=1$, and $H=\Z_m{\times}\{e\}$.
Even simpler is the special case discussed in
Proposition~\ref{proposition:productset}. We state and prove this
result to illustrate the main idea used to prove
Proposition~\ref{proposition:meshulam-induction-argument}.

\begin{proposition}\label{proposition:productset}
Let $A_1\subseteq G_1$ and $A_2\subseteq G_2$ and $f\in
\C^{G_1{\times}G_2}$ be given with $\supp f\subseteq
A_1{\times}A_2$. Then $\displaystyle
        \|\widehat{ f}\|_0\geq \theta(G_1,|A_1|)\, \theta(G_2,|A_2|)
$.
\end{proposition}
\begin{proof}
We picture $f$ as a $|G_1|{\times}|G_2|$ matrix and note that $\supp
f \subseteq A_1 {\times} A_2$ implies that $f$ has exactly
$|G_2{\setminus} A_2 |$ zero columns and $|A_2|$ columns with at
least  $|G_1{\setminus} A_1 |$ zeros.

The function $\mathcal F_1 f$ is obtained by applying the
$G_1$--Fourier transformation to each column. Hence, $\mathcal F_1
f$ has $|G_2 {\setminus} A_2 |$ zero columns and, at most,
$|G_1|-\theta(G_1, |A_1|)$ zeros in the remaining $A_2$ columns. It
is easy to see that in the scenarios which leads to the weakest
bound for $\|\widehat f\|_0$, we have $|G_1|-\theta(G_1, |A_1|)$
zeros in each of these $|A_2|$ columns and that they are lined up to
form $|G_1|-\theta(G_1, |A_1|)$ zero rows in $\mathcal F_1 f$.  In
this case, the remaining $\theta(G_1, |A_1|)$ rows contain exactly
$|G_2{\setminus} A_2 |$ zeros, i.e., $|A_2|$ nonzero elements.

Now, we calculate $\mathcal  F f$ by taking a $G_2$--Fourier
transform along each row of  $\mathcal F_1 f$.  As a result,
$|G_1|-\theta(G_1, |A_1|)$ zero rows remain, and in the other
$\theta(G_1, |A_1|)$ rows, at least $\theta(G_2,|A_2|)$ zeros are
present. We conclude that
\begin{eqnarray*}
        \|\widehat{f}\|_0 &\geq& \theta(G_1,|A_1|)\,\theta(G_2,|A_2|).
\end{eqnarray*}

\vspace{-1cm}
\end{proof}

The property that the $G= G_1{\times}G_2$--Fourier transformation
``splits" into a $G_1$--Fourier transformation and a $G_2$--Fourier
transformation is the basis of the simple proof of
Proposition~\ref{proposition:productset}. In the proof of
Proposition~\ref{proposition:meshulam-induction-argument} we shall
see that the general case follows from small adjustments to the
arguments used to prove Proposition~\ref{proposition:productset}.


{\it Proof of
Proposition~\ref{proposition:meshulam-induction-argument}.}  Let
$H=\{x_i\}$ be a subgroup of $G$ and, abusing notation, we let
$\{x_j\}$ be a set of coset representatives of the quotient group
$G/H$. Then each element in $G$ has a unique representation as
$x_i{+}x_j$. We let $H^\perp$ denote the characters $\{\xi_j \in
\widehat G : \xi_j (H) =1\}$. $H^\perp$ is a subgroup of $\widehat
G$, and we denote by $\{\xi_i\}$ a set of coset representatives of
the quotient group $\widehat{G}/H^\perp$. Every element $\xi \in
\widehat G$ has a unique decomposition as $\xi_i{+} \xi_j$.

The Pontryagin duality theorem implies  $\widehat{G}/H^\perp \cong
\widehat H$. This allows us to assign a character $\xi_i'\in
\widehat H$ to each $\xi_i\in \widehat{G}/H^\perp$ with
$\xi_{i_1}'{+}\xi_{i_2}'= (\xi_{i_1}{+}\xi_{i_2})'$
\cite{Kat76}.\footnote{In particular, in the case
$G=\mathbb{Z}_{mn}$, $\gcd(m,n)=1$, $\mathbb{Z}^\perp_m \cong
\mathbb{Z}_n$ and $\mathbb{Z}^\perp_m \cong \mathbb{Z}_n$.} Further,
$\langle \xi_i,x_i\rangle_G=\langle \xi_i',x_i\rangle_H$ for all
$x_i\in H$ and all $\xi_i\in \widehat{G}/H^\perp$. Similarly, we use
$\widehat{G/H}\cong {H^\perp}$ to assign to each $\xi_j$ an element
$\xi_j'\in \widehat{G/H}$ with $\langle \xi_j,x_j\rangle_G= \langle
\xi_j', x_j{+}H \rangle_{G/H}$ for all $x_j$.

For  $f\in \C^G$ and any $\xi=\xi_i{+}\xi_j\in \widehat G$, we
calculate
\begin{eqnarray*}
\widehat{f}(\xi)=\widehat{f}(\xi_i{+}\xi_j) &=&
\displaystyle\sum_{x_j }\displaystyle\sum_{x_i} f(x_i{+}x_j)
\overline{\langle\xi_i{+}\xi_j, x_i{+}x_j\rangle}_G\\
 &=&
\displaystyle\sum_{x_j }\displaystyle\sum_{x_i } f(x_i{+}x_j)
\overline{\langle \xi_i,x_i\rangle}_G\overline{\langle\xi_i,
x_j\rangle}_G\overline{\langle\xi_j,x_i\rangle}_G\overline{\langle\xi_j,x_j\rangle}_G\\
&=&\displaystyle\sum_{x_j }\biggl({\displaystyle\sum_{x_i }
f(x_i{+}x_j)\overline{\langle\xi_i',x_i\rangle}_H\biggr)
\overline{\langle\xi_i, x_j\rangle}_G}\ \
\overline{\langle\xi_j',x_j{+}H\rangle}_{G/H}
\end{eqnarray*}
where the last equality follows since $\xi_j \in H^\perp$ implies
$\langle\xi_j, x_i\rangle_G =1$.

\begin{figure}[th]
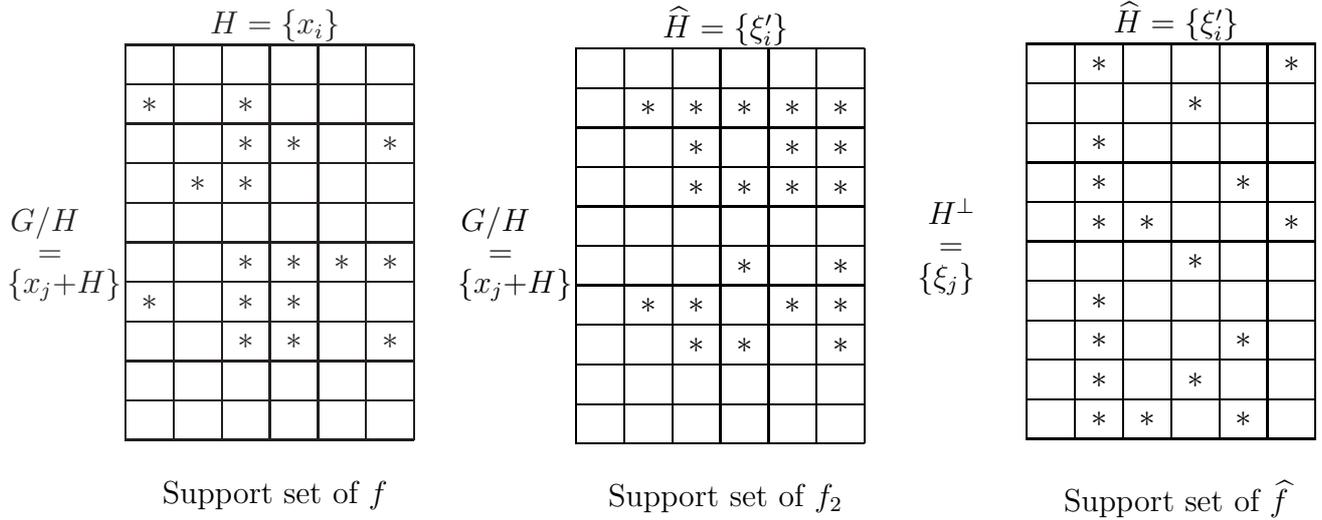

    \begin{center}
 \begin{minipage}[c]{1cm}
   \begin{center}

   \vspace{-.8cm}
 $\ G/H$

\vspace{-.1cm}

  $\ =$

  \vspace{-.1cm}

  $\{x_j{+}H\}$
 \end{center}
 \end{minipage}
        \begin{minipage}[c]{4.8cm}
        \begin{center}
        $H=\{x_i\}$

        \begin{tabular}{|c|c|c|c|c|c|}
          \hline
&&&&&\\
\hline
           $\ast$ &&$\ast$&&  &\\
           \hline
           &  &$\ast$ &$\ast$&&$\ast$\\
           \hline
             &$\ast$&$\ast$&&&    \\
          \hline
          &&&&&\\
          \hline
           &&$\ast$&$\ast$& $\ast$ &$\ast$\\
           \hline
               $\ast$& & $\ast$&$\ast$&&    \\
               \hline
               &&$\ast$&$\ast$&&$\ast$\\
\hline
&&&&&\\
          \hline
           &  &  & && \\
          \hline
        \end{tabular}
        \color{black}\\[.5cm]

        Support set of $f$
        \end{center}\color{black}
        \end{minipage}\color{black}
 \begin{minipage}[c]{1cm}
   \begin{center}

   \vspace{-.8cm}
 $\ G/H$

\vspace{-.1cm}

  $\ =$

  \vspace{-.1cm}

  $\{x_j{+}H\}$
 \end{center}
 \end{minipage}
        \begin{minipage}[c]{4.8cm}
        \begin{center}
         $\widehat H=\{\xi_i'\}$

        \begin{tabular}{|c|c|c|c|c|c|}
          \hline
          &&&&&\\
          \hline
           & $\ast$ & $\ast$ & $\ast$&$\ast$  &  $\ast$  \\
          \hline
          &&$\ast$&&$\ast$&$\ast$\\
          \hline
         \color{white} $\ast$\color{black} & & $\ast$ &  $\ast$&  $\ast$ &  $\ast$   \\
          \hline
          &&&&&\\
          \hline
 &  &  &  $\ast$&   &  $\ast$   \\
 \hline
           & $\ast$ & $\ast$& \color{white}$\ast$\color{black} & $\ast$ & $\ast$  \\
          \hline
           & & $\ast$ &  $\ast$&   &  $\ast$   \\
           \hline
           &&&&&\\
          \hline
         &&&&&\\
          \hline
        \end{tabular}
        \color{black}\\[.5cm]

        Support set of $f_2$
        \end{center}\color{black}
        \end{minipage}\color{black}
 \begin{minipage}[c]{1cm}
   \begin{center}

   \vspace{-.8cm}
 $\ H^\perp$

\vspace{-.1cm}

  $\ =$

  \vspace{-.1cm}

  $\{\xi_j\}$
 \end{center}
 \end{minipage}
        \begin{minipage}[c]{4.8cm}
        \begin{center}
        $\widehat H=\{\xi_i'\}$

        \begin{tabular}{|c|c|c|c|c|c|}
          \hline
           \color{white}$\ast$\color{black}&$\ast$&&&&$\ast$\\
           \hline
            \color{white}$\ast$\color{black} &  &\color{white}$\ast$\color{black} &$\ast$ & &\\
            \hline
            \color{white}$\ast$\color{black}& $\ast$ &&&& \\
          \hline
         \color{white}$\ast$\color{black} & $\ast$ &  &  &  $\ast$  &   \\
           \hline
          \color{white}$\ast$\color{black}&$\ast$&$\ast$&&&$\ast$\\
          \hline
           \color{white}$\ast$\color{black}&&&$\ast$& &\\
           \hline
           \color{white}$\ast$\color{black}&$\ast$&&&&\\
           \hline
            \color{white}$\ast$\color{black}   & $\ast$  &  & & $\ast$ &\\
\hline
            \color{white}$\ast$\color{black}&$\ast$&&$\ast$ &&   \\
          \hline
           \color{white}$\ast$\color{black}&$\ast$& $\ast$ &   & $\ast$ &  \\
          \hline
        \end{tabular}
        \color{black}\\[.5cm]
        Support set of $\widehat f$
        \end{center}\color{black}
        \end{minipage}\color{black}
    \end{center}\color{black}\caption{\color{black} Illustration of the proof of
    Proposition~\ref{proposition:meshulam-induction-argument} for
      $G=\Z_{10} \times \Z_6$ and $k=17$.
The function $f_2$ is obtained by the application of $H$--Fourier
transformations to the rows of $f$ which is succeeded by an unitary
multiplicatiton operator . To calculate $\fhat f$ we apply
$G/H$--Fourier transformations to the columns of $f_2$. For clarity,
we choose synthetic support sets of $f$, $f_2$, and $\widehat{f}$.
Here $t=6$ and $s=\lfloor\tfrac {17} 6 \rfloor=2.$
}\color{black}\label{figure:exampleFT}
\end{figure}

We set $f_1( \xi_i',x_j) := \displaystyle\sum_{x_i \in H} f(x_i
+x_j)\overline{\langle\xi_i', x_i\rangle}_H$, which, for fixed
$x_j$, is the $H$--Fourier transform $\mathcal{F}_H$ on the coset
$x_j{+}H$ in $G$, and $f_2( \xi_i',x_j)=f_1(
\xi_i',x_j)\overline{\langle\xi_i, x_j\rangle}_G$. Further $f_1$ and
$f_2$ have the same support sets. We summarize that $\widehat f$ can
be obtained from $f$ via two partial Fourier transformations and an
enclosed  unitary multiplication operator, as illustrated in
Figure~\ref{figure:exampleFT}.

Let us now fix $f\in \C^G$ with $\|f\|_0 \le k$ and $\|\widehat
f\|_0= \theta(G,k)$.

Let $t := | \{x_j: \supp f \cap \left(x_j{+}H \right)\ne
\varnothing\}|$. Note that the support of $f$ contains at most $k$
elements which are distributed among $t$ cosets of $H$. Hence, there
must be a coset $x_{j_0}{+} H$ which contains  $s'\leq
s=\lfloor\tfrac  k t \rfloor$ elements of $\supp f$. Therefore,
$$
    \|f_2(\cdot, x_{j_0})\|_0
        =\|\mathcal{F}_H f (\cdot + x_{j_0})\|_0\ge \theta (H, s') \ge \theta (H, s)
$$
This implies that $\Xi=\{\xi_i\in \fhat G/H^\perp
:f_2(\xi_i',\cdot)\not\equiv 0\}$ satisfies $|\Xi|\geq \theta(H,s)$.
In fact, the definition of $t$ implies that for $\xi_i\in \Xi$, we
have $0<\supp f_2(\xi_i',\cdot)\leq t$.
We conclude
$$
    \theta(G, k)
        = \|\widehat f\|_0
        = \displaystyle\sum_{\xi_i } \|\mathcal{F}_{G/H} f_2(\xi_i',\,\cdot )\|_0
        \geq \sum_{\xi_i \in \Xi } \|\mathcal{F}_{G/H} f_2(\xi_i' ,\,\cdot )\|_0
        \ge \theta (H,s) \theta (G/H, t).
$$

\vspace{-1.47cm} \hfill $\Box$\\[.1cm]


Next, we  discuss the question whether the inequality
(\ref{equation:meshulamMainResult}) in
Theorem~\ref{theorem:meshulam} is sharp, or, more precisely, we
shall check whether for some given Abelian group $G$ and $(k,l)$
chosen with $l\geq \theta(G,k) \geq \tfrac{|G|}{d_1d_2}(d_1+d_2-k)$
there exists a function $f\in \C^G$ with $\|f\|_0=k$ and
$\|\widehat{f}\|_0=l$. This question has been discussed earlier for
$G=\Z_6$ and $G=\Z_8$ in \cite{FKLM05}.


The following affirmative partial result follows from  the proof of
Proposition 4.5 in \cite{Kut03}.

\begin{proposition}\label{proposition:Gitta}
  If  $0<k,l\leq|G|$ satisfy $l+k\geq |G|{+}1$, then there exists a function $f\in \C^G$ with $\|f\|_0=k$ and
    $\|\widehat{f}\|_0=l$.
\end{proposition}

\begin{figure}\begin{center}
  \includegraphics[width=9cm]{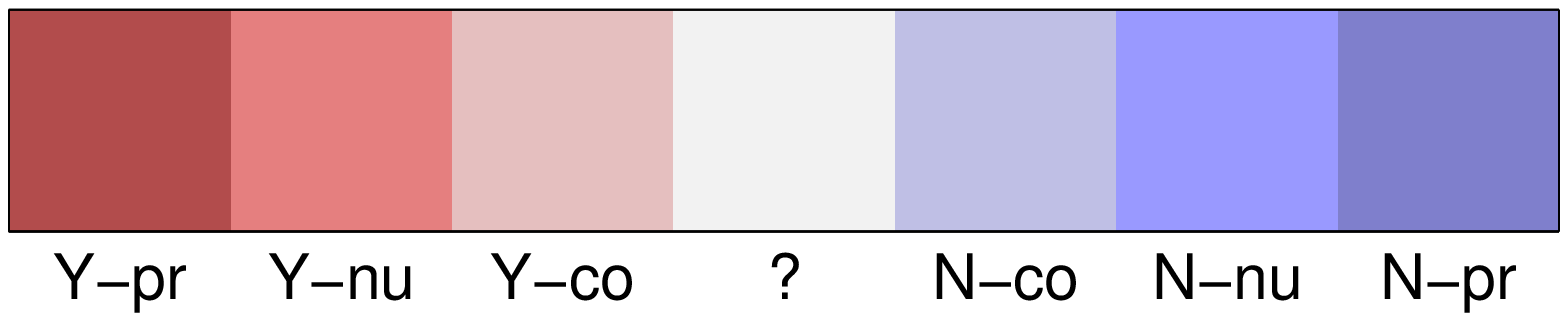}
\end{center}

\vspace{-.6cm}
\begin{caption}
 {Color coding which is used in
 Figures~\ref{figure:possiblepairsnoncyclic}--\ref{figure:possiblepairsnoncyclicVg} to describe subsets of
 $\N^2$ or $\N^3$. The color determines whether a given
 value is in the set under discussion.
{\sf Y-pr} indicates that their is proof that the corresponding
value is in the set considered. {\sf Y-nu} implies that their is
numerical evidence that the value is in the set and {\sf Y-co}
indicates that we conjecture that the value is in the set. {\sf
N-pr} indicates that their is proof that the corresponding value is
not in the set, and {\sf N-nu} and {\sf N-co} are defined
accordingly. The color adjacent to {\sf ?} implies that no judgement
is made here.
  }\label{figure:ColorMap}
   \end{caption}

\end{figure}

\begin{figure}
\begin{center}
{   \includegraphics[width=5cm]{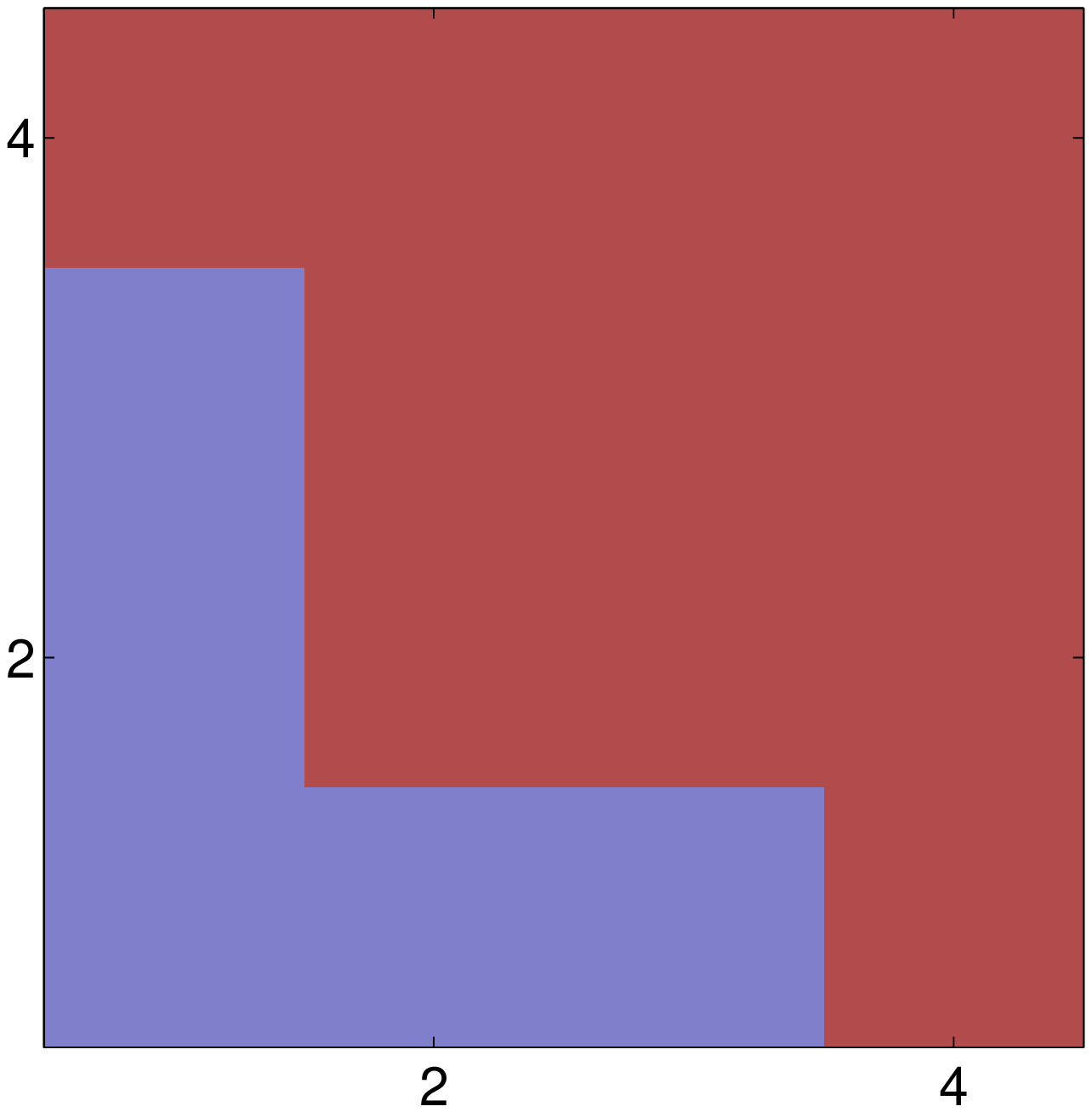} \
   \includegraphics[width=5cm]{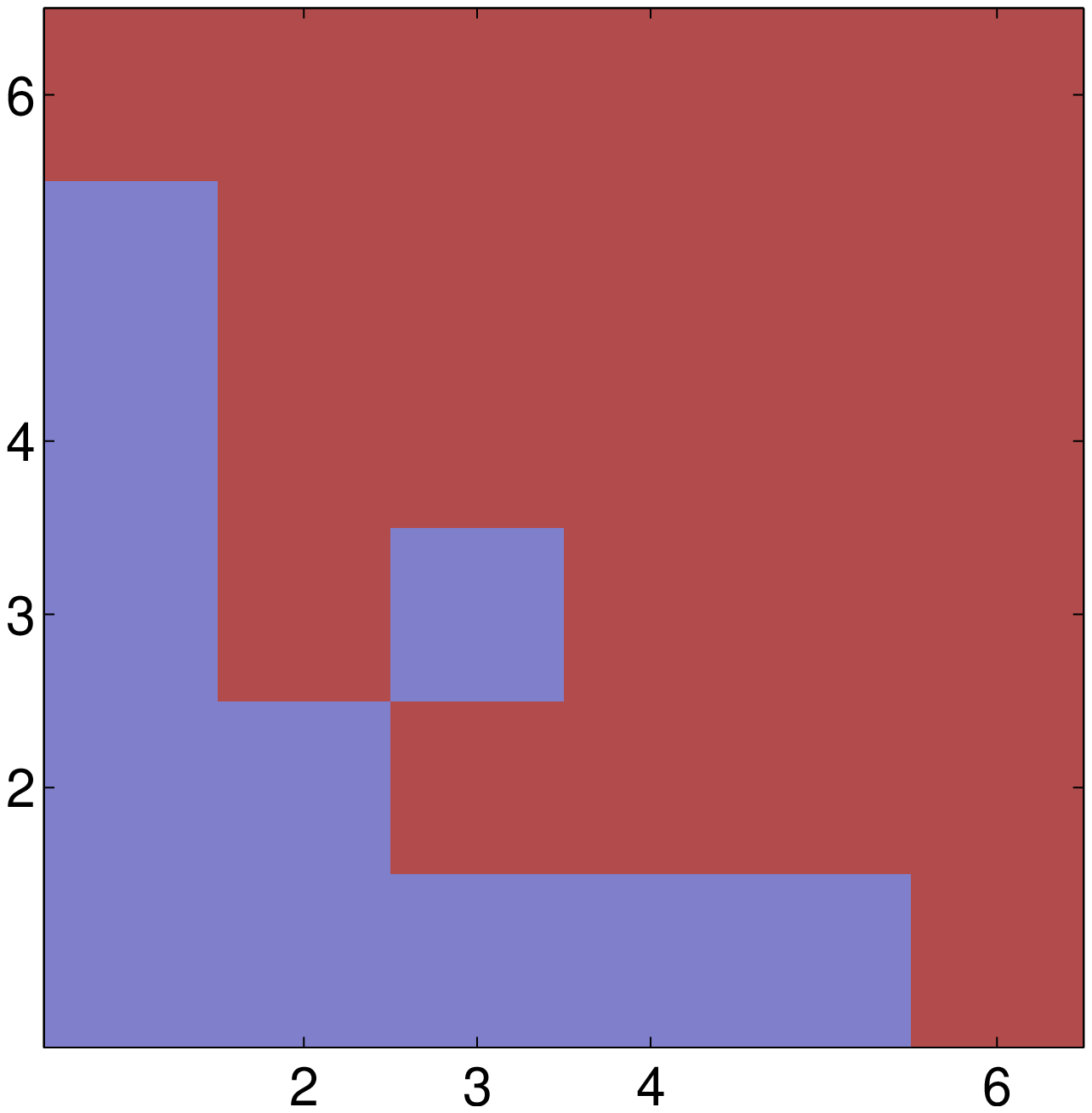} \
      \includegraphics[width=5cm]{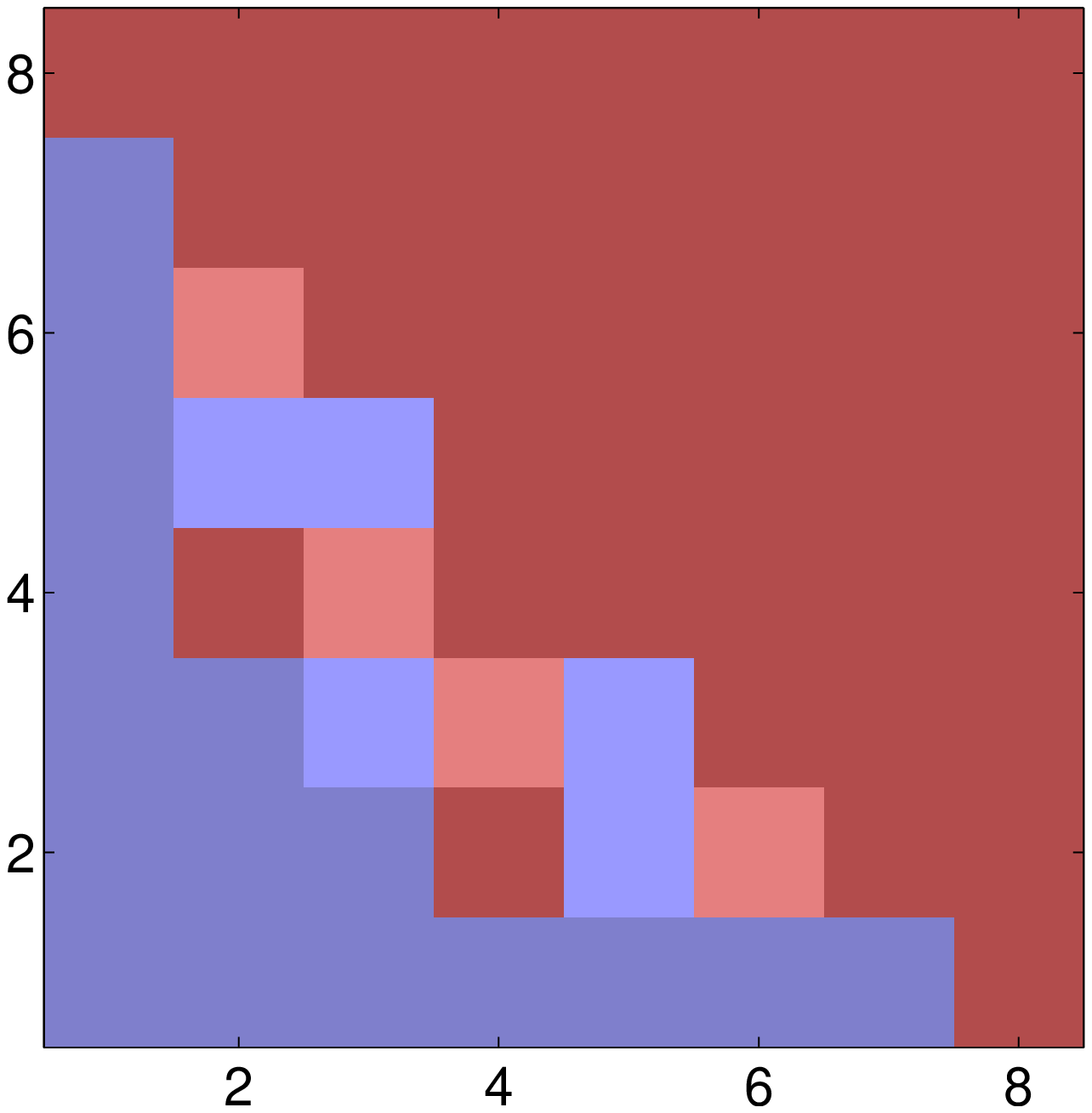}

\vspace{-4.7cm}

\hspace{1.2cm} \color{white}
  \hspace{1cm}  \begin{tabular}{p{4.65cm}p{4.65cm}p{4.65cm}}
       $\Z_4$ & $\Z_6\cong\Z_2{\times}\Z_3$ & $\Z_8,\ \Z_2{\times}\Z_4$
    \end{tabular}\color{black}

\vspace{4.2cm}

   \includegraphics[width=5cm]{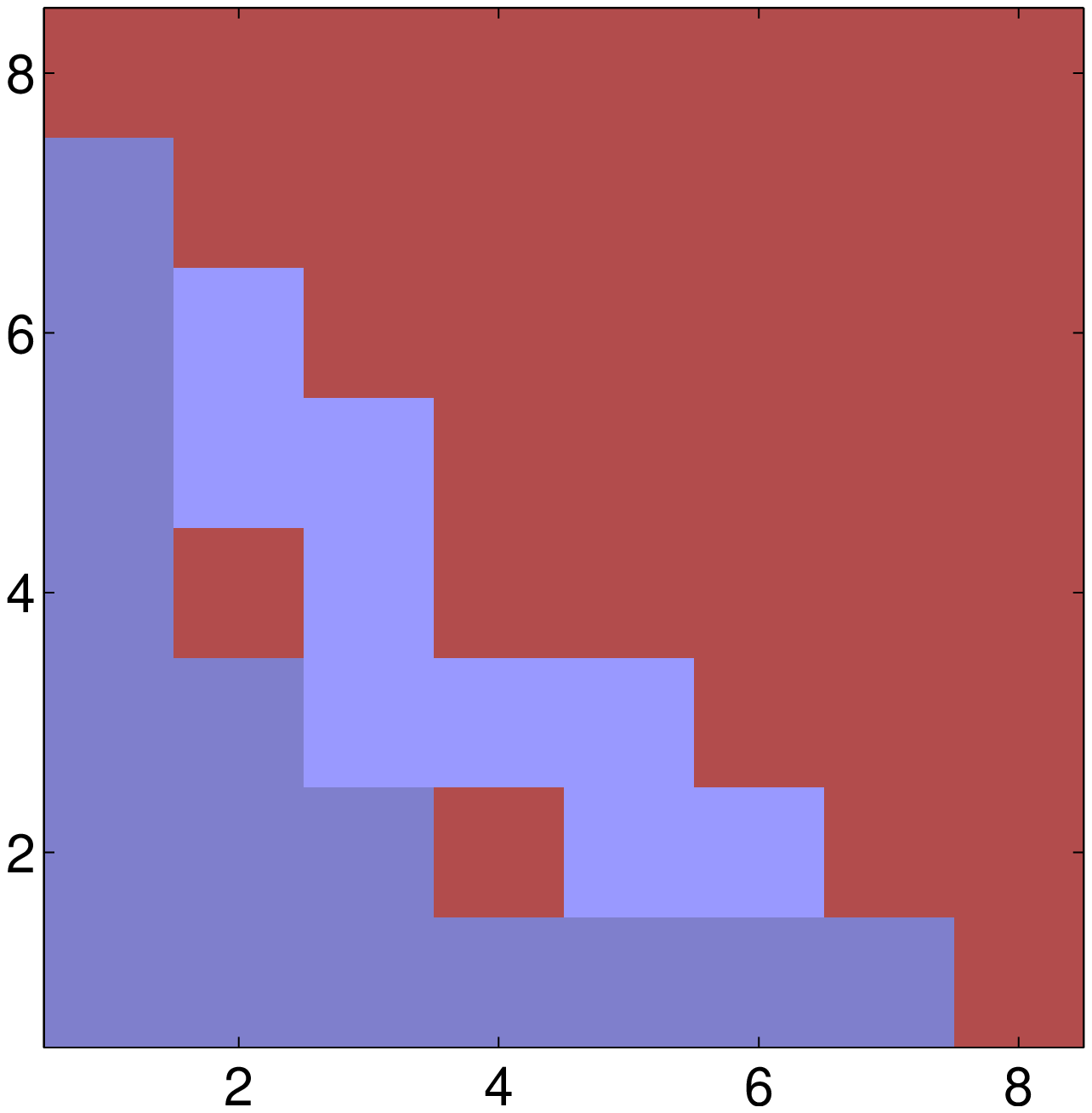} \
   \includegraphics[width=5cm]{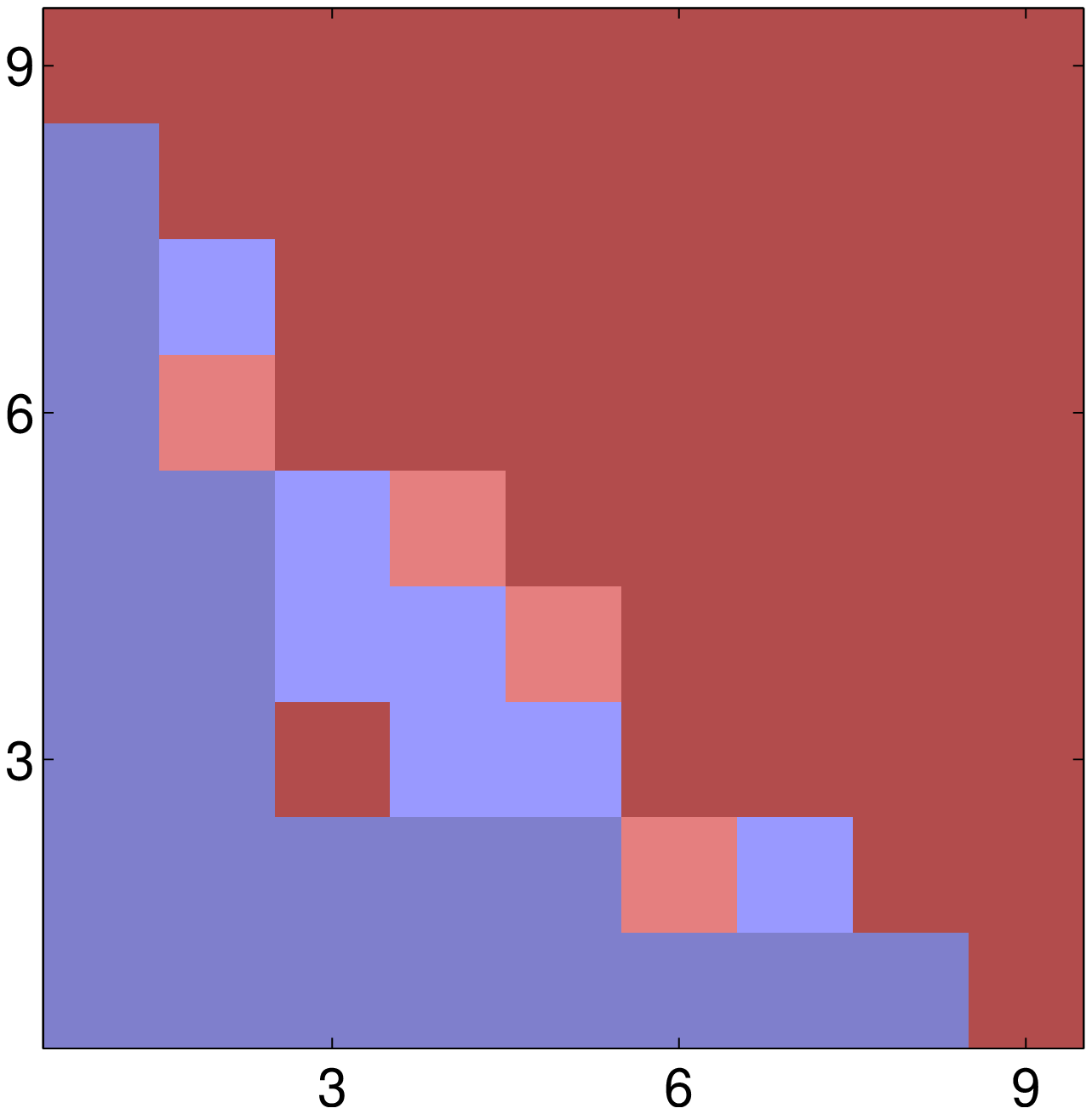} \
   \includegraphics[width=5cm]{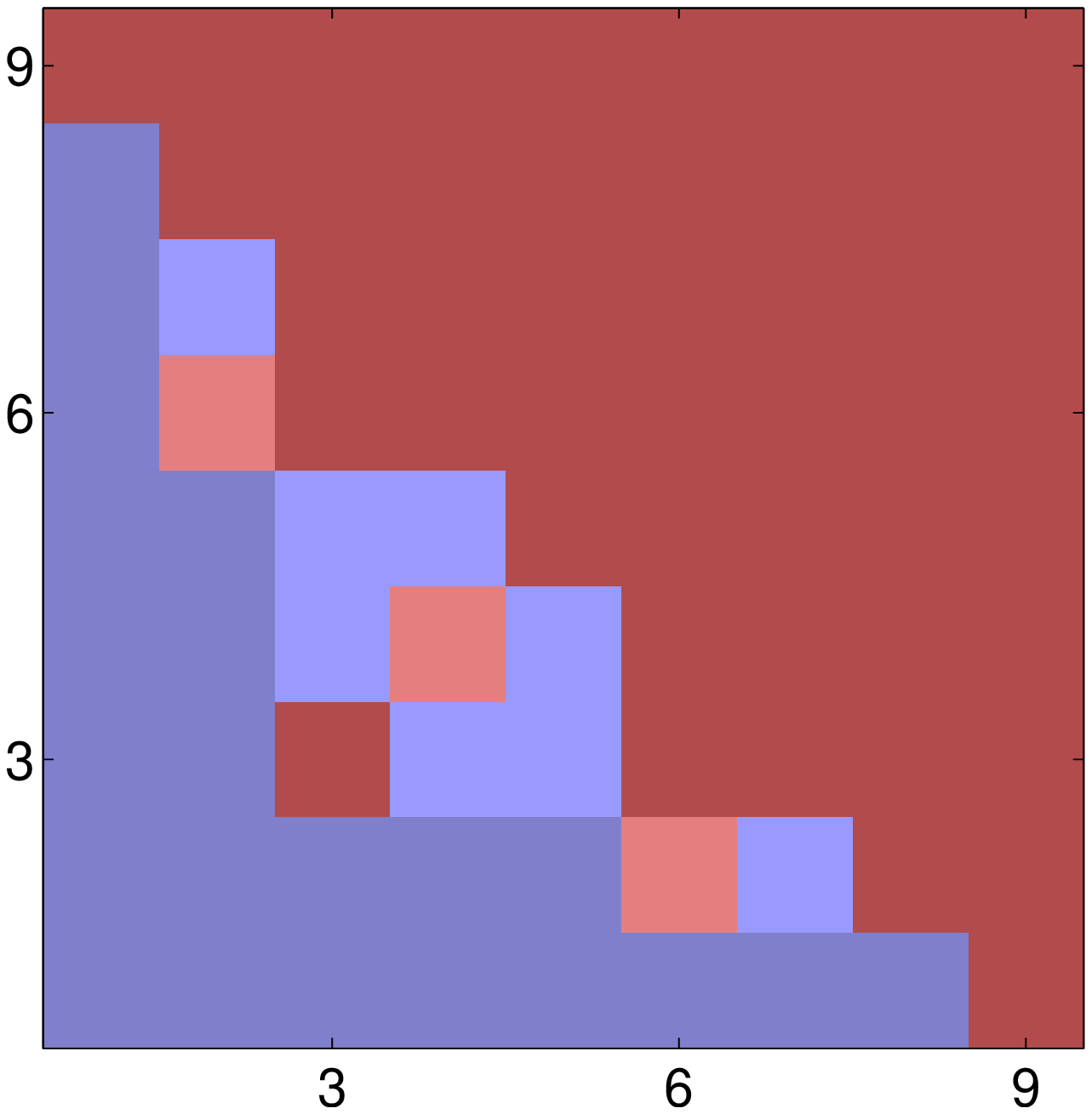}

\vspace{-4.7cm}

\hspace{1.2cm}\color{white}
  \hspace{1cm}  \begin{tabular}{p{4.65cm}p{4.65cm}p{4.65cm}}
      $\Z_2{\times}\Z_2{\times}\Z_2$ & $\Z_9$& $\Z_3{\times}\Z_3$
    \end{tabular}\color{black}

\vspace{4.2cm}

   \includegraphics[width=5cm]{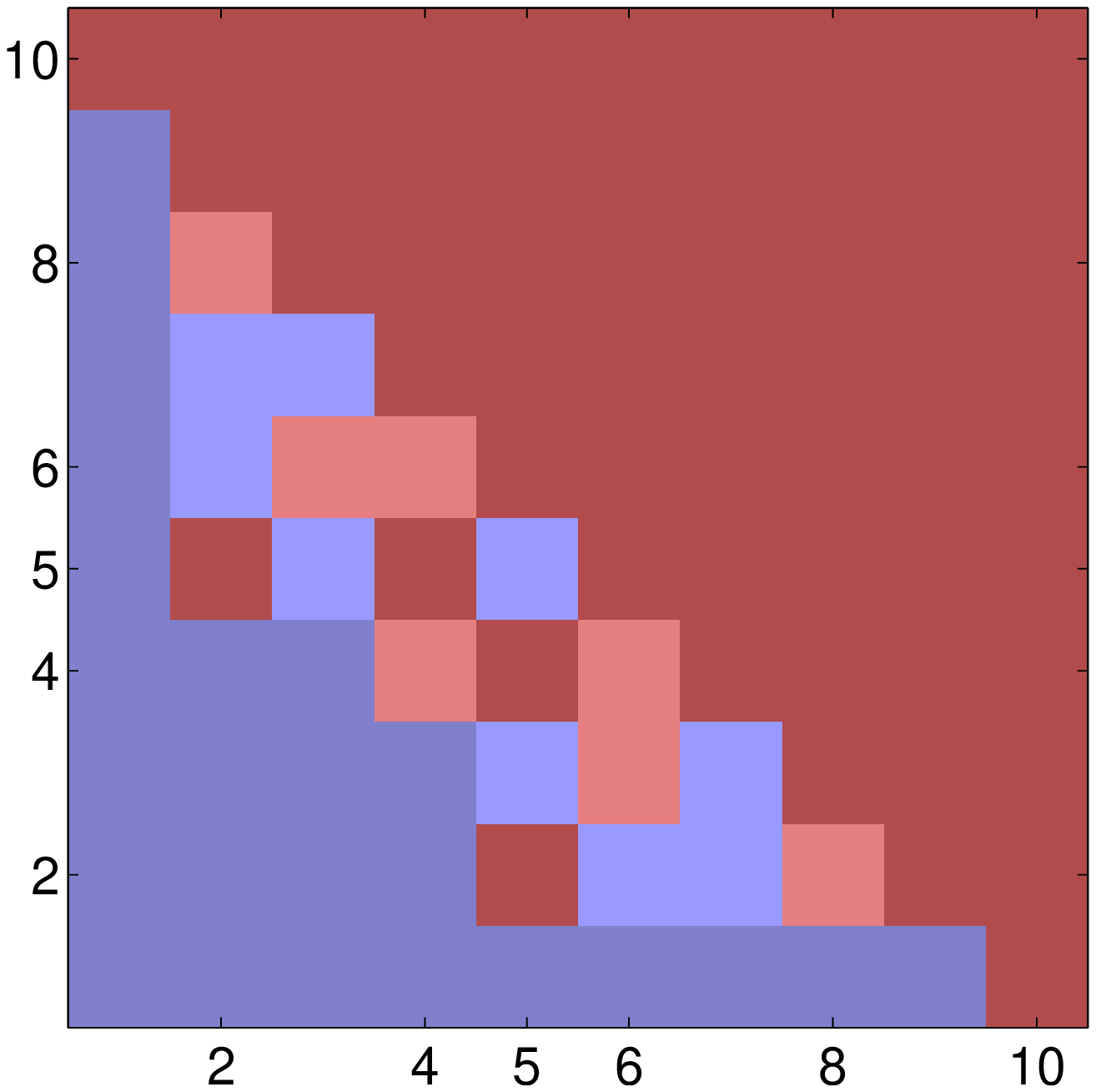} \
      \includegraphics[width=5cm]{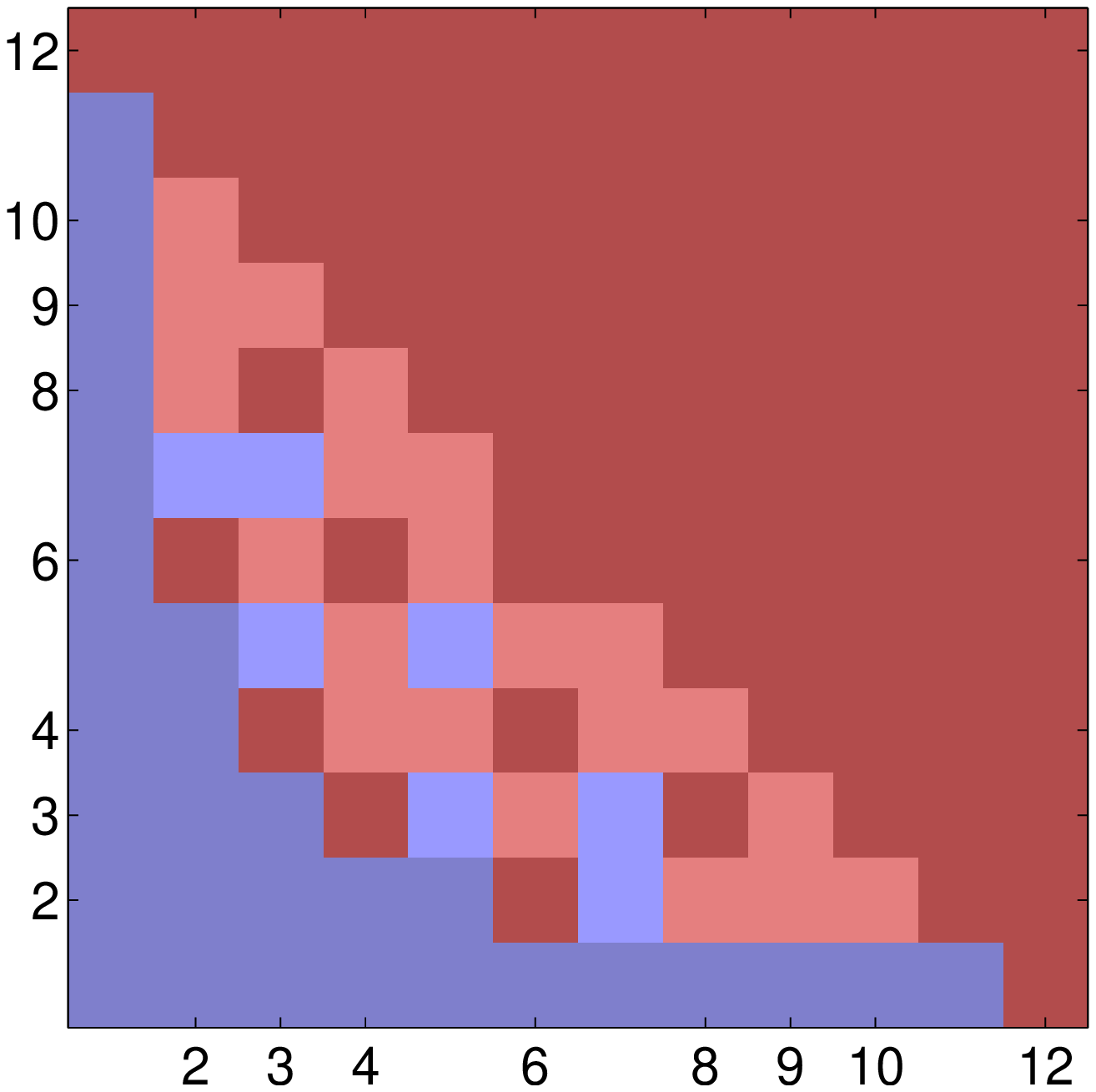} \
         \includegraphics[width=5cm]{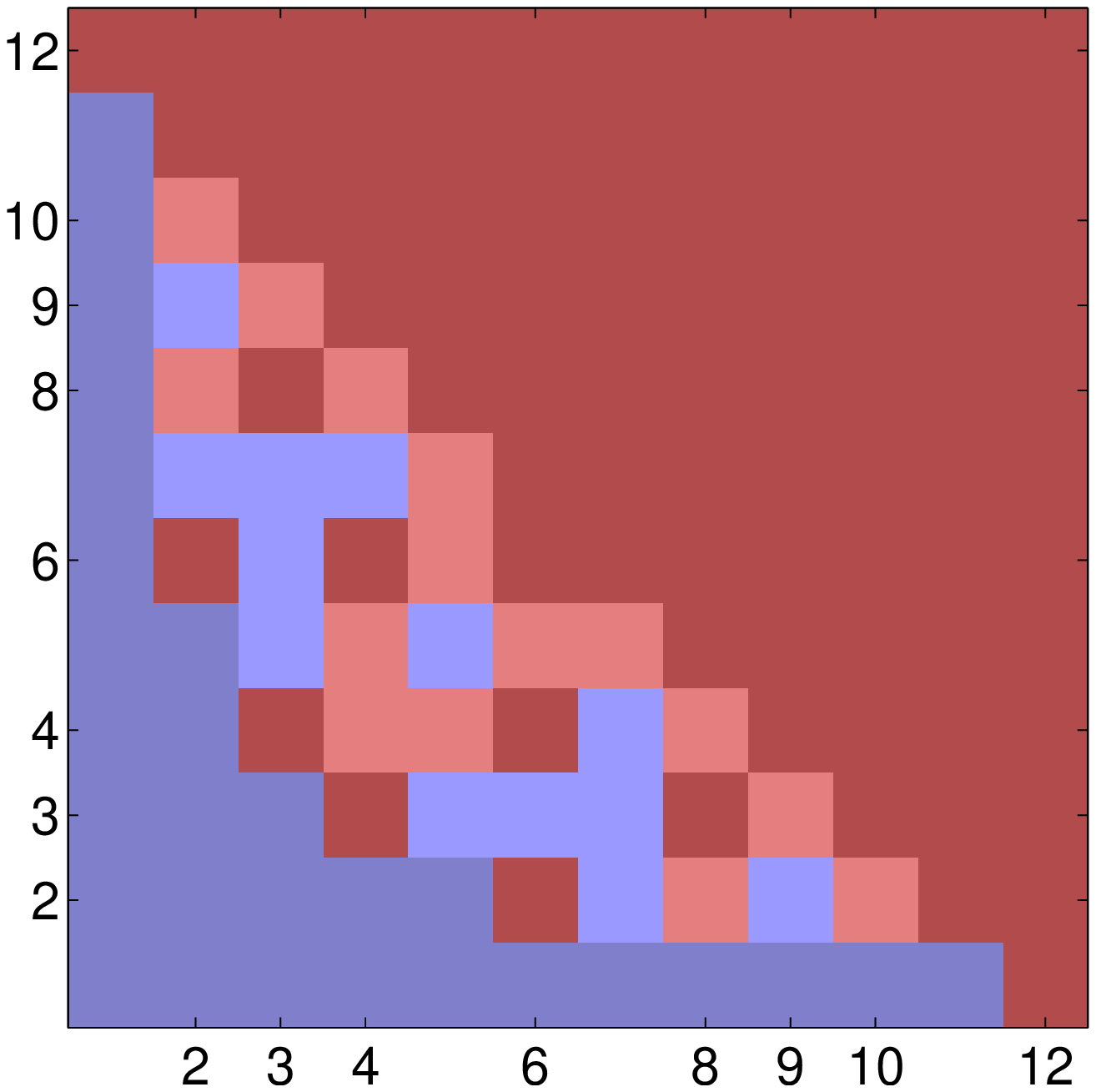}

 \vspace{-4.7cm}
 \hspace{.8cm}
 \color{white}
   \hspace{1cm} \begin{tabular}{p{4.65cm}p{4.65cm}p{4.65cm}}
      $\Z_{10}\cong\Z_2{\times}\Z_5$ & $\Z_{12}$ & $\Z_2{\times}\Z_6\cong\Z_2{\times}\Z_2{\times}\Z_3$
    \end{tabular}\color{black}
    \vspace{4.2cm} \label{figure:possiblepairsnoncyclic}
    }
\end{center}
\begin{caption}
 {
    The set $\big\{ (\|f\|_0, \|\fhat f\|_0),\ f\in\C^G{\setminus}\{0\}\big\}$ for all Abelian groups
    of non-prime order
    less than or equal to $12$.
    If $kl<|G|$, then no $f$ exists with $(k,l)=(\|f\|_0, \|\fhat f\|_0)$ by
    Theorem~\ref{theorem:classicalUncertainty}.
    If $|G|$ divides $kl$, or if $k+l\geq |G|{+}1$ then
    exists $f$ with $(k,l)=(\|f\|_0, \|\fhat f\|_0)$  by Proposition~\ref{proposition:DonohoStark} and
    Proposition~\ref{proposition:nominors}. The color code used is described
    in     Figure~\ref{figure:ColorMap}.
  }\label{figure:possiblepairsnoncyclic}
   \end{caption}

\end{figure}

\begin{figure}
\begin{center}
   \includegraphics[width=5cm]{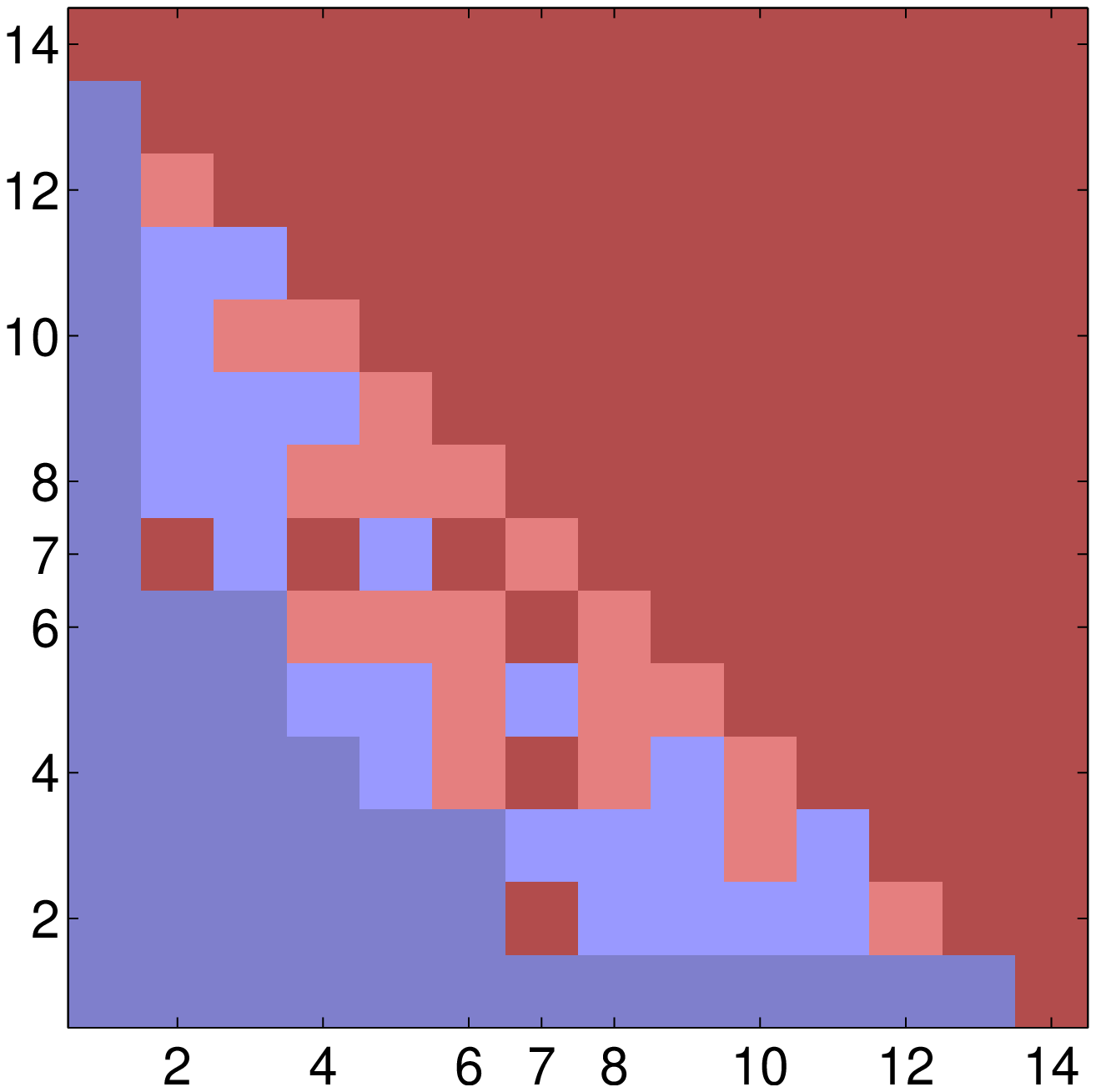} \
      \includegraphics[width=5cm]{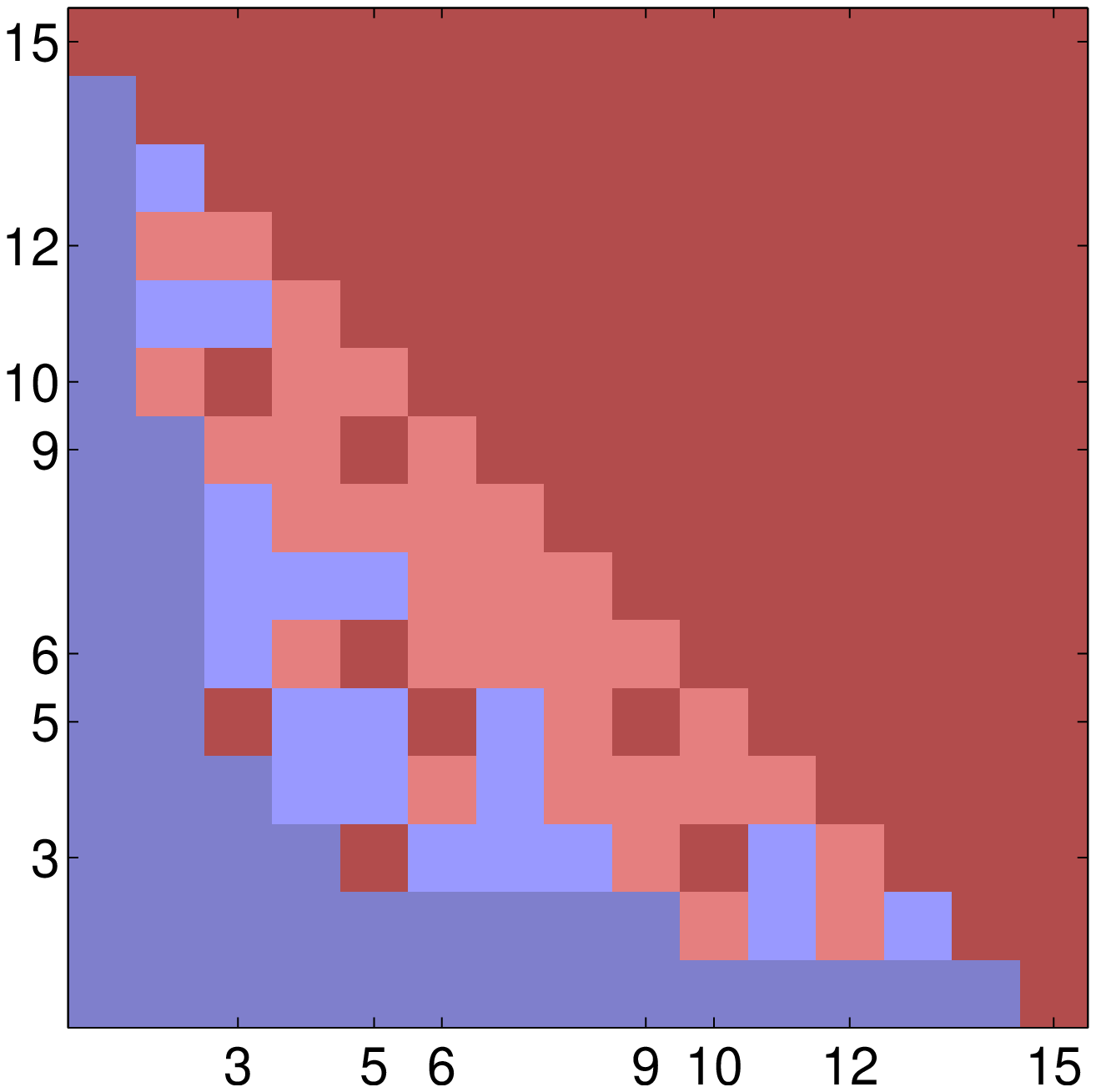} \
         \includegraphics[width=5cm]{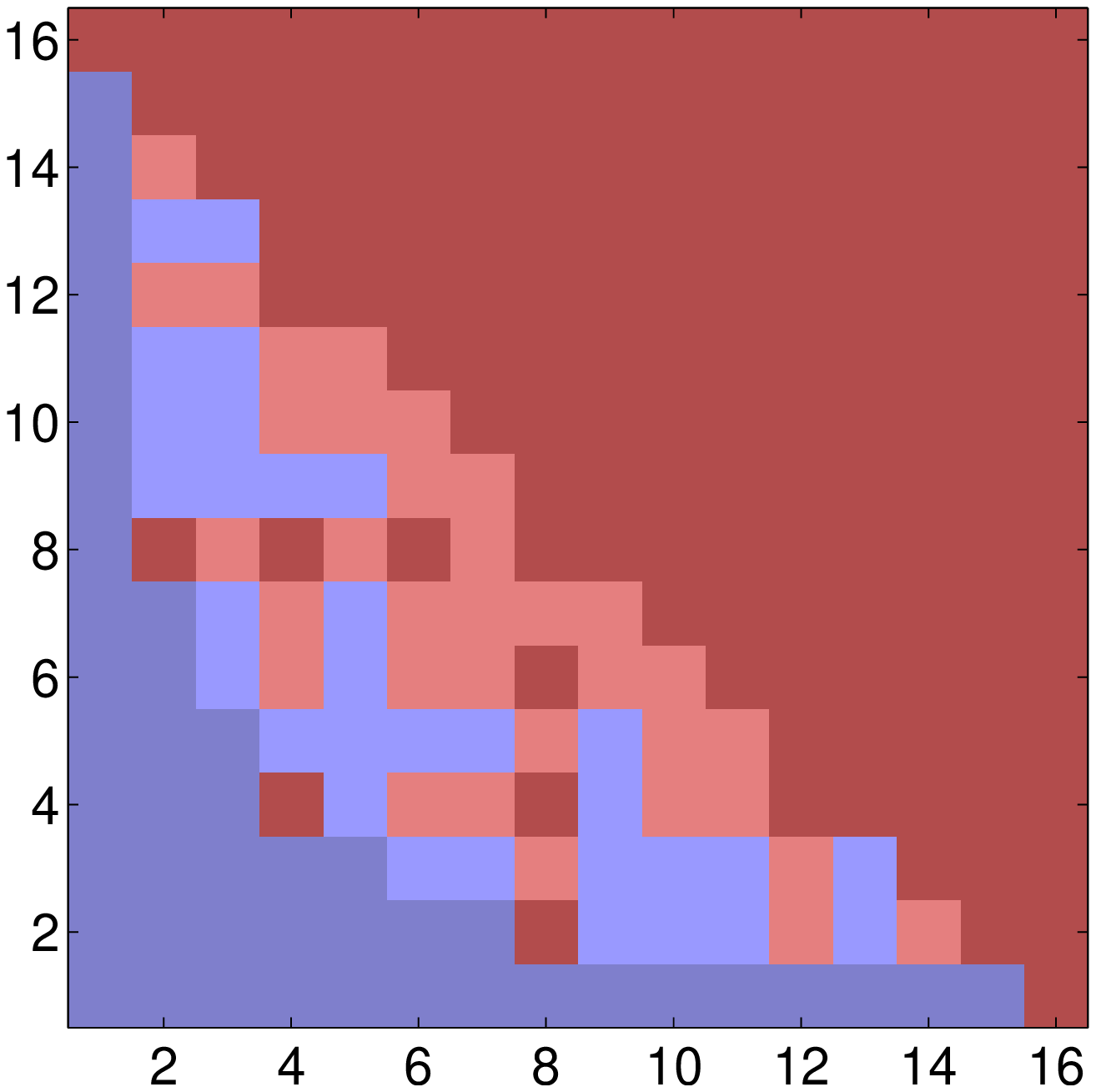}

    \vspace{-4.7cm}

\hspace{1.2cm}\color{white}
 \hspace{1cm}   \begin{tabular}{p{4.65cm}p{4.65cm}p{4.65cm}}
       $\Z_{14}\cong\Z_2{\times}\Z_7$ & $\Z_{15}\cong\Z_3{\times}\Z_5$ & $\Z_{16},\ \Z_2{\times}\Z_8$
    \end{tabular}\color{black}

\vspace{4.2cm}

   \includegraphics[width=5cm]{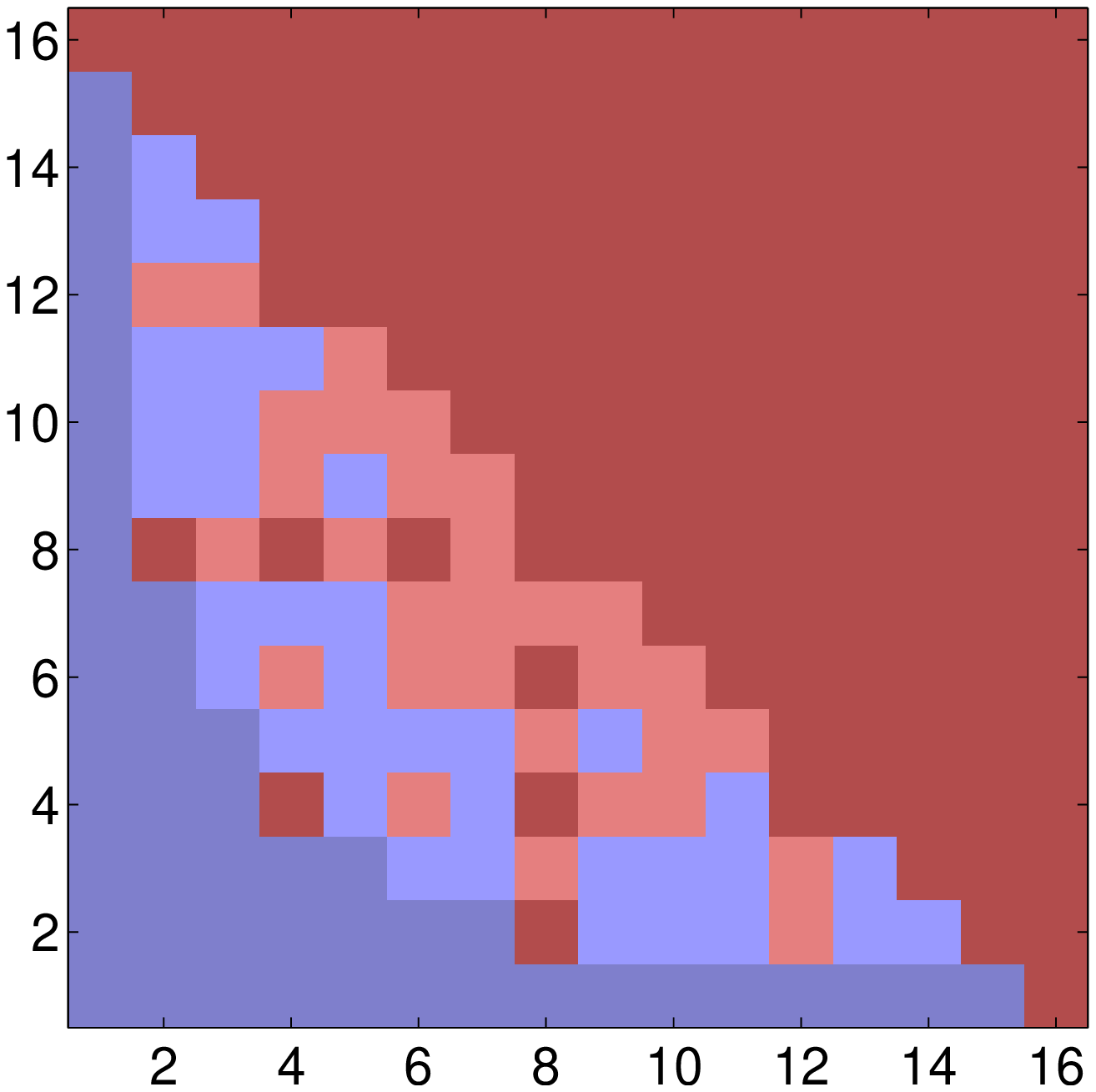} \
   \includegraphics[width=5cm]{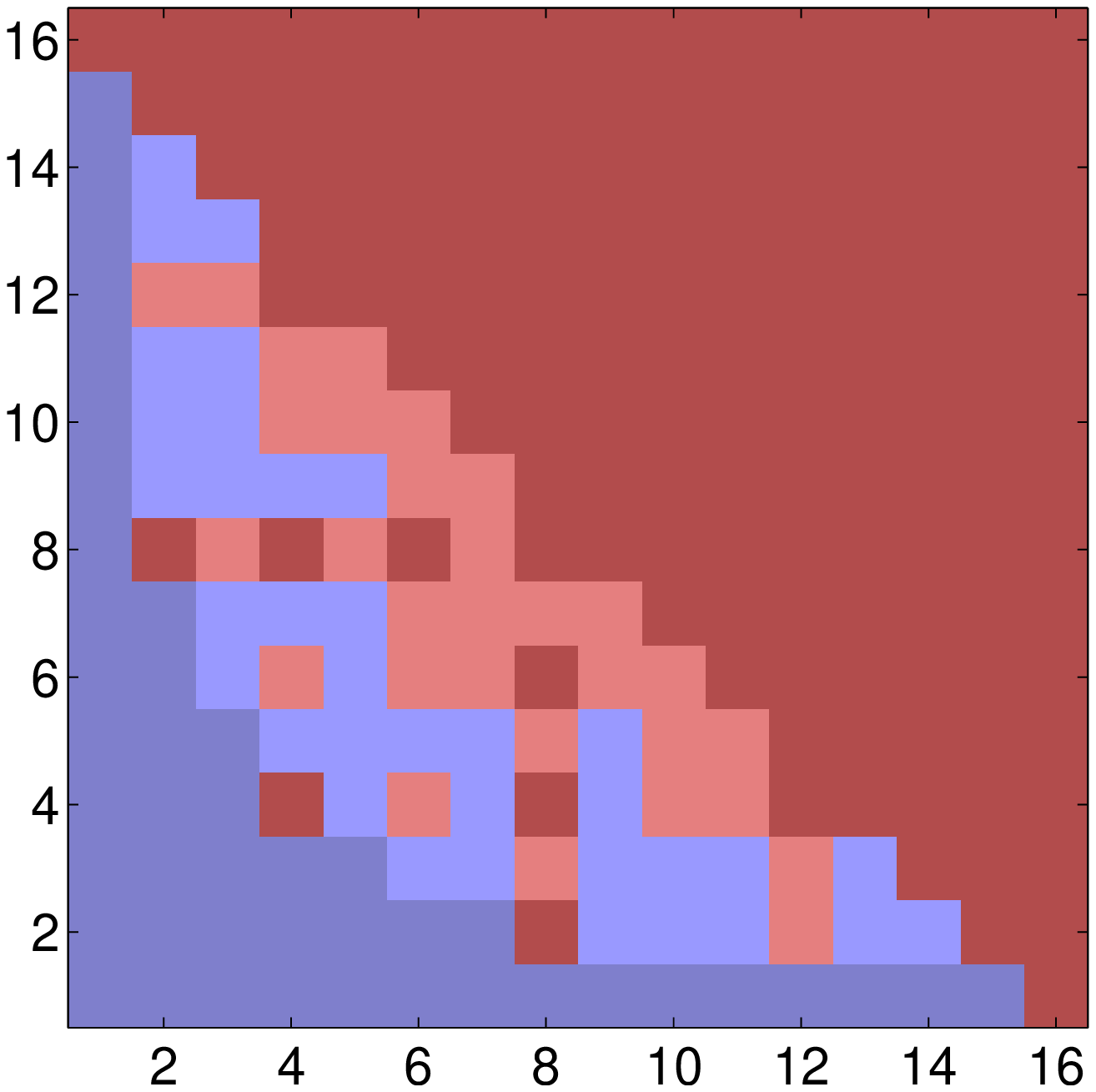} \
   \includegraphics[width=5cm]{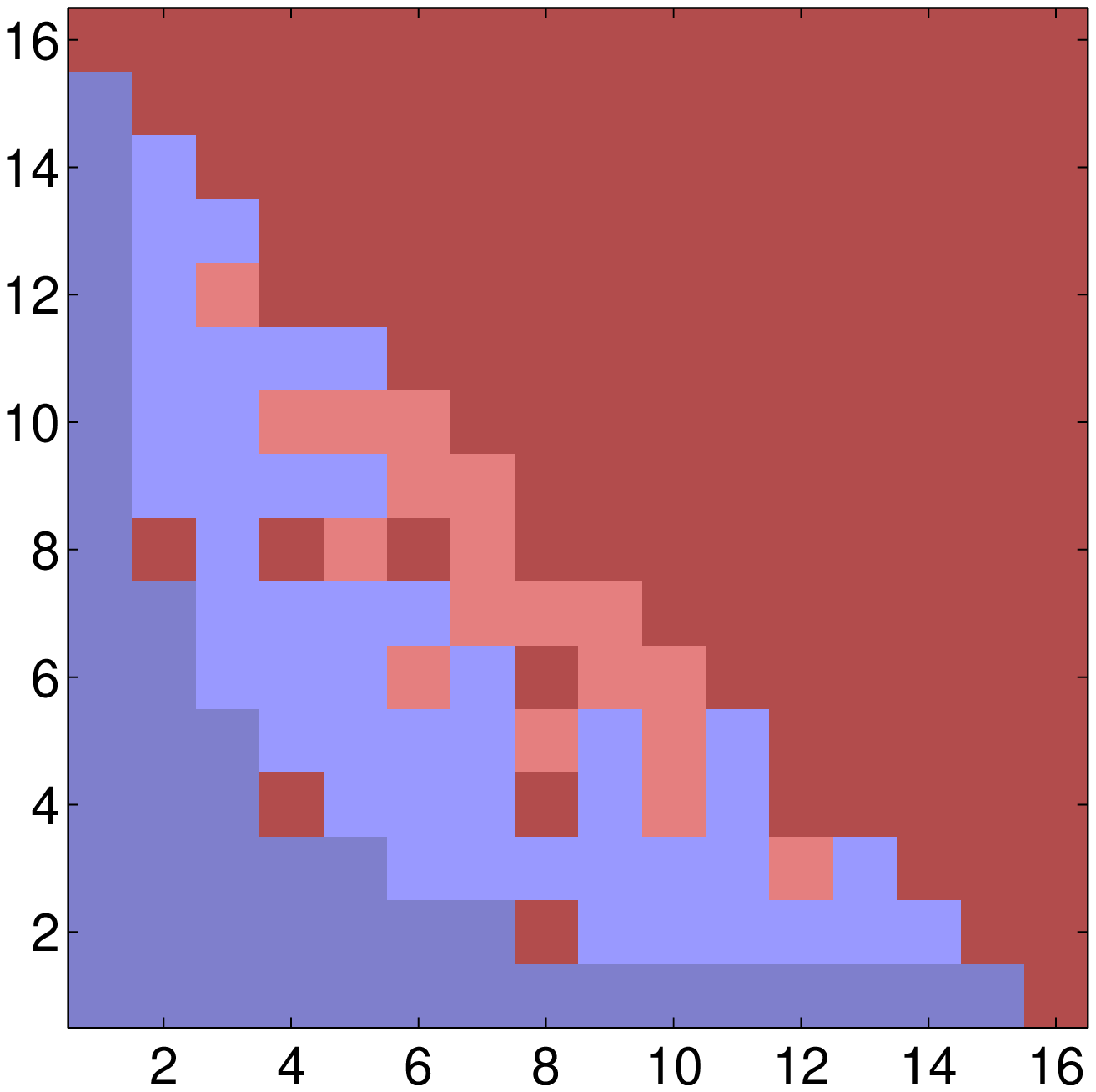}

\vspace{-4.7cm} \hspace{1.2cm}\color{white}
 \hspace{1cm}    \begin{tabular}{p{4.65cm}p{4.65cm}p{4.65cm}}
      $\Z_4{\times}\Z_4$ & $\Z_2{\times}\Z_2{\times}\Z_4$ & $\Z_2{\times}\Z_2{\times}\Z_2{\times}\Z_2$
    \end{tabular}\color{black}
    \vspace{4cm}

\end{center}

\begin{caption}{Same as Figure~\ref{figure:possiblepairsnoncyclic} for
Abelian groups of order 14, 15, and
16.}\label{figure:possiblepairsnoncyclic2}
\end{caption}
\end{figure}

The numerical results collected in
Figure~\ref{figure:possiblepairsnoncyclic} and
Figure~\ref{figure:possiblepairsnoncyclic2} are based on an idea in
\cite{FKLM05} and on Lemma~\ref{lemma:minorranks}. They show that
the set of all possible pairs $(\|f\|_0,\|\widehat{f}\|_0)$ is
nontrivial in general. The computations that lead to these results
are quite involved. For example, the computations showing that there
is no function (vector) on $\Z_{16}$ with six nonzero entries and
whose Fourier transform has nine nonzero entries include the
calculation of the singular values of $ \mysmallmatrixCENTER{1}{16\\
9}\mysmallmatrixCENTER{1}{16\\ 5}=49969920$ nine by six matrices.

In addition, we give all possible pairs
$(\|f\|_0,\|\widehat{f}\|_0)$ for the group $G=\Z_6$ and give a
partial result for the groups $G=\Z_{2p}$ for, $p\geq 5$ prime.
Their proofs are included in the appendix.

\begin{proposition}\label{proposition:Z6} For $1\leq k,l\leq 6$
exists $f\in \C^{\Z_6}$ with $\|f\|_0=k$ and $\|\widehat{f}\|_0=l$
if and only if $kl\geq 6$ and $(k,l)\neq (3,3)$.
\end{proposition}

The following result for $\Z_{2p}$, $p\geq 5$ prime,  shows that the
bound in Theorem~\ref{theorem:meshulam} is not sharp, a fact that
was  observed for the case $G=\Z_8$ in \cite{FKLM05}.


\begin{proposition}\label{proposition:FourierSetpq}
For $p\ge5$  prime there exists no $f \in \C^{\Z_{2p}}$ with
$\|f\|_0=3$ and $\|\widehat f\|_0=p{-}1$.
\end{proposition}

\section{Uncertainty principles for short--time Fourier transforms on finite Abelian groups}
\label{section:STFTuncertainty}

We now turn to discuss minimum support conditions on time-frequency
representations of elements in $\C^G$, in particular, for the
short--time Fourier transform of a function $f\in\C^G$ with respect
to a window $g\in \C^G$.

\begin{figure}[t]
\begin{center}
%
%
%
%
%
%
%
%
%
%
%
%
%
%
%
%
%
%
%
%
%
%
%
%
%
%
%
%
%

\includegraphics[width=8cm]{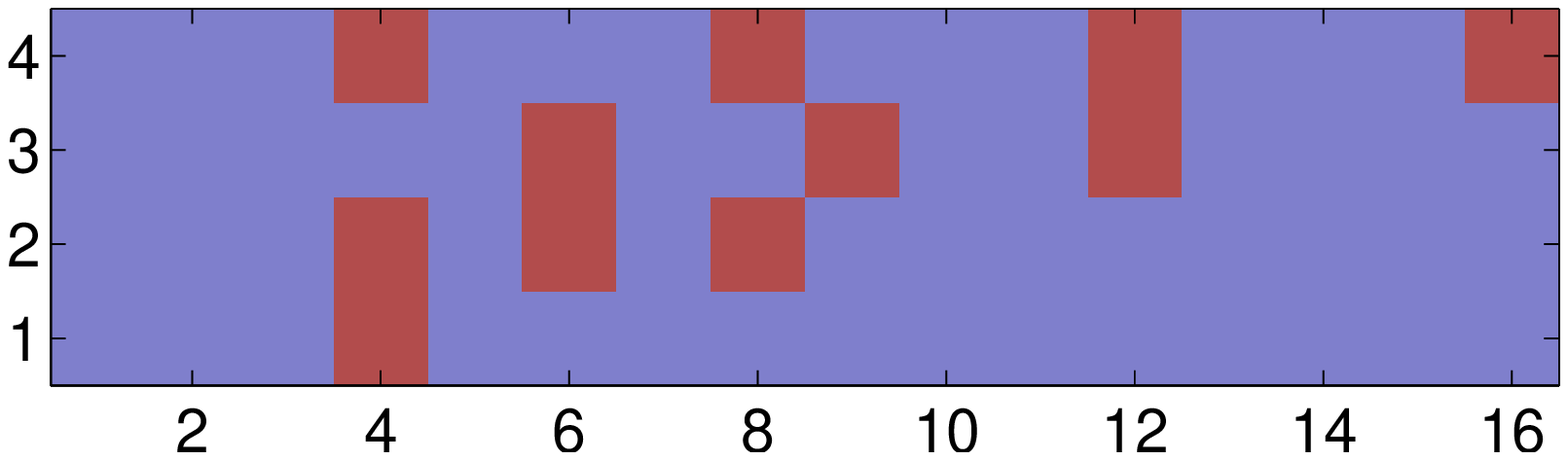}
     \caption{\color{black} For the Abelian group $G=\Z_4$ all
     possible pairs
      $(\|f\|_0, \|R f\|_0)$  are colored
     red, those pairs that are not achieved by some $f\in
     \C^{\Z_4}$ are colored blue in accordance with the color code given in Figure~\ref{figure:ColorMap}.}
     \color{black}\label{figure:tensor}
\end{center}
\end{figure}

The simplest joint time-frequency representation of $f$ is given by
the tensor product $f\otimes \fhat f$. Similarly, in electrical
engineering the so-called Rihaczek distribution, $R:G\times\widehat
G\longrightarrow \C$, which is given by
$Rf(x,\omega)=f(x)\overline{\fhat f(\omega)}\,\overline{\langle
\omega,x\rangle}$, is considered.
Theorem~\ref{theorem:classicalUncertainty} implies that
$\|Rf\|_0=\|f\otimes\widehat{f}\|_0=\|f\|_0\|\widehat{f}\|_0\geq
|G|$. Figure~\ref{figure:tensor} lists all possible pairs $(\|f\|_0,
\|Rf\|_0)$ for $f\in\C^{\Z_4}$.

Using the technique used to obtain
Theorem~\ref{theorem:classicalUncertainty}, we obtain the well-known
result
\begin{proposition}\label{propsition:uncertainty-for-Vgf}
$\displaystyle \|V_g f\|_0 \geq |G|$ for $f,g\in \C^G {\setminus}
\{0\}$ with equality for $f=g=\delta$.
\end{proposition}

\begin{proof}Clearly $\|V_{\delta}\delta\|_0 = |G|$. For  $f,g\in \C^G {\setminus} \{0\}$,
\begin{equation}
    |G| \ \|f\|_2^2 \ \|g\|_2^2 = \|V_g f\|_2^2 \leq \|V_g f\|_0\ \|V_g f\|_\infty^2
        \leq \|V_g f\|_0\ \|f\|_2^2\ \|g\|_2^2\notag
 \end{equation}
 and the result follows.
\end{proof}

We shall now seek lower bounds on $\|V_gf\|_0$ depending on
$\|f\|_0$, $\|\widehat f\|_0$, $\|g\|_0$, and $\|\widehat g\|_0$.

%



\begin{proposition}
  For $f,g\in \C^G{\setminus}\{0\}$, we have
  \begin{eqnarray}
    \|V_g f\|_0 \geq \max\{ \ \theta(G,\|g\|_0) \, \theta(G,\|\widehat{f}\|_0)\, ,\,
    \theta(G,\|f\|_0) \, \theta(G,\|\widehat{g}\|_0)\
    \}\label{equation:Vgf1general}\,,
  \end{eqnarray}
  and, therefore,
   \begin{eqnarray}
    \|V_g f\|_0 \geq \tfrac 1 2 \left( \, \theta(G,\|g\|_0) \, \theta(G,\|\widehat{f}\|_0)+\,
    \theta(G,\|f\|_0)
    \,
    \theta(G,\|\widehat{g}\|_0)\right)\label{equation:Vgf2general}\,,
  \end{eqnarray}
  and
  \begin{eqnarray}
    \|V_g f\|_0 \geq \sqrt{\, \theta(G,\|f\|_0)  \, \theta(G,\|\widehat{f}\|_0)\, \theta(G,\|g\|_0)
    \, \theta(G,\|\widehat{g}\|_0)}\label{equation:Vgf3general}\,.
  \end{eqnarray}
\end{proposition}
\begin{proof}
We shall prove $\displaystyle \|V_g f\|_0\geq
\theta(G,\|f\|_0)\theta(\widehat G, \|\widehat  g\|_0)$. Then
(\ref{equation:Vgf1general}) follows from   $\|V_g
f\|_0=\|V_{\widehat{g}} \widehat{f}\|_0$ and
$\theta(G,k)=\theta(\widehat{G},k)$ for any $k$, or, alternatively
from $\|V_gf\|_0=\|V_f  g\|_0$. Further,
  (\ref{equation:Vgf1general}) implies
  (\ref{equation:Vgf2general}) and (\ref{equation:Vgf3general})
  since the
  maximum of two positive numbers dominates their arithmetic and geometric
  means.

  To see (\ref{equation:Vgf1general}), observe first that the so-called symplectic Fourier transformation
  $\mathcal F_s = R \circ \mathcal F_{\fhat G}^{-1}\circ \mathcal F_G$,
  i.e.,
  the composition of a Fourier transformation $\mathcal F_G$ on $G$, an inverse Fourier
  transformation
  $\mathcal F_{\fhat G}^{-1}$ on $\widehat{G}$, and the axis transformation
  $R:F\mapsto F\circ \mysmallmatrix{2}{0 & 1 \\ 1& 0}$ obeys the same uncertainty principle
  as the Fourier transformation on the
  group $G{\times}\widehat G$.
  For $f,g\in \C^G$, we calculate
\begin{eqnarray*}
    \mathcal F_s V_g f(r,\rho)
        &=&     \sum_{x \in G}  \sum_{\xi \in \widehat G} V_g f(x, \xi)
                    \overline{\langle \rho, x\rangle} \langle \xi, r \rangle
        =       \sum_{x \in G} \sum_{\xi \in \widehat G}
                    \sum_{t \in G} f(t)\overline{ g (t-x)}\,\overline{\langle \xi,t
                    \rangle}
                    \ \overline{\langle \rho, x\rangle} \langle \xi, r \rangle\\
        &=&       \sum_{x \in G}
                    \sum_{t \in G} f(t)\overline{ g (t-x)}\,
                    \overline{\langle \rho, x\rangle} \sum_{\xi \in \widehat G}
                    \langle \xi,r-t \rangle
        = |G|\sum_{x \in G}
                     f(r)\overline{ g (r-x)}\,
                    \overline{\langle \rho, x\rangle} \\&=& |G|
                     \overline{\langle \rho, r\rangle} f(r)\overline{ \widehat g (\rho)}
\end{eqnarray*}
and note that $\supp \mathcal F_s V_g f = \supp f\times\supp
\widehat g$. Proposition~\ref{proposition:productset} implies that
$\displaystyle \|V_g f\|_0 = \| \mathcal F_s^{-1} \big( \mathcal F_s
V_g f\big)\|_0 \geq \theta(G,\|f\|_0)\theta(\widehat G, \|\widehat
g\|_0).$
\end{proof}

\begin{table}[t]
    \begin{center}
{\tiny
  \begin{tabular}{|r|r||r|r|r|r|r|r|r|r|r|r|r|r|r|r|r|r|r|r|r|r|r|r|r|r|r|}
\hline
 & &1&2&2&2&2&3&3&3&3&4&4&4&4&4&5&5&5&5&5&6&6&6&6&6&6\\\hline
 & &6&3&4&5&6&2&4&5&6&2&3&4&5&6&2&3&4&5&6&1&2&3&4&5&6\\\hline\hline
1&6&6&24&18&12&6&30&18&12&6&30&24&18&12&6&30&24&18&12&6&36&30&24&18&12&6\\\hline
2&3&24&20&20&20&20&25&16&16&16&25&20&15&12&12&25&20&15&10&8&30&25&20&15&10&5\\\hline
2&4&18&20&15&15&15&25&15&12&12&25&20&15&10&9&25&20&15&10&6&30&25&20&15&10&5\\\hline
2&5&12&20&15&10&10&25&15&10&8&25&20&15&10&6&25&20&15&10&5&30&25&20&15&10&5\\\hline
2&6&6&20&15&10&5&25&15&10&5&25&20&15&10&5&25&20&15&10&5&30&25&20&15&10&5\\\hline
3&4&18&16&15&15&15&20&12&12&12&20&16&12&9&9&20&16&12&8&6&24&20&16&12&8&4\\\hline
3&5&12&16&12&10&10&20&12&8&8&20&16&12&8&6&20&16&12&8&4&24&20&16&12&8&4\\\hline
3&6&6&16&12&8&5&20&12&8&4&20&16&12&8&4&20&16&12&8&4&24&20&16&12&8&4\\\hline
4&4&18&15&15&15&15&15&12&12&12&15&12&9&9&9&15&12&9&6&6&18&15&12&9&6&3\\\hline
4&5&12&12&10&10&10&15&9&8&8&15&12&9&6&6&15&12&9&6&4&18&15&12&9&6&3\\\hline
4&6&6&12&9&6&5&15&9&6&4&15&12&9&6&3&15&12&9&6&3&18&15&12&9&6&3\\\hline
5&5&12&10&10&10&10&10&8&8&8&10&8&6&6&6&10&8&6&4&4&12&10&8&6&4&2\\\hline
5&6&6&8&6&5&5&10&6&4&4&10&8&6&4&3&10&8&6&4&2&12&10&8&6&4&2\\\hline
6&6&6&5&5&5&5&5&4&4&4&5&4&3&3&3&5&4&3&2&2&6&5&4&3&2&1\\\hline
\end{tabular}
}
\begin{caption}{Numerical representation of \eqref{equation:Vgf1general} for  $G=\Z_6$.
Rows represent possible pairs $(\|f\|_0,\|\fhat f\|_0)$, columns
possible pairs $(\|g\|_0,\|\fhat g\|_0)$, and the table entries give
the lower bound on $\|V_g\|_0$. }\label{table:LowerBoundVgfZ6-v2}
\end{caption}
\end{center}
\end{table}

For $G=\Z_6$, we list in Table~\ref{table:LowerBoundVgfZ6-v2} the
lower bounds on $\|V_g f\|_0$ given by (\ref{equation:Vgf1general})
for  different values of  $\|f\|_0$, $\|\widehat f\|_0$, $\|g\|_0$
and $\|\widehat g\|_0$.


\begin{corollary}\label{corollary:primeorder} For $f,g \in \C^{\Z_p}{\setminus}\{0\}$, $p$ prime,
  \begin{eqnarray}
    \|V_g f\|_0 \geq \max\{ \ (p{+}1{-}\|g\|_0)  (p{+}1{-}\|\widehat{f}\|_0)\, ,\,
    (p{+}1{-}\|f\|_0)  (p{+}1{-}\|\widehat{g}\|_0)\
    \}\notag
  \end{eqnarray}
  and \quad $\displaystyle
    \|V_g f\|_0 \geq (p{+}1)^2    - \tfrac 1 2
                         (p{+}1)(\|f\|_0 +\|\widehat{f} \|_0 + \|g\|_0+\|\widehat{g}\|_0 )
                            + \tfrac 1 2\left(\|\widehat{f} \|_0\|g\|_0+\|f\|_0 \|\widehat{g}\|_0\right)
                            \notag
 $.
\end{corollary}

Now,  we give an improvement to the lower bound on $\|V_gf\|_0$ that
is given in Corollary~\ref{corollary:primeorder}.


\begin{proposition}
For $f,g \in \C^{\Z_p}{\setminus}\{0\}$, $p$ prime,
\[ \|V_g f\|_0 \geq \left\{ \begin{array}{ll}
         p(p{+}1)- \|f\|_0\|g\|_0 & \quad  \mbox{if $\|f\|_0+\|g\|_0>p$};\\
        p(p{+}1) - (p {+}1{-} \|f\|_0)(p{+}1{-} \|g\|_0) &\quad  \mbox{if $\|f\|_0+\|g\|_0\leq p$}.\end{array}
        \right. \]
\label{proposition:primevgf}
\end{proposition}

\vspace{-.3cm}
\begin{proof}
Note that for all $x\in G$, $V_g f (x,\cdot) =  \langle f, \pi(x
,\cdot) g\rangle$ represents the Fourier transform of a vector of
the form $f\,T_x \bar g$, i.e.,
$$
   V_gf(x,\xi)= \langle f, \pi(x,\xi) g \rangle
        =\displaystyle\sum_{y} f(y)\overline{g(y{-}x)} \overline{\langle \xi,x\rangle}
        =\widehat{f\,T_x \bar g}(\xi)\, \quad x\in G,\xi\in \widehat G \,.
$$

\vindex \noindent As long as $f\,T_x \bar g \neq 0$,
Theorem~\ref{theorem:tao} applies and so $\|f\,T_x \bar
g\|_0+\|\widehat{f\,T_x \bar g}\|_0\ge p + 1$. For $K:=\{x: f\,T_x
\bar g\neq 0\}$ we get
$$
    \|V_g f\|_0{=} \sum_{x\in K} \|\widehat{f\,T_x \bar g}\|_0 \,
        \ge\, |K|(p {+} 1) {-} \sum_x \|f\,T_x \bar g\|_0
        = |K|(p {+} 1) {-} \|f\|_0\|g\|_0,
$$
where $\displaystyle\sum_x\|f\,T_x \bar g\|_0= \|f\|_0\|g\|_0$
follows from a simple counting argument.

We shall now estimate $|K|$ using the Cauchy-Davenport inequality,
which states that for non-empty subsets $A$ and $B$ of
$\mathbb{Z}_p$, $p$ prime, $|A{+}B|\geq \min(|A|{+}|B|{-}1, p)$,
where $A+B=\{a{+}b:a\in A, b\in B\}$ \cite{K05}. Now $K = \{x:
f\,T_x \bar g \neq 0\}=\{x: \{(\supp \bar g) {+}x\} \cap \supp
f\neq\varnothing \} = \supp f {-}\supp \bar g$.
We set $A=\supp f, B=\supp \bar g$, and obtain $|K|=|\supp f{-}\supp
\bar g|\ge \min(\|f\|_0{+}\|g\|_0{-}1,p)$.

If $\|f\|_0{+}\|g\|_0 \ge p{+}1$, then $|K| =p$ and, hence, $\|V_g
f\|_0 \ge p(p{+}1) {-}\|f\|_0\|g\|_0$. If $\|f\|_0{+}\|g\|_0 \le p$,
then $|K| \ge \|f\|_0{+}\|g\|_0{-}1$ and so
$$
    \|V_g f\|_0 \geq (\|f\|_0{+}\|g\|_0{-}1)(p{+}1) - \|f\|_0\|g\|_0
        = p(p{+}1) - (p {+} 1 {-} \|f\|_0)(p{+} 1 {-} \|g\|_0)\,.
$$

\vspace{-.8cm}
\end{proof}

\begin{figure}[th]
\begin{center}
%
\includegraphics[height=3.6cm]{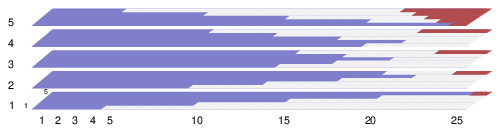}
\end{center}
\vspace{-3.2cm} \hspace{1.2cm}{ $\|f\|_0$}

\vspace{1.4cm} \hspace{1.1cm}$\|g\|_0$

\vspace{-.2cm}\hspace{8cm}$\|V_g f\|_0$
 \caption{\color{black} The set $\big\{
(\|f\|_0,\|g\|_0,\|V_g f\|_0),\ f,g\in\C^{G}{\setminus}\{0\}\big\}$
for $G=\Z_5$. The color coding is based on
Figure~\ref{figure:ColorMap} and justified by
Proposition~\ref{proposition:primevgf} and
Theorem~\ref{theorem:p2+1}.
    \label{figure:LowerBoundsVgfZ5}
               }

\end{figure}

The lower bound on $\|V_gf\|_0$ given in
Proposition~\ref{proposition:primevgf} is illustrated for $G=\Z_5$
in Table~\ref{figure:LowerBoundsVgfZ5}. To establish results similar
to Proposition~\ref{proposition:Z6} for the short--time Fourier
transformations for a given group $G$ is quite tedious since it
requires to check all combinations of $\|f\|_0$ and $\| g\|_0$. For
the case $G=\Z_3$, however, we have assembled all possible and
impossible combinations in Figure~\ref{figure:vgf}. A derivation of
the entries can be found in the appendix.

\begin{figure}[th]
\begin{center}
\includegraphics[height=3.5cm]{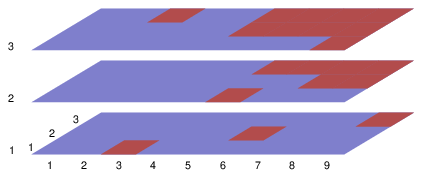}

\end{center}
\vspace{-3.2cm} \hspace{2.9cm}{ $\|f\|_0$}

\vspace{1.4cm} \hspace{2.8cm}$\|g\|_0$

\vspace{-.2cm}\hspace{8cm}$\|V_g f\|_0$
     \caption{\color{black} Same as Figure~\ref{figure:LowerBoundsVgfZ5} for $G=\Z_3$.  }\color{black}\label{figure:vgf}

\end{figure}

\subsection{Groups of prime order}\label{section:taobiroSTFT}

In the following, we shall fix the window $g$ and vary only the
analyzed function $f$. The main result in this section is
\begin{theorem}\label{theorem:p2+1}
  There exists $g \in \C^{\Z_p}$, $p$ prime,  such that for all $f\in
  \C^{\Z_p}$
  \begin{eqnarray}
  \displaystyle \|f\|_0 + \|V_g f\|_0\geq p^2{+}1 .\label{equation:p2+1}
  \end{eqnarray} Moreover, for $1 \le  k \le p$ and $1 \le l \le p^2$ with
  $ k+l\ge p^2{+}1$ there exists $f$ with $\|f\|_0=k$ and $\|V_g
  f\|_0=l$.
\end{theorem}

\begin{figure}[th]
\begin{center}

\includegraphics[height=1.7cm]{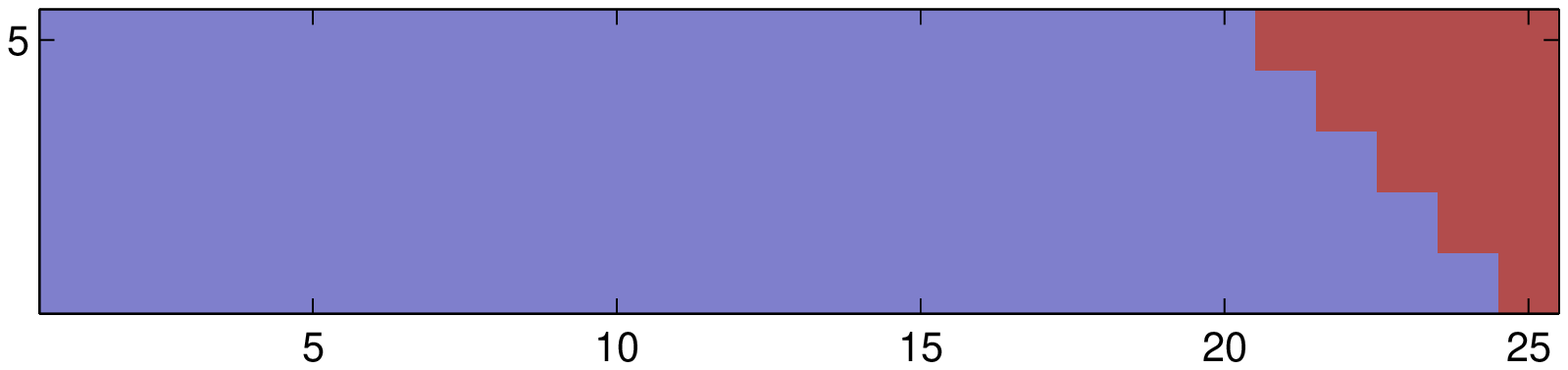}
\includegraphics[height=1.8cm]{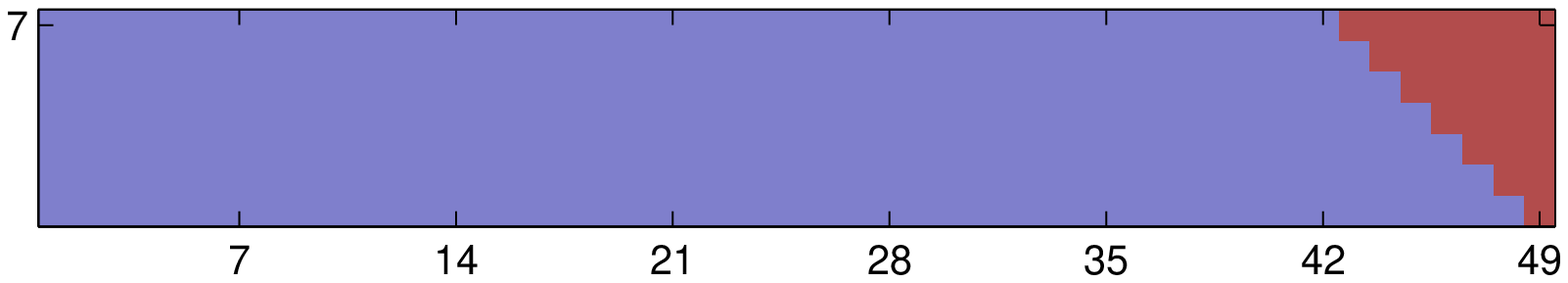}

\color{black}
\end{center}
\vspace{-2cm}\hspace{.7cm}$\Z_5$ \hspace{6.7cm} $\Z_{7}$
\color{black}

    \vspace{1.2cm}
\begin{caption}
  {The set $S_g=\big\{ (\|f\|_0,\|V_g f\|_0),\ f\in\C^G {\setminus} \{ 0 \}\big\}$ for
  appropriately chosen
  $g\in \C^G{\setminus}\{0\}$
   and
$G=\Z_5$ or $G=\Z_7$. The color coding is based on
Figure~\ref{figure:ColorMap} and justified by
    Theorem~\ref{theorem:p2+1}.}\label{figure:possiblepairsnoncyclicVgprime}
   \end{caption}
\end{figure}

We picture this result for $G=\Z_5$ and $G=\Z_7$ in
Figure~\ref{figure:possiblepairsnoncyclicVgprime}. Note that
Theorem~\ref{theorem:p2+1} follows from
Proposition~\ref{proposition:nominors} together with Theorem 4 from
\cite{LPW05} which we state as

\begin{theorem} \label{theorem:LaPfaWa} For almost every $g\in
\C^{\Z_p}$, $p$ prime, we have that every minor of $A_{\Z_p,g}$ is
nonzero.
\end{theorem}

{\it Outline of a proof of Theorem~\ref{theorem:LaPfaWa}.} It
suffices to show that each square submatrix  $(A_{\Z_p,g})_{A,B}$
has determinant nonzero for almost every $g$.

To this end, choose $A\subseteq G$ and $B\subseteq G{\times}\widehat
G$ with $|A|=|B|$ and set $P_{A,B}(z)= \det (A_{\Z_p,z})_{A,B}$,
$z=(z_0, z_1,\ldots, z_{p{-}1})$.  To show that $P_{A,B}\neq 0$, we
shall locate a term in the polynomial in standard form which has a
nonzero coefficient. To construct this term, we determine first the
maximal possible exponent of $z_0$ in one of the terms of $P$ that
are not trivially zero. Next, we determine the maximal exponent that
$z_1$ can have in a monomial where the maximal exponent of $z_0$ is
attained and so on.

Using generalized Vandermonde determinants, it can then be shown
that the coefficient of this ``maximal" term within $P_{A,B}$ can be
expressed as a product of different minors of the discrete Fourier
matrix $W_{\Z_p}$. For $p$ prime, all these minors are nonzero, so
the polynomial $P$ has a nonzero coefficient for this ``maximal
term", hence is not identically $0$, and nonzero almost everywhere.
We have $P=\prod\limits_{A,B:\,|B|=|A|}P_{A,B}\not\equiv 0$, which
implies that for $g\notin Z_P=\{z:\ P(z)=0 \}$, every minor of
$A_{\Z_p,g}$ is nonzero. Clearly, since $P\not\equiv 0$, $Z_P$ has
Lebesgue measure 0. \hfill $\Box$

Clearly, this proof of Theorem~\ref{theorem:LaPfaWa} is based on
Chebotarev's Theorem. Also, Chebotarev's Theorem and therefore
Theorem~\ref{theorem:tao} can be obtained as a corollary to
Theorem~\ref{theorem:LaPfaWa} as shown in the Appendix.

It is easy to see that if $g\in \C^{\Z_p}$ satisfies
(\ref{equation:p2+1}) then $\|g\|_0=\|\fhat g\|_0=p$, i.e.,
$g(x)\neq 0$ for all $x\in G$ and $\widehat g(\xi)\neq 0$ for all
$\xi \in \widehat G$ \cite{LPW05}. Further, we have

\begin{proposition}\label{unitg}
There exists a $g\in\C^{\Z_p}$, $p$ prime, with $|g(x)|=1$ for all
$x\in G$ and which satisfies the conclusions of
Theorem~\ref{theorem:p2+1}.
\end{proposition}
\begin{proof}
Theorem~\ref{theorem:LaPfaWa} implies that all minors of $A_{G,g}$
are nonzero polynomials in the polynomial ring $\C[z_0, ...
,z_{n-1}]$. Let $P$ be the product of all these minor polynomials,
which, by assumption, is nonzero. We have to show that $P(g)\neq 0$
for some $g\in\C^{\Z_p}$ with $|g(x)|=1$ for all $x\in G$.

This follows since the only polynomial $P$ with $P(g)=0$ whenever
$|g(x)|=1$ for all $x\in G$ is trivial, i.e., $P\equiv 0$, which we
show below using induction over the number of variables $n$.

The case $n=1$ follows since any nonzero polynomial in one variable
has only finitely many zeros, i.e., only $P\equiv 0$ vanishes for
all $z\in S^1=\{z:\ |z|=1\}$. Next, we consider a polynomial $P$ of
$n$ variables which we regard as a polynomial in $z_{n{-}1}$ with
coefficients in the polynomial ring $\C[z_0, ... ,z_{n{-}2}]$, i.e.,
\[
P(z_{n{-}1})=Q_m(z_0,...,z_{{n{-}2}}) z_{n{-}1}^m +
Q_{m{-}1}(z_0,...,z_{n{-}2}) z_{n{-}1}^{m{-}1}+\cdots +
Q_0(z_0,...,z_{n{-}2})
\]
For any fixed $(c_0,\ldots,c_{n{-}2})\in (S^1)^{n-1}$ we have
\[
Q_m(c_0,...,c_{{n{-}2}}) z_{n{-}1}^m + Q_{m{-}1}(c_0,...,c_{n{-}2})
z_{n{-}1}^{m{-}1}+\cdots + Q_0(c_0,...,c_{n{-}2})=0
\]
for all $z_{n-1}\in S^1$, hence, all its coefficients
$Q_k(c_0,...,c_{n-2})$, $k=0,\ldots, m$ vanish. In other words, we
have that $Q_{k}\in \C[z_0, ... ,z_{n-2}]$, $k=0,\ldots, m$ vanish
on $(S^1)^{n-1}$, which, by induction hypothesis, implies that all
$Q_k\equiv 0$ and therefore $P\equiv 0$.\end{proof}

Table~\ref{table:ranksOfMinorsVg} together with
Lemma~\ref{lemma:minorranks} show that the condition ``$G=\Z_p$ with
$p$ prime" is necessary for the existence of $g\in\C^G$ satisfying
(\ref{equation:p2+1}).

\begin{proposition}
If $|G|$ is not prime, then $A_{G,g}$ has zero minors for all $g \in
\C^G$.
\end{proposition}
\begin{proof}
Let $|G|=k\cdot m$, $k,m\neq 1$. We consider only $G=\Z_{km}$, the
general case follows since the Fourier matrix $W_G$ for any
non-cyclic $G$ is a Kronecker product of Fourier matrices of cyclic
groups.

For a primitive $|G|$-th root of unity $\omega$, we have
$(\omega^{k})^m=\omega^{|G|}=1$, so the discrete Fourier matrix
$W_G$ has a $1$ in
 its $(k, m)$-entry. Now the matrix given by the first $|G|$ columns of $A_{G,g}$ results
from $W_G$ by multiplying the $i$-th row by $c_i$. So the minor
given by the columns $0$ and $k$ and the rows $0$ and $m$ of $A$ is
$\det \mysmallmatrixLEFT{2}{ c_{0\,}& c_{0\,}\\ c_m& c_m}=0$. Hence
$A_{G,g}$ has a zero minor.
\end{proof}


%
%
%
%
%

%

\subsection{Groups of non-prime order}\label{section:meshulamSTFT}

Recall Proposition~\ref{propsition:uncertainty-for-Vgf}, namely, the
fact that for any $G$ the estimates $|G|\leq \|V_g f\|_0\leq |G|^2$
are sharp. In other words, for all $G$ and $0<k\leq |G|$ we have
$$
\min_{g\in\C^G{\setminus}\{0\}}\, \min \big\{\|V_g f\|_0:\ f\in \C^G\text{ and } 0<  \|f\|_0\leq k\big\}
    = |G|\, ,
$$
and
$$
\max_{g\in\C^G{\setminus}\{0\}}\, \max \big\{ \|V_g f\|_0: \
f\in \C^G\text{ and } 0<  \|f\|_0\leq k\big\} = |G|^2 \, .
$$

Certainly, $\|V_g f\|_0 = |G|$ is a rare event.  In fact, it is
reasonable to assume that  $\|V_g f\|_0 = |G|^2$ for almost every
pair $(f,g)$. We shall now address the question whether for an
appropriately chosen window $g$, we can achieve $\|V_g f\|_0\geq l$
for some $|G|<l\leq |G|^2$.

\begin{table}[th]
\begin{minipage}[t]{7.5cm}
{\footnotesize        \begin{tabular}[t]{|c||c|c|c|c|c|}
 \hline
            & 1 & 2 &3&4 &5   \\
          \hline
          \hline
           1& 125 &0  & 0 &  0&    0 \\
          \hline
           2&  0& 3000 & 0 &0  &   0 \\
          \hline
           3& 0 &  0&  23000& 0 &   0 \\
          \hline
           4& 0 & 0 & 0 &  63250 &  0   \\
          \hline
           5&  0&  0 &  0& 0 & 53130   \\
          \hline
        \end{tabular}  }
        \end{minipage}
\hspace{.2cm}
        \begin{minipage}[t]{8cm}
 {\footnotesize       \begin{tabular}[t]{|c||c|c|c|c|c|c|}
\hline
          & 1 & 2 &3&4 &5 &6  \\
          \hline
          \hline
           1& 216 & 216 & 0 &  0& 0 &  0 \\
          \hline
           2& 0 & 9234 & 1368 & 0 & 0&0   \\
          \hline
           3& 0 & 0 &  141432& 2106  & 0&  0 \\

          \hline
           4&0  &  0&  0&  881469 & 0 & 0  \\
          \hline
           5& 0 & 0  & 0 & 0 & 2261952& 0  \\
          \hline
           6& 0 & 0 &  0&  0 & 0& 1947792  \\
           \hline
        \end{tabular}}
        \end{minipage}
    \caption{\color{black} Count of numerically computed ranks of minors of
    $A_{\Z_5, g}$ and $A_{\Z_6,g}$
    for randomly generated $g$.
    Columns correspond to the dimension of square submatrices and rows to the
    rank of submatrices considered. \label{table:ranksOfMinorsVg} }\color{black}
\end{table}

To this end, we define for $1\leq k\leq|G|$
\begin{equation}\label{equation:tensorphi}
\phi(G,k):=\max\limits_{g\in \C^G{\setminus}\{0\}} \min \big\{ \|V_g
f\|_0: \ f\in \C^G\text{ and } 0<  \|f\|_0\leq k\big\}.
\end{equation}
Using this notation, Theorem~\ref{theorem:p2+1} indicates that
$\phi(\Z_p,k)=p^2-k+1$ for $p$ prime. Taking $\max$ and $\min$ is
justified due to the compactness of the unit ball in $\C^G$. In
fact, we have
\begin{proposition}\label{proposition:AlwaysAlmostEvery}
For almost every $g\in \C^G$, $\displaystyle
\min\limits_{0<\|f\|_0\leq k} \|V_g f\|_0 =\phi(G,k)$ for all $k\leq
|G|$.
\end{proposition}

In the following, we set $Q_{A,B}(z)=\det (A_{G,z})_{A,B}^\ast
(A_{G,z})_{A,B}$, $z=(z_0,z_1,\ldots, z_{|G|{-}1})$, for $A\subseteq
G$ and $B\subseteq G{\times}\widehat G$. $Q_{A,B}$ is a homogeneous
polynomial in $z_0,z_1,\ldots,z_{|G|{-}1}$ of degree $2|A|$.

\begin{lemma}\label{lemma:AEwindowAlwaysGood}
The vector $g\in \C^G$ satisfies $\min\limits_{0<\|f\|_0\leq k}
\|V_g f\|_0 \geq  l$ if and only if $Q_{A,B}(g)\neq 0$ for all
$A\subseteq G$ with $|A|=k$ and all $B\subseteq G{\times}\widehat G$
with $|B|=|G|^2-l+1$.
\end{lemma}

\begin{proof}
Fix $A\subseteq G$ with $|A|=k$ and $g \in \C^G$. Then $g$ satisfies
$\|V_g f\|_0\geq l$ for all $f$ with $\supp f\subseteq A$ if and
only if $\langle f|_A,\pi(\lambda) g|_A \rangle= \langle
f,\pi(\lambda) g \rangle\neq0$ for at least $l$ elements $\lambda
\in G{\times}\widehat G$ for all $f$ with $\supp f\subseteq A$,
i.e., for at most $|G|^2-l$ vectors in $\{\pi(\lambda)g\}$  we have
$\langle f,\pi(\lambda) g \rangle=0$ for $\supp f\subseteq A$. This
is equivalent to $\{\pi(\lambda)g|_A\}_{\lambda\in B}$ spans $\C^A$
whenever $|B|=|G|^2-l+1$. That is, if and only if $\rank
(A_{G,g})_{A,B}=|A|$ for all $B$ with $|B|=|G|^2-l+1$. But this is
equivalent to $Q_{A,B}(g)\neq 0$ for all $|B|=|G|^2-l+1$. The result
follows since for each $f$ with $\|f\|_0\leq k$ exists $A\subseteq
G$ with $|A|=k$ and $\supp f\subseteq A$.
\end{proof}

{\it Proof of Proposition~\ref{proposition:AlwaysAlmostEvery}.}
Lemma~\ref{lemma:AEwindowAlwaysGood} and $\min\limits_{0<\|f\|_0\leq
k} \|V_{g_k} f\|_0 \geq \phi(G,k)$, $k\leq |G|$, for some
$g_k\in\C^G{\setminus}\{0\}$ imply that $Q_{A,B}\not\equiv 0$ for
all pairs $A\subseteq G$ and $B\subseteq G{\times}\widehat G$ with
$|B|=|G|^2-\phi(G,|A|) + 1$. Hence, $
Q=\prod\limits_{A,B:\,|B|=\phi(G,|A|){+}1}Q_{A,B}\not\equiv 0$. This
implies that $Q(g)\neq 0$ for almost every $g\in\C^G$ and therefore,
for almost every $g\in\C^G$ we have $\displaystyle
\min\limits_{0<\|f\|_0\leq k} \|V_g f\|_0\geq \phi(G,k)$ for all
$k\leq |G|$. \hfill $\Box$

To obtain  bounds on $\phi(G,k)$ for groups of non-prime order, we
shall apply Meshulam's strategy to the function $\phi$.

\begin{proposition}\label{proposition:phi}
Let $H$ be a subgroup of the finite Abelian group $G$. For $k\in\N$
exist $s,t\in\N$ with $st \le k$ such that
\begin{equation}\label{equation:MeshulamSTFT}
\phi (G,k) \ge \phi (H,s)\phi(G/H, t)\end{equation}
\end{proposition}

\vpropend
\begin{proof}
In the following, we express the short--time Fourier transformation
for functions defined on $G$ as two consecutive short--time Fourier
transformations. We apply again the notation from the proof of
Theorem~\ref{theorem:meshulam}, i.e., $H=\{x_i\}=\{y_i\}$ and
$\{x_j\}=\{y_j\}$ is a set of coset representatives of the quotient
group $G/H$.  As before $H^\perp = \{\xi_j \in \widehat{G}:
\xi_j(H)=1\}$ and $\{\xi_i\}$ is a set of coset representatives of
$\widehat{G}/ H^\perp$.

Set
$$
    \displaystyle \phi_H (G,k)=\max_{g_1\in \C^H,\ g_2\in \C^{G/H}}\,
        \min \big\{ \|V_{g_1\otimes g_2}
f\|_0: \ f\in \C^G\text{ and } 0<  \|f\|_0\leq k\big\} \, ,
$$
where $g_1{\otimes} g_2(x_i+x_j)=g_1(x_i)g_2(x_j+H)$. Clearly
$\phi(G,k) \ge \phi_H(G,k)$, so \eqref{equation:MeshulamSTFT}
follows from $\phi_H (G,k) \ge \phi (H,s)\phi(G/H, t)$, which we
shall show below. First, note that a similar argument as is used in
Proposition~\ref{proposition:AlwaysAlmostEvery} gives that for
almost every pair $(g_1,g_2)$,
$$
    \phi_H (G,k)  =  \min_{0<\|f\|_0\leq k} \|V_{g_1\otimes g_2}
    f\|_0,\quad 1\le k\le |G|.
$$
Therefore, we can pick $g_1$ and $g_2$ so that for all possible
$k,s,t$,
\begin{eqnarray}
\phi_H (G,k)  =  \min_{0<\|f\|_0\leq k} \|V_{g_1\otimes g_2} f\|_0,
    \ \ \phi (H,s)  =  \min_{0<\|f_1\|_0\leq s} \|V_{g_1} f_1\|_0,\ \ \phi (G/H,t)
    =   \min_{0<\|f_2\|_0\leq t} \|V_{g_2} f_2\|_0\,
    .\label{equation:FelixHeart}
\end{eqnarray}

We  fix $x =x_i + x_j$ and $\xi =\xi_i + \xi_j$, and compute as in
the proof of
Proposition~\ref{proposition:meshulam-induction-argument}
\begin{eqnarray*}
    V_{g_1\otimes g_2} f(x,\xi)
        &=& \sum_{y_j}\sum_{y_i}
            f(y_i {+} y_j)\, \overline{g_1(y_i{-}x_i)}\, \overline{g_2(y_j {-} x_j +
            H)}\,
            \overline{\langle \xi_i',y_i\rangle}_H \overline{\langle\xi_i, y_j\rangle}_G
            \overline{\langle\xi_j',y_j+H\rangle}_{G/H}\\
        &=& \sum_{y_j}\overline{g_2(y_j {-} x_j + H)}\,  \overline{\langle\xi_i, y_j\rangle}_G
            \overline{\langle\xi_j',y_j+H\rangle}_{G/H}\sum_{y_i}
            f(y_i {+} y_j)\, \overline{g_1(y_i{-}x_i)}
            \overline{\langle \xi_i',y_i\rangle}_H
\end{eqnarray*}
where we used $\xi_j \in H^\perp$, i.e., $\langle\xi_j, y_i
\rangle_G =1$. For
$$
    F_H(x_i,\xi_i,y_j):= \overline{\langle\xi_i, y_j\rangle}_G \sum_{y_i}
            f(y_i {+} y_j)\, \overline{g_1(y_i{-}x_i)}\,
            \overline{\langle \xi_i',y_i\rangle}_H
$$
we have
\[
F_H(x_i,\xi_i,y_j)=
    \overline{\langle\xi_i, y_j\rangle}_G
        V_{g_1}T_{-y_j}f(x_i,\xi_i')
\]
and $ V_g f(x, \xi)= \big( V_{g_2} F_H(x_i,\xi_i, \, \cdot)
\big)(x_j{+}H,\xi_j') .$

We fix now $f$ such that $\|f\|_0 \le k$. Let $t = |\{y_j: \supp f
\cap y_j {+} H \ne \varnothing\}|$. If for some $y_j,\, \supp f \cap
y_j + H = \varnothing$, then $F_H(\cdot\,,\cdot\,, y_j) \equiv 0$
too. Therefore, $\|F_H(x_i,\xi_i,\cdot\,)\|_0\le t$ and using
(\ref{equation:FelixHeart})
 we obtain $\|V_{g_2}
F_H(x_i,\xi_i,\cdot\,,\cdot\,)\|_0 \ge \phi (G/H, t)$. Also, by
distributing $\supp f$ over  $t$ cosets of $H$ in $G$, there is a
coset $y_{j_0}{+}H$ with $| \supp f \cap\,  y_{j_0} {+} H|=s \le
k/t$. Because $F_H(\cdot\,,\cdot\,, y_{j_0})$ is, up to a nonzero
factor, the partial short--time Fourier transform of $T_{-y_{j_0}}f$
with window $g_1$ on that coset,
$$
\|F_H(\cdot\,, \cdot\,,y_{j_0})\|_0=\|V_{g_1}T_{-y_{j_0}}f\|_0 \ge
\phi (H, s).
$$
We have obtained that the set $\Lambda=\{(x_i,\xi_i')\in H{\times}
\widehat H :F_H(x_i, \xi_i,y_{j_0}) \ne 0\}$ has at least
$\phi(H,s)$ elements so
\begin{eqnarray*}
  \|V_g f(x_i{+}x_j, \xi_i {+} \xi_j)\|_0
    &=& \sum_{(x_i,\xi_i')\in H{\times}\widehat H} \|V_{g} f(x_i,\xi_i,\cdot\,, \cdot)\|_0
    \ge \sum_{(x_i,\xi_i)\in  \Lambda} \|V_{g_2} F_H(x_i,\xi_i,\cdot\,, \cdot)\|_0
    \\
    &\ge& \phi(H,s)\phi(G/H, t)\,.
\end{eqnarray*}
This inequality holds for all $V_g f$ with $0 < \|f\|_0 \le k$ and
therefore, $\phi_H(G,k) \ge \phi(H,s) \phi(G/H, t)$.
\end{proof}
\begin{theorem} \label{theorem:OurMainResult} For any finite Abelian group $G$ and $k\leq |G|$, let
$d_1$ be the largest divisor of $|G|$ which is less than or equal to
$k$ and let $d_2$ be the smallest divisor of $|G|$ which is larger
than or equal to $k$. Then
\begin{eqnarray}\label{equation:OurMainResult}
\phi(G,k) \geq \frac{|G|^2}{d_1d_2}(d_1+d_2-k).
\end{eqnarray}
\end{theorem}

\vspace{-.3cm}
\begin{proof}
The function $v(n, k) = n\, u(n,k) =
\tfrac{n^2}{d_1d_2}(d_1+d_2-k)$,
 is submultiplicative since $u$ is \cite{Mes05}, i.e., $v(a,b)v(c,d) \ge
v(ac,bd)$. We proceed by induction on $|G|=n$. Suppose
(\ref{equation:OurMainResult}) holds for $|G|=1, \ldots, n{-}1$. If
$n$ is prime, then Proposition \ref{theorem:p2+1} implies $v(n,k) =
n(1+n-k) < n^2-k+1 = \phi (\Z_p, k)$ for all $k$.  Else, we choose a
nontrivial divisor $d$ of $n$, and let $H$ be a subgroup of $G$ of
order $d$. By Proposition \ref{proposition:phi}, there exist $s, t$
with $1 {\le} s {\le} d,\ 1 {\le} t{\le} \min\{\tfrac k s, \tfrac n
d\}$ such that $\phi(G, k) \ge \phi(H, s)\phi(G/H, t)$. Therefore,
$\phi(G, k) \ge v(d, s)v(\tfrac n d, t) \ge v (n, st) \ge v (n, k)$.
\end{proof}

For the case $G=\Z_{pq}$, we can improve this estimate by finding
the convex hull of all pairs $(|H|,|G/H|)$ for all subgroups $H$ of
$G$ as in \cite{Mes05}.


\begin{proposition}\label{proposition:Felix}
Let $G=\Z_{pq}$ with $q < p$ and $p,\,q$ prime. Then
\begin{eqnarray}
 \phi(G,k) \geq \left\{ \begin{array}{ll}
         p^2(q^2-k+1)  & \quad  \mbox{if $k<q$};\\
        (p^2-\tfrac k q +1)(q^2-q+1) &\quad  \mbox{else}.\end{array}
        \right.  \label{equation:FelixBound}
\end{eqnarray}
\end{proposition}

\vspace{-.3cm} The proof of Proposition~\ref{proposition:Felix} is
included in the appendix. At $k=q$, the two lower bounds in
(\ref{equation:FelixBound}) coincide and lead to what a geometric
argument shows to be the optimal value that can be obtained using
$g=g_1\otimes g_2$. So the two straight lines give a convex hull
similar to \cite{Mes05}. However, as expected, the computational
results are better than those given in (\ref{equation:FelixBound}),
since the tensor approach cannot be used to find optimal bounds for
$\phi(G,k)$. See Table~\ref{table:LowerBoundVgfZ6} for an
illustration of (\ref{equation:FelixBound}) for $G=\Z_6$.

\begin{table}[th]
\begin{center}

   {\footnotesize     \begin{tabular}{|c||c|c|c|c|c|c|}

\hline

         $\|f\|_0$   & 1 & 2 &3&4 &5 &6  \\

          \hline \hline

           Theorem~\ref{theorem:OurMainResult} & 36 &18  & 12 &  10&    8 & 6\\

          \hline

           Proposition~\ref{proposition:Felix} &  36& 26 & 25 &23  &   22 & 20\\

          \hline
        \end{tabular}  }

    \caption{\color{black} Lower bounds for $\|V_g f\|_0$ given by Theorem~\ref{theorem:OurMainResult}
    and Proposition~\ref{proposition:Felix} for $G=\Z_6$ and almost every $g\in \C^{\Z_6}$.
    \label{table:LowerBoundVgfZ6}
}
\end{center}
\end{table}

\subsection{Outlook}\label{section:outlook}

For $|G|$ prime, Theorem~\ref{theorem:p2+1} characterizes all pairs
$(\| f\|_0, \|V_g f\|_0)$, $f\in \C^G$ which are achieved for almost
every window function $g\in \C^G$. Below, we conjecture a similar
classification result for general finite Abelian Groups.

\begin{conjecture}\label{conjecture:STFT}
For every finite Abelian group $G$ and almost every $g \in \C^G$, we
have $$\big\{(\|f\|_0, \|V_g f\|_0),\
f\in\C^G{\setminus}\{0\}\big\}=\big\{(\,\|f\|_0\, ,\, \|\widehat
f\|_0{+} |G|^2{-}|G|\,),\ f\in\C^G{\setminus}\{0\}\big\}.$$
\end{conjecture}

\begin{figure}[th]
\begin{center}

\includegraphics[height=1.9cm]{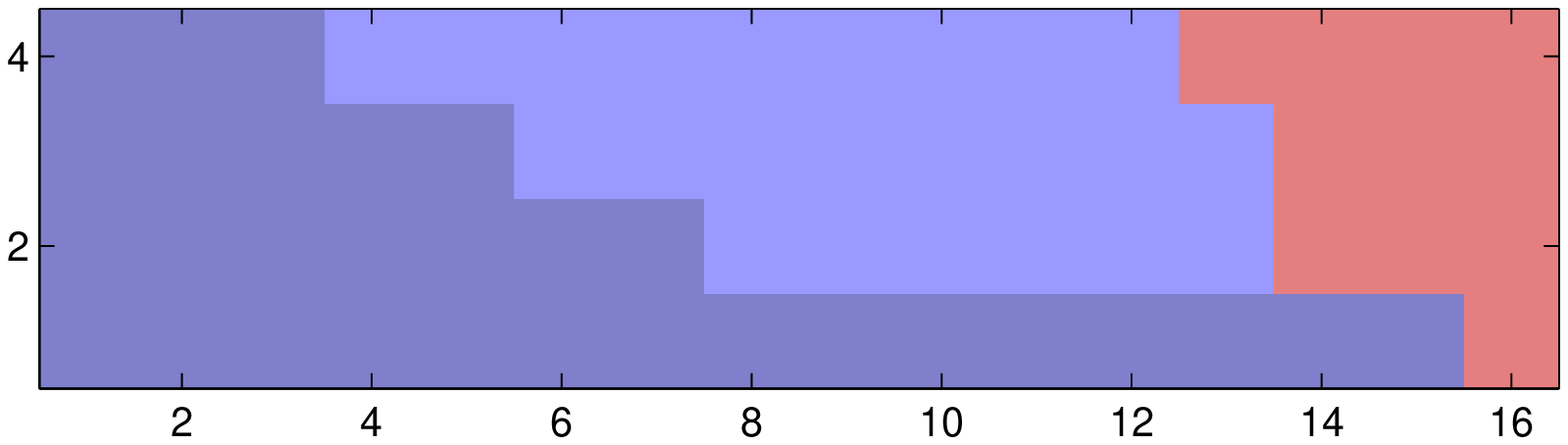}
 \includegraphics[height=2.1cm]{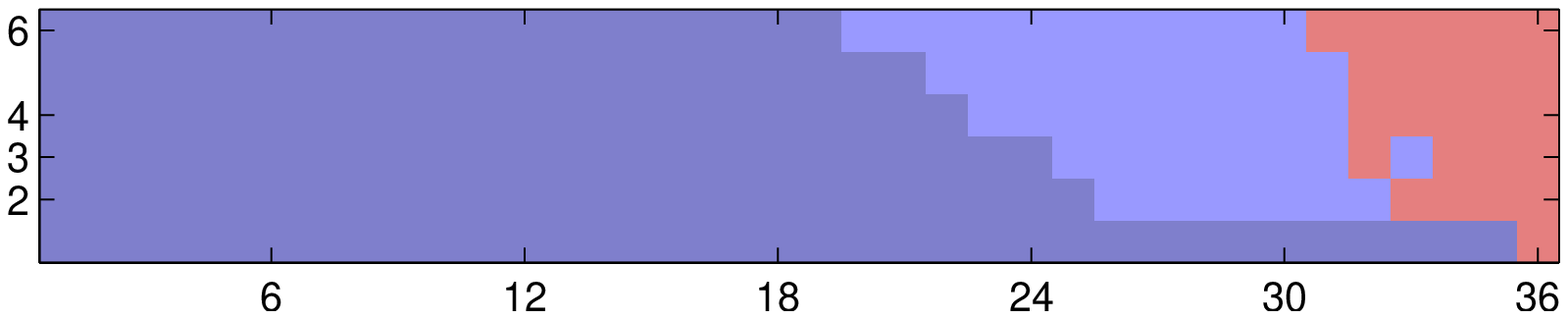}
 \includegraphics[height=2.9cm]{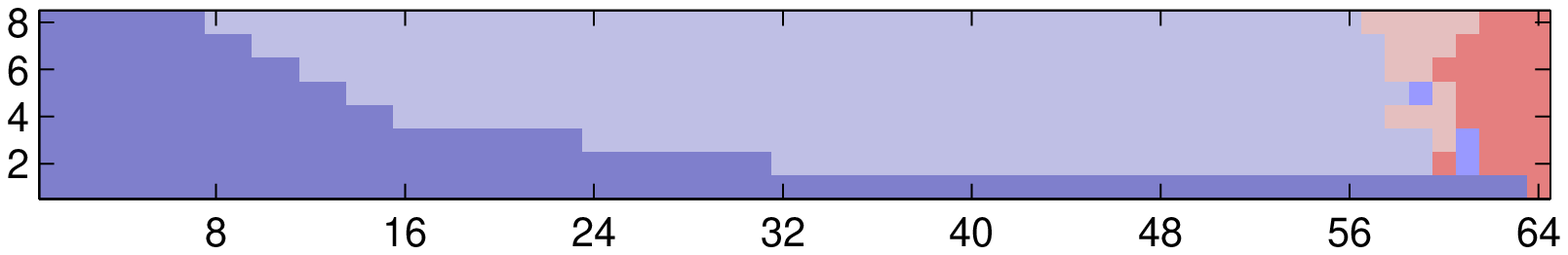}
\end{center}

 \vspace{-5.2cm}
 \hspace{0.7cm}
\color{black}
      $\Z_4$\hspace{6.5cm}
       $\Z_6$

       \vspace{1.6cm}
       \hspace{0.7cm}
        $\Z_8$
\color{black}

    \vspace{2.2cm}
\begin{caption}
  {Same as Figure~\ref{figure:possiblepairsnoncyclicVgprime} for $G=\Z_4$, $G=\Z_6$, and
  $G=\Z_8$. The color coding from Figure~\ref{figure:ColorMap} is
  applied in accordance with Conjecture~\ref{conjecture:STFT} and
  numerical experiments based on Lemma~\ref{lemma:minorranks}.}\label{figure:possiblepairsnoncyclicVg}
   \end{caption}
\end{figure}

This conjecture is illustrated in
Figure~\ref{figure:possiblepairsnoncyclicVg}. As noted earlier, the
numerical testing based on the rank of submatrices of $A_{G,g}$ is
very cost intensive since the number of submatrices that have to be
considered grows combinatorially.


Consequences of this conjecture are discussed in
Section~\ref{section:applications}. Here, we state some preliminary
observations regarding Conjecture~\ref{conjecture:STFT}.


For example, the technique used to prove Theorem~\ref{theorem:p2+1}
possesses certain degrees of freedom, that is, we only need to show
that a particular product of minors is nonzero. Nevertheless, these
degrees of freedom do not allow us to prove
Conjecture~\ref{conjecture:STFT}. For example, for $G=\Z_4$, we can
choose the $4\times 4$ submatrix
$$
   M(z)=\big(A_{\Z_4},(z_0,z_1,z_2,z_3)\big)_{\{0,1,4,12\},\{0,1,2,3\}}= \mysmallmatrix{4}{
        z_0 & z_0 & z_3 & z_1\\
        z_1 &-z_1 & z_0 & z_2\\
        z_2 & z_2 & z_1 & z_3\\
        z_3 & -z_3 & z_2 & z_0
    }
$$
In this submatrix, none of the monomials that is ``maximal" in the
sense described above, namely the monomials $z_0^3z_2$, $z_1^3z_3$,
$z_2^3z_0$, and  $z_3^3z_1$,  has a nonzero coefficient
in the polynomial  $P(z_0,z_1,z_2,z_3)= \det M(z) = -2 z_0^2 z_1^2-2
z_1^2 z_2^2-2 z_0^2z_3^2-2z_2^2 z_3^2 -4 z_0z_1z_2z_3 \not\equiv 0$.


Using Proposition~\ref{proposition:DFTminors},  we  derive a partial
result on nonzero minors of $A_{\Z_n,g}$.

\begin{proposition}\label{proposition:V_gAdjacent}
For every $n$, any minor of the full Gabor system matrix
$A_{\Z_n,g}$, where the columns corresponding to each fixed
translation are adjacent with respect to modulation is nonzero for
almost every $g$. The same holds for a minor corresponding to a
submatrix of size $n{\times}n$, where the columns corresponding to
each fixed modulation are adjacent with respect to translation.
\end{proposition}





%

\section{Applications}\label{section:applications}

We shall now turn to applications of the results stated in
Section~\ref{section:STFTuncertainty} to communications engineering
and, in the subsequent section, to the problem of recovering sparse
signals from incomplete data.

\subsection{Gabor frames, erasures, and the
identification of operators}\label{subsection:applications-Gabor}



We are interested in transmitting information in the form of the
entries of a vector $f\in \C^G$ over a channel in such a way that
recovery of the information at the receiver is robust to errors
introduced by the channel. In particular, we will focus on two
problems. First, we shall discuss transmission over a channel with
erasure, i.e., some of the vector entries may be lost during
transmission.  Second, we discuss the so-called identification
problem for another class of operators, namely, of linear
time--variant operators which play a central role in wireless and
mobile communications. Clearly, knowledge of the operator at hand
would help to counteract disturbances that were caused during
transmission.

We begin with a brief discussion of the recovery of information from
a vector that suffered erasures. Rather then sending the information
in raw form, i.e., sending vector entries one-by-one, information is
being coded prior to transmission. For example, we can choose a
frame $\{\varphi_k\}_{k\in K}$ for $\C^G$ and send the coefficients
$\langle f,\varphi_k\rangle$, $k\in K$. If none of the transmitted
coefficients are lost, the receiver can  use a dual frame
$\{\widetilde \varphi_k\}$ of $\{\varphi_k\}$ and recover $f$ using
(\ref{equation:FrameInversion}). In fact, even if some coefficients
are lost and only $\langle f,\varphi_k\rangle$ is received for $k\in
K'\subset K$, then the information can still be recovered if and
only if $\{\varphi_k\}_{k\in K'}$ remains a frame. This necessitates
that $|K'|\geq |G|=\dim \C^G$.

\begin{definition}\label{def:maxrobust}
A frame ${\cal F}=\{\varphi_k\}_{k\in K}$ in $\C^G$ is {\em
maximally robust to erasures} if the removal of any $l\le |K|-|G|$
vectors from ${\cal F}$ leaves a frame.
\end{definition}

Similarly, we give
\begin{definition}
A set of $m$ vectors in $\C^G$ is {\em in general position}, if any
collection of at most $|G|$ of these vectors are linearly
independent.
\end{definition}

Before giving slight generalizations of results from \cite{LPW05} on
Gabor frames that are maximally robust to erasure in
Theorem~\ref{theorem:equivalence}, we introduce some  vocabulary and
notation regarding the previously mentioned operator identification
problem.

%

\begin{definition}
A linear space of operators $\mathcal H$ mapping $\C^A$ to $\C^B$ is
called {\em identifiable with identifier $g\in \C^A$} if the linear
map $\varphi_g:\mathcal H\longrightarrow \C^B,\ H\mapsto Hg$ is
injective, i.e., if $Hg\neq 0$ for all $H\in \mathcal
H{\setminus}\{0\}$.
\end{definition}

Time--variant communication channels, for example, multipath
channels in wireless telephony, are often modeled through a
combination of translation operators (time--shift, delay) and
modulation operators (frequency shifts that are caused by the
Doppler effect). Therefore, identification of $\displaystyle
\mathcal H_\Lambda= \{ \sum_{\lambda\in\Lambda}c_\lambda
\pi(\lambda),\
  c_\lambda\in\C
  \}$ for $\Lambda\subseteq G{\times}\widehat G$  is
quite relevant.


\begin{theorem}\label{theorem:equivalence}
For $g\in\C^G{\setminus}\{0\}$, the following are equivalent:
\begin{enumerate}

\item Every minor of $A_{G,g}$ of order $|G|$ is nonzero.

\item \label{equivalence:Haar}
The vectors from the Gabor system $\{\pi(\lambda)g\}_{\lambda\in
G{\times}\widehat G}$ are in general position.

\item \label{equivalence:maxrobust} The Gabor system $\{\pi(\lambda)g\}_{\lambda\in G{\times}\widehat G}$
is an equal norm tight frame which is maximally robust to erasures.

\item \label{equivalence:zerosets} For all $f\in \C^G{\setminus} \{0\}$ we have
$\|V_g f\|_0\geq |G|^2{-}|G|{+}1$.

\item \label{equivalence:sampling} For all $f\in \C^G$, $V_g f(\lambda)$, and, therefore, $f$, is completely
determined by its values on any set $\Lambda$ with $|\Lambda|=|G|$.

\item \label{equivalence:identify}$\mathcal H_\Lambda$ is identifiable by $g$ if and only if
$|\Lambda|\leq |G|$

\end{enumerate}
\end{theorem}

For $|G|$ prime, Theorem~\ref{theorem:p2+1} ensures the existence of
$g\in\C^G$ which satisfy parts {\it 1-5} in
Theorem~\ref{theorem:equivalence}.  A verification of
Conjecture~\ref{conjecture:STFT} would also confirm the existence of
$g\in \C^G$ satisfying Theorem~\ref{theorem:equivalence} part
{\it\ref{equivalence:zerosets}}, and therefore
Theorem~\ref{theorem:equivalence} parts {\it 1-5} for general finite
Abelian groups.

\begin{remark}\rm
To our knowledge, the only known {\it equal norm tight frames that
are maximally robust to erasures} are so-called harmonic frames (see
Conclusions in \cite{CK03}). Harmonic frames for $\C^n$ with $m\geq
n$ elements are obtained by deleting uniformly $m-n$ components of
the characters of $\Z_m$ \cite{CK03}. Similarly,
Theorem~\ref{theorem:LaPfaWa} together with Theorem~\ref{unitg}
provides us with a large class of equal norm tight frames with $p^2$
elements in $\C^n$ for $n\leq p$. Namely, we can choose $g\in
(S^1)^p$ and remove $p-n$ components of the equal norm tight frame
$\{\pi(\lambda)g \}_{\lambda\in G{\times}\widehat G}$ in order to
obtain an equal norm tight frame which is maximally robust to
erasure. Note that this frame is not a Gabor frame proper. Reducing
the number of vectors in the frame to $m\leq p^2$ vectors leaves an
equal norm frame which is maximally robust to erasure but which
might not be tight. This holds for harmonic frames too and with the
restriction to $p$ prime, we have shown the existence of Gabor
frames which share the usefulness of harmonic frames when it comes
to transmission of information through erasure channels.
\end{remark}

Background and more details on frames and erasures can be found in
\cite{CK03,GK01,SH03} and the references cited therein.

\subsection{Signals with sparse representations}\label{subsection:applications-sparsity}

In Section~\ref{subsection:applications-Gabor} we discussed the
recovery of signals or operators from $|G|$ known complex numbers.
Here, we will use the functions $\displaystyle\phi$ and
$\displaystyle \theta$ which were defined in
Section~\ref{section:meshulamFT} and
Section~\ref{section:meshulamSTFT} to refine some of these findings.
That is, we show that a function/signal which can be represented as
a linear combination of a small number of pure frequencies or of a
small number of time--frequency shifts of a fixed function $g$, can
be recovered from fewer than $|G|$ of its values. Our brief
discussion is based on the most basic ideas and results from the
theory of sparse signal recovery.

There exist a number of entry points to the theory of sparse signal
recovery. Here, we shall consider  dictionaries
$\Dcal=\{g_0,g_1,\ldots,g_{N-1}\}$ of $N$ vectors in $\C^n$, or
equivalently, in $\C^G$. For $k\leq n=|G|$ we shall examine the sets
$$
 \Sigma^\Dcal_k=\{f\in\C^n:\ f=M_\Dcal\, c=\sum_r c_r g_r, \text{ with } \|c\|_0\leq k
 \}\,.
$$
The central question is: {\it  how many values of $f\in
\Sigma^\Dcal_k$ need to be known (or stored), in order that
$c\in\C^N$ with $f=\sum_r c_r g_r$ and $\|c\|_0\leq k$, and
therefore $f$,  is uniquely determined by the known data?}

To this end, we set
$$
\psi(\Dcal,k)= \min\big\{\|f \|_0:\ f \in \Sigma^\Dcal_k\big\}\,,
$$
and observe the following well known result.
\begin{proposition}\label{proposition:2k}
  Any $f\in \Sigma^\Dcal_k$ is fully determined by any choice of
  $n{-}\psi(\Dcal,\,2k\,){+}1$ values of $f$.
\end{proposition}

Note that unlike in Theorem~\ref{theorem:equivalence},  we do not
assume knowledge of the set $\supp c$ for $c$ with $M_\Dcal c=f$,
$\|f\|_0$ in Proposition~\ref{proposition:2k} and in the following.

\begin{proof}
Assume that for some $B\subset \C^n$ with
$|B|=n{-}\psi(\Dcal,\,2k\,){+}1$, two coefficient vectors
$c_1,c_2\in\C^N$ exist that satisfy $r_B M_\Dcal c_1=r_B f=r_B
M_\Dcal c_2$ and $\|c_1\|_0,\|c_2\|_0\leq k$. Then
$\|c_2-c_1\|_0\leq 2k$ with $\|M_\Dcal (c_2-c_1)\|_0\leq n -|B|=
n-(n{-}\psi(\Dcal,\,2k\,){+}1)=\psi(\Dcal,\,2k\,){-}1$, a
contradiction.
\end{proof}

A classical dictionary for $\C^G$ is $\Dcal_G =\{ \xi \}_{\xi \in
\widehat G}$, where $G$ is a finite Abelian group. Then
$$
\psi(\Dcal,k)= \min\big\{\|f \|_0:\ f \in \Sigma^\Dcal_k\big\}=
\min\big\{\|\widehat f \|_0:\ \|f\|_0 \leq k \big\}=\theta(G,k)\,.
$$
This equality together with Proposition~\ref{proposition:2k}
demonstrates the relevance of the results cited in
Section~\ref{section:FourierUncertainty} for the recovery of signals
with limited spectral content. For example,
Theorem~\ref{theorem:meshulam} shows that for any finite Abelian
group of order $16$ we have $\theta(G,6)\geq 3$. In fact, our
computations that are illustrated in
Figure~\ref{figure:possiblepairsnoncyclic2} show that
$\theta(G,6)=4$ for $|G|=16$, and, hence,  any
$f\in\Sigma_3^{\Dcal_{G}}=\{f:\ \|\widehat f\|_0\leq 3 \}$ can be
recovered from any choice of $|G|-\theta(G,2\cdot 3)+1=16-4+1=13$
values of $f$.  For $f\in\Sigma_3^{\Dcal_{\Z_{17}}}$ on the other
side, Theorem~\ref{theorem:tao} implies that $f$ is already fully
determined by $|\Z_{17}|-\theta(\Z_{17},2\cdot 3)+1=17-(17-6+1)+1=6$
of its values.

The results in Section~\ref{section:STFTuncertainty} which involve
the function $\phi$ are relevant to determine vectors which have
sparse representations in the dictionary $\Dcal_{A_{G,g}}$ which
consists of the columns of $A_{G,g}$. In fact, we have $F\in
\Sigma^{\Dcal_{A_{G,g}}}_k$ if and only if $F=V_gf$ for some
$f\in\C^G$ with $\|f\|_0\leq k$ and, therefore,
$$
\psi( \Dcal_{A_{G,g}} ,k)= \min\big\{\|V_g f \|_0:\ \|f\|_0\leq
k\big\}=\phi(G,k)\,.
$$
For $|G|$ prime for example, this leads to the following short--time
Fourier transform version of Theorem 1.1 in \cite{CRT04}.

\begin{theorem}\label{theorem:robustUncertainty} Let $g\in\C^{\Z_p}$,
$p$ prime, satisfy the conclusion of Theorem~\ref{theorem:p2+1}.
Then any $f\in\C^{\Z_p}$ with $\|f\|_0\leq \tfrac 1 2 |\Lambda|$,
$\Lambda\subset \Z_p{\times}\widehat{\Z_p}$ is uniquely determined
by $\Lambda$ and $r_\Lambda V_g f$.
\end{theorem}

In terms of sparse representations, the Gabor frame dictionary
$\{\pi(\lambda) g\}_{\lambda\in G{\times}\widehat G}$ of
time--frequency shifts of a prototype vector $g$, i.e., the
dictionary consisting of the rows of $A_{G,g}$, appears to be more
interesting.  Rudimentary numerical experiments based on
Lemma~\ref{lemma:minorranks} give some indication that for any
Abelian group $G$, and almost every $g\in\C^G$, we have for $k\leq
|G|$
$$
\psi(\{\pi(\lambda) g\}_{\lambda\in G{\times}\widehat
G},k)=\theta(G,k).
$$
For $|G|$ prime,  Theorem~\ref{theorem:LaPfaWa} implies that
$\psi(\{\pi(\lambda) g\}_{\lambda\in G{\times}\widehat G},k)=p-k+1
=\theta(G,k)$, and analogous to
Theorem~\ref{theorem:robustUncertainty}, we obtain
\begin{theorem}\label{theorem:robustUncertainty2} Let $g\in\C^{\Z_p}$,
$p$ prime, satisfy the conclusion of Theorem~\ref{theorem:p2+1}.
Then any $f\in\C^{\Z_p}$ with $f=\sum_{\lambda\in\Lambda} c_\lambda
\pi(\lambda)g$, $\Lambda \subset \Z_p{\times}\widehat{\Z_p}$ is
uniquely determined by $B$ and $r_B f$ whenever $|B|\geq
2|\Lambda|$.
\end{theorem}
Note that similar to before,  the recovery of $f$ from $2|\Lambda|$
samples of $f$ in Theorem~\ref{theorem:robustUncertainty2} does not
require knowledge of $\Lambda$.

\section{Appendix}\label{section:appendix}

\subsection{Proof of Lemma~\ref{lemma:minorranks}}

If $f\in \C^n$ with $\|f\|_0=k$ and $\|Mf\|_0=l$, then $A=\supp f$
and $B^c= \supp Mf$ satisfy $0\neq r_A f\in \ker M_{A,B}$, so $\rank
M_{A,B}<|A|$. Moreover, for $a \in A$, $\supp f = A$ implies $f
\notin \{g: \|g\|_0<|A|\} \supset i_{A{\setminus}\{a\}} \ker
M_{A{\setminus}\{a\},B}$ and, hence,
$$f \in i_{A}\ker M_{A,B}{\setminus} i_{A{\setminus}\{a\}}\ker M_{A{\setminus}\{a\},B}\, .$$
So $\dim \ker M_{A,B}\ge\dim \ker M_{A{\setminus}\{a\},B} + 1$. We
conclude that for all $a \in A$,
$$\rank M_{A{\setminus}\{a\},B} \le \rank M_{A,B} =
|A|-\dim\ker M_{A,B} \leq |A|-\dim\ker
M_{A{\setminus}\{a\},B}-1=\rank M_{A{\setminus}\{a\},B}$$ which
implies $\rank M_{A{\setminus}\{a\},B}=\rank M_{A,B}.$
Also, $\supp Mf = B^c$, so for $y \in B^c$, $Mf(y)\neq 0$.
Therefore, $f\notin \ker M_{A,B\cup\{y\}}$and so $f \in i_A \ker
M_{A,B}{\setminus} i_A \ker M_{A,B\cup\{y\}}$. This implies
\begin{eqnarray*}
      \ \rank M_{A,B}\ =\ |A| - \ker M_{A,B} <|A| - \ker
      M_{A,B\cup\{y\}}=
       \rank M_{A,B\cup\{y\}}\, .
    \end{eqnarray*}
The submatrices considered differ only by one column, so the rank
can increase at most by one and we get $\rank M_{A,B}\ =\ \rank
M_{A,B\cup\{y\}}-1$.

Suppose now  that $A\subseteq\{0,\ldots,n{-}1\}$ and
$B\subseteq\{0,\ldots,m{-}1\}$ with $|A|=k$ and $|B|=m-l$ satisfy
(\ref{equation:minorranks}). This implies $\dim \ker M_{A,B}\geq 1$
and that for any $a\in A$,
  $$
    \dim \ker M_{A{\setminus}\{a\},B} = |A|-1-\rank
    M_{A{\setminus}\{a\},B}= |A|-1-\rank
    M_{A,B}=\dim \ker M_{A,B}-1.
  $$
So $\ i_{A{\setminus}\{a\}} \ker
M_{A{\setminus}\{a\},B}\varsubsetneqq i_{A}\ker M_{A,B}$, and there
exists $f_a \in i_{A}\ker M_{A,B}{\setminus}
i_{A{\setminus}\{a\}}\ker M_{A{\setminus}\{a\},B}$, so $f_a(a)\neq
0$, $f_a(x)=0$  for $x\notin A$ and $\supp M f_a \cap B =
\varnothing. $

Similarly, (\ref{equation:minorranks}) implies also that for any $y
\in B^c$ we have $i_A \ker M_{A,B\cup\{y\}}\varsubsetneqq i_A \ker
M_{A,B}$, so there exists $g_y$ such that $Mg_y(y)\neq 0$ while
$Mg_y(b)=0$ for all $b\in B$.

To conclude this proof, we enumerate the vectors $f_a$, $a\in A$ and
$g_y$, $y\in B^c$ and choose a linear combination
  \begin{equation}
    f=\sum_{a\in A}c_a f_a+\sum_{y \in B^c}c_y g_y = \sum_{r=0}^{k+l-1} d_r h_r \label{equation:examplemr}
  \end{equation}
with the property that $\supp f=\bigcup_{a\in A}\supp f_a=A$ and
$\supp Mf=\bigcup_{y\in B^c}\supp Mg_y =B^c$.

By construction we have  $\supp f\subseteq A$ and $\supp Mf
\subseteq B^c$. To get the reverse inequality, we assume without
loss of generality that $\min\limits_{x\in \supp h_r} |h_r(x)|=1$
for all $r$, and choose $d_r=N^{2r}$, where $ N{-}1 \ge
\|h_r\|_\infty, \|M h_r \|_\infty,\|M h_r\|_\infty^{-1} $ for
$r=0,1,\ldots,k{+}l{-}1$. Since $f_{a_0}(a_0)\neq 0$ we can find $s=
\max\{r:\ h_r(a_0)\neq 0\}$. Then
$$
    |f(a_0)|    =   \big| \sum_{r=0}^{s} d_r h_r(a_0) \big|
                \geq |N^{2s}h_s(a_0)| - \big| \sum_{r=0}^{s-1} N^{2r} h_r(a_0) \big|
                \geq N^{2s} - (N{-}1) \sum_{r=0}^{s-1} (N^2)^r  =N^{2s}-\tfrac{{N^{2s}}-1}{N+1}>0,
$$
so $a_0\in \supp f$.

Similarly, $M g_{y_0}(y_0)\neq 0$ for fixed $y_0\in B^c$ implies
that for $s= \max\{r:\ Mh_r(y_0)\neq 0\}$ we have
$$
    |Mf(y_0)|    =   \big| \sum_{r=0}^{s} d_r M h_r(y_0) \big|
                \geq |N^{2s} M h_s(y_0)| - \big| \sum_{r=0}^{s-1} N^{2r} M h_r(y_0) \big|
                \geq \tfrac{N^{2s}}{N{-}1} -
                \tfrac{{N^{2s}}-1}{N+1}>0.
$$
We conclude that  $\supp f=A$ and $\supp Mf=B^c$.

\subsection{Proof of Proposition~\ref{proposition:Z6}}
Theorem~\ref{theorem:classicalUncertainty} and
Proposition~\ref{proposition:Gitta} cover all cases but
$(k,l)=(2,4),(3,3),(4,2)$. For  $\omega =e^{2\pi i/6}$, we have
$\mathcal F (1, -1, 0, 1, -1, 0 )= (0, 0, 1{-}\omega^2, 0,
1{-}\omega^4, 0)$, and only the case $(k,l)=(3,3)$ remains to be
excluded.

The assumption $\|f\|_0=3$ leads to three different cases.

\vspace{.1cm}\noindent Case 1. If  $f=(c_0,0,c_2,0,c_4,0)$ then
$\widehat{f}(\xi)= \widehat{f}(\xi+3)$ and if
$f=(0,c_1,0,c_3,0,c_5)$ then $\widehat{f}(\xi)=
-\widehat{f}(\xi+3)$. In either case, $\|\widehat{f}\|_0$ is even
and cannot be $3$.

\vspace{.3cm}\noindent Case 2. If two entries whose indices differ
by $3$ are both nonzero, then the support of the Fourier transform
cannot be $3$ either. To see this, consider without loss of
generality, $f=(c_0,*,*,c_3,*,*)$. Then, for $c_k$, located at
position $k$, being the third nonzero entry, we have
\begin{eqnarray}
\widehat f = (c_0{+}c_3{+}c_k,\,c_0{-}c_3{+}\omega^k c_k,\,
c_0{+}c_3 {+}\omega^{2k} c_k,\,c_0{-}c_3{+}\omega^{3k}c_k,\,
c_0{+}c_3{+} \omega^{4k}c_k,\,c_0{-}c_3{+}\omega^{5k}c_k)\, .
\label{equation:Z6case2}
\end{eqnarray}

If three coordinates of $\widehat f$ are $0$, then two of the
respective sums in (\ref{equation:Z6case2}) contain either both
$c_0+c_3$ or both $c_0-c_3$. Without loss of generality, we assume
that $\widehat f (l_1)= c_0{+}c_3{+}\omega^{l_1 k}c_k\neq 0\neq
c_0{+}c_3{+}\omega^{l_2 k}c_k=\widehat f (l_2)$, $l_1<l_2$. Since
$c_k\neq 0$ we have $\omega^{l_1 k}=\omega^{l_2 k}$ and
$\omega^{(l_2-l_1)k}=1$. Since $k=1,2,4$ or 5, we must have $3$
divides $l_1-l_2$, but that is a contradiction, as of two entries
with distance $3$, one must contain the summand $c_3-c_0$ and one
$c_0+c_3$.

\vspace{.3cm}\noindent Case 3. If all three nonzero entries are
adjacent, then $\widehat f$ must have three adjacent entries as
well, as otherwise, we could just exchange the roles of $f$ and
$\widehat f$ and return to Case 1 or Case 2. Without loss of
generality we assume $f=(c_0,c_1,c_2,0,0,0)$. A modulation in $f$
results in a translation in $\widehat f$, so without loss of
generality, we can also assume the first three entries of $\widehat
f$ to be $0$. Hence,
$$
\mysmallmatrix{3}{
1 &1&1\\
1& \omega& \omega^2\\
1& \omega^2& -\omega } \mysmallmatrix{1}{
c_0\\
c_1\\
c_2}=0 \quad \text{but}\quad \det\mysmallmatrix{3}{
1 &1&1\\
1& \omega& \omega^2\\
1& \omega^2& -\omega }=-1\neq 0$$ and, therefore, $f=0$.


\subsection{Proof of Proposition~\ref{proposition:FourierSetpq}}

The group $\Z_{pq}$ has $(p{-}1)(q{-}1)$ automorphisms, each of them
mapping one of the $(p{-}1)(q{-}1)$ elements of order $pq$ to 1. The
$p{-}1$ automorphisms on the group $\Z_{2p}=\{0,1,2,\ldots,
2p{-}1\}$ will allow us to consider only $f$ with
well-``concentrated'' nonzero entries.

Every automorphism $\sigma$ on $\Z_{pq}$ induces an automorphism
$\tilde\sigma$ on the character group $\widehat{\Z}_{pq}$, which
satisfies $\langle \tilde\sigma(\xi), x\rangle = \langle \xi,
\sigma^{-1}(x)\rangle$. Further,
$$
    \widehat{f{\circ}\sigma}(\xi)
        = \tfrac 1{pq} \sum \limits_{x\in \Z_{pq}}f(\sigma(x))\overline{\langle
        \xi,x\rangle}
        =\tfrac 1{pq} \sum \limits_{y\in \Z_{pq}}f(y)\overline{\langle
        \xi,\sigma^{-1}(y)\rangle}
        =\tfrac 1{pq} \sum \limits_{y\in \Z_{pq}}f(y)\overline{\langle
        \tilde\sigma(\xi),y\rangle}
        = \widehat{f}( \tilde\sigma(\xi))\
$$ 

Let $f\in \C^{\Z_{2p}}$, $p\geq 5$ prime,  be given with
$\|f\|_0=3$. Then at least two of the addresses of the non-zero
elements have the same parity. By a translation of $f$ we can move
those elements to positions $0,2k$, where $k\in\Z_{2p}$. The support
of $\widehat f$ is not affected by this. If $k$ is odd, then $k$ is
a generator of $\Z_{2p}$ and we choose $\sigma_1$ with
$\sigma_1(k)=1$. If $k$ is even, then $p+k$ is odd and we pick
$\sigma_1$ with $\sigma_1(p{+}k)=1$. In either case
$\sigma_1(2k)=2$.  The corresponding automorphism $\tilde\sigma_1 $
in $\fhat\Z_{2p}$ will affect $\supp \widehat f$, but $\|\fhat
f\|_0$ does not change.

Let the third non-zero element have address $r$. If $\sigma_1(r)\neq
p{+}1$, then there are either $p{-}1$ adjacent zeroes among the
addresses $3, \ldots, p{+}1$ or among $p{+}1, \ldots, 2p{-}1$.

In case that $\sigma_1(r)=p{+}1$, then we apply another automorphism
$\sigma_2$ in a similar way as above. If $\tfrac{p{+}1}2$ is a
generator for $\Z_{2p}$, then $\sigma_2(\tfrac{p{+}1}2)=1$,
$\sigma_2(2)=\sigma_2(4\tfrac{p{+}1}2)=4\sigma_2(\tfrac{p{+}1}2)=4$,
and $\sigma_2(p{+}1)=\sigma_2(2\tfrac{p{+}1}2)
=2\sigma_2(\tfrac{p{+}1}2)= 2$. Otherwise, we choose $\sigma_2$ such
that $\sigma_2(p+\tfrac{p{+}1}2)=1$, so
$\sigma_2(p{+}1)=2\sigma_2(p+\tfrac{p{+}1}2)= 2$ and
$\sigma_2(2)=2\sigma_2(p{+}1)=4$. In both cases, $\supp
(f{\circ}\sigma_2{\circ}\sigma_1) =\{0,2,4\}$, so the vector
contains a string of at least $p{-}1$ consecutive zeros on addresses
$5,\ldots, 2p{-}1$.

The following lemma from \cite{DS89} implies that $\|\widehat
{f{\circ}\sigma'{\circ}\sigma}\|_0> p{-}1$ and, therefore,
$\|\widehat f\|_0\geq p$.

\begin{lemma}\label{lemma:ds5} If $\fhat f$ has $N$ nonzero
elements, then $f$ cannot have $N$ consecutive zeros.
\end{lemma}

\subsection{Justification of Figure~\ref{figure:vgf}}

Let $\omega=e^{2\pi i/3}$. For $\|f\|_0=1$, we calculate
$$V_{(a,b,c)}(d,0,0)=(d\overline{a}, \omega^2 d\overline{a}, \omega
d\overline{a},d\overline{c}, \omega^2 d\overline{c}, \omega
d\overline c, d\overline b, \omega^2 d\overline b, \omega d\overline
b)$$ So in any case, $\|V_g f\|_0=3\|g\|_0$, which justifies all
cases involving $\|f\|_0=1$ or $\|g\|_0=1$.

For the case $\|f\|_0=2$ and $\|g\|_0=2$, we note
$\|V_{(1,1,0)}(1,-1,0)\|_0=8$ and $\|V_{(1,1,0)}(1,10,0)\|_0=9$,
which justifies the two red fields. Now assume that there are $f$
and $g$ with $\|f\|_0=\|g\|_0=2$ and $\|V_g f\|_0\leq 7$. Then $V_g
f$ has at least two zero entries. Note that the scalar product of
$f$ and another vector with support size $2$ can only vanish, if
$\supp f=\supp g$. So the zero entries in $V_g f$ must correspond to
the same translation. If we set without loss of generality
$f=(a,b,0), g=(c,d,0)$, then zeros at two different modulations
$M_{j_1}$ and $M_{j_2}$ imply $a\overline c+
\overline\omega^{j_1}b\overline d=0= a\overline c+
\overline\omega^{j_2}b\overline d$, which clearly admits no
nontrivial solution.

For the case $\|f\|_0=2$ and $\|g\|_0=3$ which is equivalent to the
case $\|f\|_0=3$ and $\|g\|_0=2$, we note that
$\|V_{(1,1,1)}(1,-1,0)\|_0=6$, $\|V_{(2,-4,8)}(2,1,0)\|_0=7$,
$\|V_{(1,2,3)}(2,-1,0)\|_0=8$ and $\|V_{(1,2,3)}(1,2,0)\|_0=9$,
which justifies the four red fields. Now assume, there are $f$ and
$g$ with $\|f\|_0=2$, $\|g\|_0=3$ and $\|V_g f\|_0\leq 5$. Then $V_g
f$ has at least four zero entries, in particular two that correspond
to the same translation. Without loss of generality, we assume that
this is the zero-translation and that $f$ is supported in the first
two coordinates, i.e., $f=(a,b,0),\ g=(c,d,e)$. Then we get as
before $a\overline c+ \overline\omega^{j_1}b\overline d=0=
a\overline c+ \overline\omega^{j_2}b\overline d$ which has no
nontrivial solutions.

For the case $\|f\|_0=3$ and $\|g\|_0=3$, we note that
$\|V_{(1,1,1)}(1,1,1)\|_0=3$, $\|V_{(1,1,1)}(1,1,-2)\|_0=6$,
$\|V_{(1,2,5)}(10,5,2)\|_0=7$, $\|V_{(1,2,3)}(-5,1,1)\|_0=8$ and
$\|V_{(1,2,3)}(1,2,3)\|_0=9$, which justifies the five red fields.
Multiplying $f$ or $g$ by a constant does not change $\|V_g f\|_0$,
so we can normalize $f(0)=g(0)=1$. Hence we can set $f=(1,a,b)$,
$g=(1,c,d)$. Then again, $\|V_g f\|_0\leq 5$ implies that $V_g f$
has two zero entries that correspond to the same translation and we
shall assume without loss of generality and for the remainder of
this section that those appear at $x=0$ and $\xi=1,2$, i.e., we have
$$1+\omega a\bar c + \omega ^2 b\bar d =0= 1+\omega ^2 a\bar c+
\omega b\bar d$$ and hence $b\bar d=a\bar c=1$ and $g=\left( 1,
\tfrac{1}{\bar a}, \tfrac{1}{\bar b}\right)$.

Before continuing, we state

\begin{lemma}\label{shearing} Let $S$ be a shearing on $\C^{\Z_3{\times} \Z_3}$, i.e., $S$
translates the ($x=1$)-row of an element in $\C^{3{\times}3}$ by $1$
and  the ($x=2$)-row by 2. Then given $f,g\in \C^{\Z_3}$, there
exist $\tilde f,\tilde g\in \C^{\Z_3}$, such that $\supp V_{\tilde
g}\tilde f$ is the image of $\supp(V_g f)$ under $S$.
\end{lemma}
\begin{proof}Suppose, two vectors $f=(u,v,w)$ and $g=(x,y,z)$ are given, and consider the vectors $\tilde f=(u,v, \omega w)$ and $\tilde g=(x,y,\omega z)$. Then $$
    V_{\tilde g}\tilde f(0,\xi)
        =u\overline x +\overline \omega^\xi v \overline y
            + \overline\omega^{ 2\xi}(\omega
w) (\overline \omega \overline z)
        = u\overline x +\overline \omega^\xi v \overline y + \overline\omega^{ 2\xi} \overline z  w
        = V_g f(0,\xi)\, ,
$$
$$
    V_{\tilde g}\tilde f(1,\xi)
        =u\overline y +\overline \omega^\xi v \overline \omega \overline z
                + \overline\omega^{ 2\xi}(\omega w) ( \overline x)
        = u\overline y
                +\overline \omega^{\xi+1} v \overline z + \overline\omega^{ 2\xi+2}
                    \overline x  w = V_g f(1,\xi+1)\, , $$ and $$
    V_{\tilde g}\tilde f(2,\xi)
        =u\overline\omega\overline z +\overline\omega^{\xi} v
            \overline x + \overline\omega^{2 \xi} \omega w\overline y
        = \overline\omega ( u\overline z +\overline\omega ^{\xi+2}v\overline x
            + \overline\omega ^{2\xi+1}w\overline y)
        = \overline\omega V_g f(2,\xi+2)\,.
$$
As a multiplication by $\overline \omega$ does not change the
support, we get the sheared image of the original support set as
desired.
\end{proof}

We now use Lemma~\ref{shearing} to show that in the case
$\|f\|_0=\|g\|_0=3$, no support size of $4$ is possible. In fact
this would imply that the short--time Fourier transform has five
zeroes, so there is a second row with two zeroes (without loss of
generality the row $x=1$). By shearing we can move them to $\xi=1,2$
without changing the first row, i.e., $$
    \tfrac 1 a +\overline\omega \tfrac a b +\overline\omega^2 b
    =   0
    =   \tfrac 1 a +\overline\omega^2 \tfrac a b +\overline\omega b\,.
$$

This implies $\tfrac 1 a = \tfrac a b = b$ and hence $a=1$,
$a=\omega$ or $a=\omega^2$, and $b=\overline a$ accordingly. This
reduces to the the example for $\|V_g f\|_0 =3$ given above. Thus,
$\|V_g f\|_0=4$ is impossible.

For a support size of $5$, we can use the same argument to exclude
that the remaining two zeroes occur at the same $x$. So in addition
to the two zeros for $x=0$, we can have zeroes at $x=1,2$ and either
$\xi=0$ for both or $\xi=1$ for both. All other combinations can be
reduced to these two by shearing and conjugation (using
$\omega^2=\bar \omega$).

These two cases correspond to solving
$$
    a+\overline\omega^k \tfrac{b}{a} +\overline\omega^{2k}\tfrac{1}{b}
    =0
    =\tfrac{1}{a}+\overline\omega^k \tfrac{a}{b} +\overline\omega^{2k}{b}
$$ for $k=0,1$.
These equations can be solved exactly using Mathematica. The only
solutions are modulations of shearings of the solution with $\|V_g
f\|_0=3$ considered above. So again, it follows that a short--time
Fourier transform with support size $5$ is not possible.

\subsection{Proof of Cheboratev's Theorem~\ref{theorem:tao} based on Thoerem\~ref{theorem:LaPfaWa}.}

Fix $A,\widetilde{A} \subseteq \Z_p$ with $|A|=|\widetilde{A}|$. We
have to show that the restricted Fourier transformation $\mathcal
F_{A\to\widetilde A}: \C^A\to \C^{\widetilde{A}}$ is an isomorphism.
For $g$ such that $A_{\Z_p,g}$ has no zero minors, define
$M_g:\C^p\longrightarrow \C^p$ to be the pointwise multiplication
operator with the vector $g$. Since $g$ has no zero components, $M$
is an isomorphism, and, moreover, $M_g$ restricts to an isomorphism
on $\C^A$. Set $B=\{0\} \times \widetilde{A}$. Therefore, $V_g :
\C^A \longrightarrow \C^B$ is an isomorphism since
$|B|=|\widetilde{A}|$. The result follows since the restricted
Fourier transformation $\mathcal F_{A\to\widetilde A}$ is nothing
but $ P \circ V_g \circ M_g$ where $P$ is the projection of $B=\{0\}
\times \widetilde{A}$ onto $\widetilde{A}$.

\subsection{Proof of Proposition~\ref{proposition:Felix}}

Proposition~\ref{proposition:meshulam-induction-argument} implies
that there exists $s, t$ such that $st\leq k$ and $ \phi(G,k)\geq
\phi(H,s)\phi(G/H,t)$. For $G=\Z_{pq}$ and $|H|=p$, we have
$\phi(H,s)=p^2-s+1$ and $\phi(G/H,t)=q^2-t+1$. As $st \leq k$, we
can find $\overline t\in\R$ such that $ q\geq \overline t \geq t$
and $p \geq \frac k{\overline t} \geq s$. Hence,
\begin{eqnarray*}
\phi(G,k) & \geq & (p^2-s+1)(q^2-t+1)  \geq  (p^2-\tfrac k{\overline
t} +1)(q^2-\overline t+1)\,.
\end{eqnarray*}

So $\phi(G,k)$ must exceed the minimum of $M(u)=(p^2-\tfrac k u
+1)(q^2-u+1)$, where $u$ ranges from $\tfrac k p$ to $q$ since
$\tfrac k u\leq p$ and $u\leq q$ is assumed.
We have $M'(u)=-(p^2+1)+\tfrac{k(q^2+1)}{u^2}=0$ if and only if
$u=\pm\sqrt{k \tfrac{q^2+1}{p^2+1}}$.

As $M(u)\rightarrow - \infty$ for $u \rightarrow 0^+$ and $u
\rightarrow \infty$, the only positive extremum is a maximum and the
minimum is attained in a boundary point. A simple calculation gives
that $M(q) \leq M\left(\tfrac k p\right)$.

For $k<q$, the condition $1\leq s$, $1\leq t$, implies that $t$
ranges only from $1$ to $k$. The same arguments as used above show
again that the minimum is attained at a boundary point and that
$M(1) \geq M(k)$.

\subsection{Proof of Proposition~\ref{proposition:V_gAdjacent}}

As in the proof of Theorem~\ref{theorem:LaPfaWa}, choose $A\subseteq
G$ and $B\subseteq G{\times}\widehat G$ with $|A|=|B|$ and set
$P_{A,B}(z)= \det (A_{\Z_n,z})_{A,B}$, $z=(z_0, z_1,\ldots,
z_{n{-}1})$. In that proof, we identified a ``maximal" term within
$P_{A,B}$, the coefficient of which can be expressed as a product of
different minors of the discrete Fourier matrix $W_{\Z_n}$. Each of
these minors arise from the columns of $P_{A,B}(z)$ that correspond
to a specific translation. By assumption, these columns are adjacent
with respect to modulation in $A_{\Z_n,z}$.

So each of these minors is a minor of the DFT matrix corresponding
to adjacent columns, where each row is multiplied by some factor
$z_i$. Using the multilinearity of the determinant, we can pull the
factors outside. By Proposition~\ref{proposition:DFTminors}, we
conclude that these minors of the DFT-matrix are nonzero, hence also
their product. So the "maximal" term has a nonzero coefficient.

To obtain the dual statement, take the Fourier transform of each
column of $A_{\Z_n,g}$.
 By linearity, the resulting matrix can have no size-$n$ zero minors either, as that would
 mean that one column of the corresponding submatrix is a linear combination of other columns.
As $\widehat{M_\xi T_x g}=T_\xi M_{-x} \widehat g$, the resulting
matrix will correspond
 to $A_{\Z_n, \widehat g}$, except that modulations and translations have exchanged their roles.
So modulation adjacency becomes translation adjacency, which implies
the dual statement.





\noindent {\bf Acknowledgment.} We would like to thank Norbert
Kaiblinger, Franz Luef, and Ewa Matusiak for sharing with us the
results of their thorough discussions on the uncertainty principle
for functions on finite Abelian groups. Further, we thank Michael
Stoll for offering advice on algebraic geometry issues that are
relevant to our work, and Dan Alistarh and Sergiu Ungureanu for
writing some of Matlab code used.


\bibliography{krapfara0510}
\bibliographystyle{alpha}


\end{document}